\documentclass[letterpaper,english,11pt,aps,letter,superscriptaddress,nofootinbib,floatfix,pre]{revtex4}
\usepackage[T1]{fontenc}
\usepackage[latin1]{inputenc}
\usepackage{color}
\usepackage{verbatim}
\usepackage{amsmath}
\usepackage{amssymb}
\usepackage{graphicx}

\makeatletter

\pdfpageheight\paperheight
\pdfpagewidth\paperwidth

\providecommand{\tabularnewline}{\\}

\@ifundefined{textcolor}{}
{%
 \definecolor{BLACK}{gray}{0}
 \definecolor{WHITE}{gray}{1}
 \definecolor{RED}{rgb}{1,0,0}
 \definecolor{GREEN}{rgb}{0,1,0}
 \definecolor{BLUE}{rgb}{0,0,1}
 \definecolor{CYAN}{cmyk}{1,0,0,0}
 \definecolor{MAGENTA}{cmyk}{0,1,0,0}
 \definecolor{YELLOW}{cmyk}{0,0,1,0}
}

\usepackage{ae,aecompl}

\renewcommand{\citet}[1]{\cite{#1}}

\makeatother

\usepackage{babel}
\begin{document}

\title{An Immersed Boundary Method for Rigid Bodies}

\author{Bakytzhan Kallemov}

\affiliation{Courant Institute of Mathematical Sciences, New York University,
New York, NY 10012}

\author{Amneet Pal Singh Bhalla}

\affiliation{Department of Mathematics, University of North Carolina, Chapel Hill,
NC 27599}

\author{Boyce E. Griffith}

\affiliation{Departments of Mathematics and Biomedical Engineering, University
of North Carolina, Chapel Hill, NC 27599}

\author{Aleksandar Donev}

\email{donev@courant.nyu.edu}

\selectlanguage{english}%

\affiliation{Courant Institute of Mathematical Sciences, New York University,
New York, NY 10012}
\begin{abstract}
We develop an immersed boundary (IB) method for modeling flows around
fixed or moving rigid bodies that is suitable for a broad range of
Reynolds numbers, including steady Stokes flow. The spatio-temporal
discretization of the fluid equations is based on a standard staggered-grid
approach. Fluid-body interaction is handled using Peskin's IB method;
however, unlike existing IB approaches to such problems, we do not
rely on penalty or fractional-step formulations. Instead, we use an
unsplit scheme that ensures the no-slip constraint is enforced exactly
in terms of the Lagrangian velocity field evaluated at the IB markers.
Fractional-step approaches, by contrast, can impose such constraints
only approximately, which can lead to penetration of the flow into
the body, and are inconsistent for steady Stokes flow. Imposing no-slip
constraints exactly requires the solution of a large linear system
that includes the fluid velocity and pressure as well as Lagrange
multiplier forces that impose the motion of the body. The principal
contribution of this paper is that it develops an efficient preconditioner
for this exactly constrained IB formulation which is based on an analytical
approximation to the Schur complement. This approach is enabled by
the near translational and rotational invariance of Peskin's IB method.
We demonstrate that only a few cycles of a geometric multigrid method
for the fluid equations are required in each application of the preconditioner,
and we demonstrate robust convergence of the overall Krylov solver
despite the approximations made in the preconditioner. We empirically
observe that to control the condition number of the coupled linear
system while also keeping the rigid structure impermeable to fluid,
we need to place the immersed boundary markers at a distance of about
two grid spacings, which is significantly larger from what has been
recommended in the literature for elastic bodies. We demonstrate the
advantage of our monolithic solver over split solvers by computing
the steady state flow through a two-dimensional nozzle at several
Reynolds numbers. We apply the method to a number of benchmark problems
at zero and finite Reynolds numbers, and we demonstrate first-order
convergence of the method to several analytical solutions and benchmark
computations.
\end{abstract}
\maketitle
\global\long\def\V#1{\boldsymbol{#1}}
\global\long\def\M#1{\boldsymbol{#1}}
\global\long\def\Set#1{\mathbb{#1}}

\global\long\def\D#1{\Delta#1}
\global\long\def\d#1{\delta#1}

\global\long\def\norm#1{\left\Vert #1\right\Vert }
\global\long\def\abs#1{\left|#1\right|}

\global\long\def\av#1{\langle#1\rangle}

\global\long\def\grad{\boldsymbol{\nabla}}
\global\long\def\div{\boldsymbol{\nabla\cdot}}
\global\long\def\curl{\grad\times}
\global\long\def\half{\frac{1}{2}}

\global\long\def\Re{\text{Re}}

\global\long\def\p{\pi}

\global\long\def\HV{h^{d}}

\global\long\def\BV{V_{m}}

\global\long\def\DV{\text{\text{vol}}}

\global\long\def\r{\V r}

\global\long\def\q{\V R}

\global\long\def\Q{\V R}

\global\long\def\X{\V X}

\global\long\def\v{\V v}

\global\long\def\u{\V V}

\global\long\def\W{\V W}

\global\long\def\U{\V V}

\global\long\def\Bv{\V U}

\global\long\def\vort{\V{\omega}}

\global\long\def\Bomega{\V{\Omega}}

\global\long\def\Y{\V Y}

\global\long\def\F{\V F}

\global\long\def\f{\V f}

\global\long\def\cforce{\V{\lambda}}

\global\long\def\lamb{\V{\Lambda}}

\global\long\def\Lamb{\V{\Lambda}}

\global\long\def\torq{\V T}

\global\long\def\ctorq{\V{\mathcal{T}}}

\global\long\def\R{\V H}

\global\long\def\sM#1{\M{\mathcal{#1}}}

\global\long\def\Id{\sM I}

\global\long\def\J{\sM J}

\global\long\def\S{\sM S}

\global\long\def\A{\sM A}

\global\long\def\Div{\sM D}

\global\long\def\Grad{\sM G}

\global\long\def\Lap{\sM L_{v}}

\global\long\def\L{\sM L}

\global\long\def\Lp{\sM L_{p}}

\global\long\def\Mob{\sM M}

\global\long\def\BMob{\sM N}

\global\long\def\K{\sM K}

\textbf{}

\section{Introduction}

A large number of numerical methods have been developed to simulate
interactions between fluid flows and immersed bodies. For rigid bodies
or bodies with prescribed kinematics, many of these approaches \citet{DirectForcing_Uhlmann,RigidIBAMR,IBM_Projection,RigidIBAMR_ZhangGuy,IBM_Rigid_Boundary}
are based on the immersed boundary (IB) method of Peskin \citet{IBM_PeskinReview}.
The simplicity, flexibility, and power of the IB method for handling
a broad range of fluid-structure interaction problems was demonstrated
by Bhalla \emph{et al.} \citet{RigidIBAMR}. In that study, the authors
showed that the IB method can be used to model complex flows around
rigid bodies moving with specified kinematics (e.g., swimming fish
or beating flagella) as well as to compute the motion of freely moving
bodies driven by flow. In the approach of Bhalla \emph{et al.}, as
well as those of others \citet{DirectForcing_Uhlmann,IBM_Projection,IBM_Rigid_Boundary,RigidIBAMR_ZhangGuy},
the rigidity constraint enforcing that the fluid follows the motions
of the rigid bodies is imposed only \emph{approximately.} Here and
throughout this manuscript, when we refer to the \emph{no slip condition},
we mean the requirement that the interpolated fluid velocity exactly
match the rigid body velocity at the positions of the IB marker points\emph{.}
In this work, we develop an effective solution approach to an IB formulation
of this problem that \emph{exactly} enforces both the incompressibility
and no-slip constraints, thus substantially improving upon a large
number of existing techniques.

A simple approach to implementing rigid bodies using the traditional
IB method is to use stiff springs to attach markers that discretize
the body to tether points constrained to move as a rigid body \citet{TetherPoint_IBM}.
This penalty-spring approach leads to numerical stiffness and, when
the forces are handled explicitly, requires very small time steps.
For this reason, a number of \textit{direct forcing} IB methods \citet{IBSE_Poisson}
have been developed that aim to constrain the flow inside the rigid
body by treating the fluid-body force as a Lagrange multiplier $\Lamb$
enforcing a no-slip constraint at the locations of the IB markers.
However, to our knowledge, all existing direct forcing IB methods
use some form of time step splitting to separate the coupled fluid-body
problem into more manageable pieces. The basic idea behind these approaches
is first to solve a simpler system in which a number of the constraints
(e.g., incompressibility, or no-slip along the fluid-body interface)
are ignored. The solution of the unconstrained problem is then \textit{projected}
onto the constraints, which yields \emph{estimates} of the true Lagrange
multipliers. In most existing methods, the fluid solver uses a fractional
time stepping scheme, such as a version of Chorin's projection method,
to separate the velocity update from the pressure update \citet{DirectForcing_Uhlmann,IBM_Projection,IBM_Rigid_Boundary}.
Taira and Colonius also use a fractional time-stepping approach in
which they split the velocity from the Lagrange multipliers $\left(\p,\Lamb\right)$.
They obtain approximations to $\left(\p,\Lamb\right)$ in a manner
similar to that in a standard projection method for the incompressible
Navier-Stokes equations. A modified Poisson-type problem (see (26)
in \citet{IBM_Projection}) determines the Lagrange multipliers and
is solved using an unpreconditioned conjugate gradient method. The
method developed in Ref. \citet{RigidIBAMR} avoids the pressure-velocity
splitting and instead uses a combined iterative Stokes solver, and
in Ref. \citet{RigidIBAMR_ZhangGuy} (see supplementary material),
periodic boundary conditions are applied, which allows for the use
of a pseudo-spectral method. In both works, however, time step splitting
is still used to separate the computation of the rigidity constraint
forces from the updates to the fluid variables. In the approach described
in the supplementary material to Ref. \citet{RigidIBAMR_ZhangGuy},
the projection step of the solution onto the rigidity constraint is
performed twice in a predictor-corrector framework, which improves
the imposition of the constraint; however, this approach does not
control the accuracy of the approximation of the constraint forces.
Curet et al. \citet{FIISPA_Patankar} and Ardekani et al. \citet{FIISPA_Ardekani}
go a step closer in the direction of exactly enforcing the rigidity
constraint by iterating the correction until the relative slip between
the desired and imposed kinematics inside the rigid body reaches a
relatively loose tolerance of 1\%. The scheme used in Ref. \citet{FIISPA_Patankar}
is essentially a fixed-point (Richardson) iteration for the constrained
fluid problem, which uses splitting to separate the update of the
Lagrange multipliers from a fluid update based on the SIMPLER scheme
\citet{CFD_Patankar}. Unlike the approach developed here, fixed point
iterations based on splitting are not guaranteed to converge, yet
alone converge rapidly, especially in the steady Stokes regime for
tight solver tolerances.

An alternative view of direct forcing methods that use time step splitting
is that they are penalty methods for the unsplit problem, in which
the penalty parameter is related to the time step size. Such approaches
inherently rely on inertia and implicitly assume that fluid velocity
has memory. Consequently, \emph{all} such splitting methods fail in
the steady Stokes limit. Furthermore, even at finite Reynolds numbers,
methods based on splitting cannot exactly satisfy the no-slip constraint
at fluid-body interfaces. Such methods can thereby produce undesirable
artifacts in the solution, such as penetration of the flow through
a rigid obstacle. It is therefore desirable to develop a numerical
method that solves for velocity, pressure, and fluid-body forces in
a single step with controlled accuracy and reasonable computational
complexity.

The goal of this work is to develop an effective IB method for rigid
bodies that does not rely on any splitting. Our method is thus applicable
over a broad range of Reynolds numbers, including steady Stokes flow,
and is able to impose rigidity constraints exactly. This approach
requires us to solve large linear systems for velocity, pressure,
and fluid-body interaction forces. This linear system is not new.
For example, (13) in Ref. \citet{IBM_Projection} is essentially the
same system of equations that we study here. The primary contributions
of this work are that we do not rely on any approximations when solving
this linear system, and that we develop an effective preconditioner
based on an approximation of the Schur complement that allows us to
solve (\ref{eq:constrained_Stokes}). The resulting method has a computational
complexity that is only a few times larger than the corresponding
problem in the absence of rigid bodies. In the context of steady Stokes
flows, a rigid-body IB formulation very similar to the one we use
here has been developed by Bringley and Peskin \citet{IBM_Sphere};
however, that formulation relies on periodic boundary conditions,
and uses a very different spatial discretization and solution methodology
from the approach we describe here. Our approach can readily handle
a broad range of specified boundary conditions. In both Refs. \citet{IBM_Sphere}
and a very recent work by Stein \emph{et al}. on a higher-order IB
smooth extension method for scalar (e.g., Poisson) equations \citet{IBSE_Poisson},
the Schur complement is formed densely in an expensive pre-computation
stage. By contrast, in the method proposed here we build a simple
\emph{physics-based} approximation of the Schur complement that can
be computed ``on the fly'' in a scalable and efficient manner.

Our basic solution approach is to use a preconditioned Krylov solver
for the fully constrained fluid problem, as has been done for some
time in the context of finite element methods for fluid flows interacting
with elastic bodies \citet{heil2008solvers,FluidStructure_FEM_AMG}.
A key difficulty that we address in this work is the development of
an efficient preconditioner for the constrained formulation. To do
so, we construct an \emph{analytical }approximation of the Schur complement
(i.e., the mobility matrix) corresponding to Lagrangian rigidity forces
(i.e., Lagrange multipliers) enforcing the no-slip condition at the
positions of the IB markers. We rely on the near translational and
rotational invariance of Peskin's IB method to approximate the Schur
complement, following techniques commonly used for suspensions of
rigid spheres in steady Stokes flow such as Stokesian dynamics \citet{BrownianDynamics_OrderNlogN,StokesianDynamics_Rigid},
bead methods for rigid macromolecules \citet{HYDROLIB,SphereConglomerate,HYDROPRO,HYDROPRO_Globular}
and the method of regularized Stokeslets \citet{RegularizedStokeslets,RegularizedStokeslets_2D,RegularizedBrinkmanlet}.
In fact, as we explain herein, many of the techniques developed in
the context of steady Stokes flow can be used with the IB method both
at zero and also, perhaps more surprisingly, finite Reynolds numbers.

The method we develop offers an attractive alternative to existing
techniques in the context of steady or nearly-steady Stokes flow of
suspensions of rigid particles. To our knowledge, most other approaches
tailored to the steady Stokes limit rely on Green's functions for
Stokes flow to eliminate the (Eulerian) fluid degrees of freedom and
solve only for the (Lagrangian) degrees of freedom associated to the
surface of the body. Because these approaches rely on the availability
of analytical solutions, handling non-trivial boundary conditions
(e.g., bounded systems) is complicated \citet{BD_IBM_Graham} and
has to be done on a case-by-case basis \citet{BrownianDynamics_DNA2,StokesianDynamics_Wall,StokesianDynamics_Slit,StokesianDynamics_Confined,BD_LB_Comparison,RegularizedStokeslets_Walls,RegularizedStokeslets_Periodic}.
By contrast, in the method developed here, analytical Green's functions
are replaced by an ``on the fly'' computation that may be carried
out by a standard finite-volume, finite-difference, or finite-element
fluid solver %
\footnote{In this work, we use a staggered-grid discretization on a uniform
grid combined with multigrid-preconditioned Stokes solvers \citet{NonProjection_Griffith,StokesKrylov}.%
}. Such solvers can readily handle nontrivial boundary conditions.
Furthermore, suspensions at small but nonzero Reynolds numbers can
be handled without any extra work. Additionally, we avoid uncontrolled
approximations relying on truncations of multipole expansions to a
fixed order \citet{BrownianDynamics_OrderNlogN,ForceCoupling_Stokes,ISIBM,IrreducibleActiveFlows_PRL},
and we can seamlessly handle arbitrary body shapes and deformation
kinematics. For problems involving active \citet{ActiveSuspensions}
particles, it is straightforward to add osmo- or electro-phoretic
coupling between the fluid flow and additional fluid variables such
as the electric potential or the concentration of charged ions or
chemical reactants. Lastly, in the spirit of fluctuating hydrodynamics
\citet{BrownianBlobs,ForceCoupling_Fluctuations,SELM}, it is straightforward
to generate the stochastic increments required to simulate the Brownian
motion of small rigid particles suspended in a fluid by including
a fluctuating stress in the fluid equations. We also point out that
our method also has some disadvantages compared to methods such as
boundary integral or boundary element methods. Notably, it requires
filling the domain with a dense uniform fluid grid, which is expensive
at low densities. It is also a low-order method that cannot compute
solutions as accurately as spectral boundary integral formulations.
We do believe, nevertheless, that the method developed here offers
a good compromise between accuracy, efficiency, scalabilty, flexibility
and extensibility, compared to other more specialized formulations.

\section{Semi-Continuum Formulation}

Our notation uses the following conventions where possible. Vectors
(including multi-vectors), matrices, and operators are bolded, but
when fully indexed down to a scalar quantity we no longer bold the
symbol; matrices and operators are also scripted. We denote Eulerian
quantities with lowercase letters, and the corresponding Lagrangian
quantity with the same capital letter. We use the Latin indexes $i,j,k,l,m$
to denote a specific fluid grid point or IB marker (i.e., physical
location with which degrees of freedom are associated), the indices
$p,q,r,s,t$ to denote a specific body in the multibody context, and
Greek superscripts $\alpha,\beta,\gamma$ to denote specific Cartesian
components. For example, $\v$ denotes fluid velocity (either continuum
or discrete), with $v_{k}^{\alpha}$ being the fluid velocity in direction
$\alpha$ associated with the face center $k$, and $\U$ denotes
the velocity of all IB markers, with $V_{i}^{\alpha}$ being the velocity
of marker $i$ along direction $\alpha$. Our formulation is easily
extended to a collection of rigid bodies, but for simplicity of presentation,
we focus on the case of a single body.

We consider a region $\mathcal{D}\subset\Set R^{d}$ ($d=2$ or $3$)
that contains a single rigid body $\Omega\subset\mathcal{D}$ immersed
in a fluid of density $\rho$ and shear viscosity $\eta$. The computational
domain $\mathcal{D}$ could be a periodic region (topological torus),
a finite box, an infinite domain, or some combination thereof, and
we will implicitly assume that some consistent set of boundary conditions
are prescribed on its boundary $\partial\mathcal{D}$ even though
we will not explicitly write this in the formulation. We require that
the linear velocity of a given reference point (e.g., the center of
mass of the body) $\Bv(t)$ and the angular velocity $\Bomega(t)$
of the body are specified functions of time, and without loss of generality,
we assume that the rigid body is at rest %
\footnote{The case of more general specified kinematics is a straightforward
generalization and does not incur any additional mathematical or algorithmic
complexity \citet{IBAMR_Fish,RigidIBAMR}.%
}. In addition to features of the fluid flow, typical quantities of
interest are the total drag force $\F(t)$ and total drag torque $\torq(t)$
between the fluid and the body. Another closely related problem to
which IB methods can be extended is the case when the motion of the
rigid body (i.e., $\Bv(t)$ and $\Bomega(t)$) is not known but the
body is subject to specified external force $\F(t)$ and torque $\torq(t)$.
For example, in the sedimentation of rigid particles in suspension,
the external force is gravity and the external torque is zero. Handling
this \emph{free kinematics} problem \citet{IBAMR_Fish,RigidIBAMR}
requires a nontrivial extension of our formulation and numerical algorithm.

\begin{figure}[h]
\begin{centering}
\includegraphics[width=0.49\textwidth]{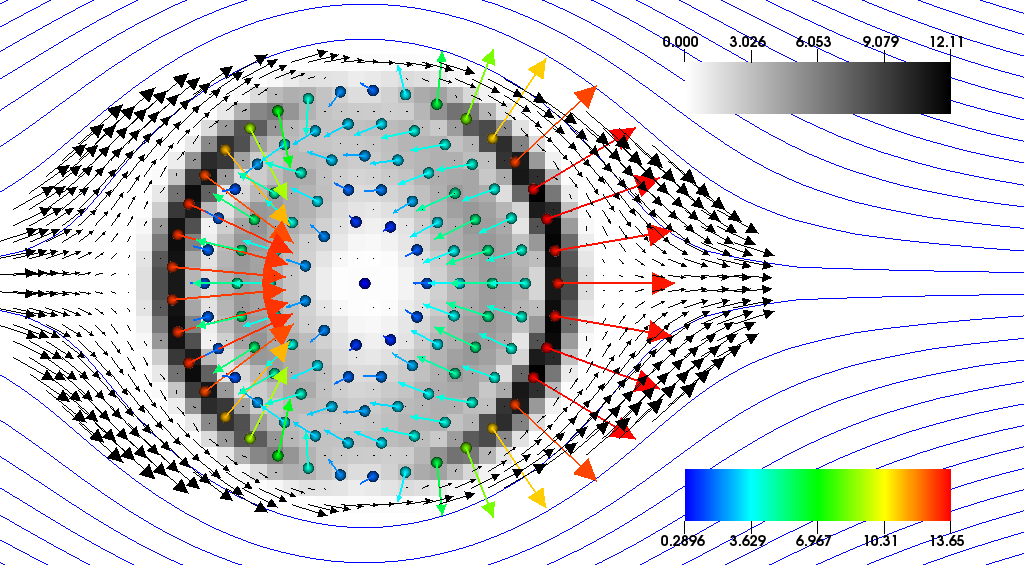}\includegraphics[width=0.49\textwidth]{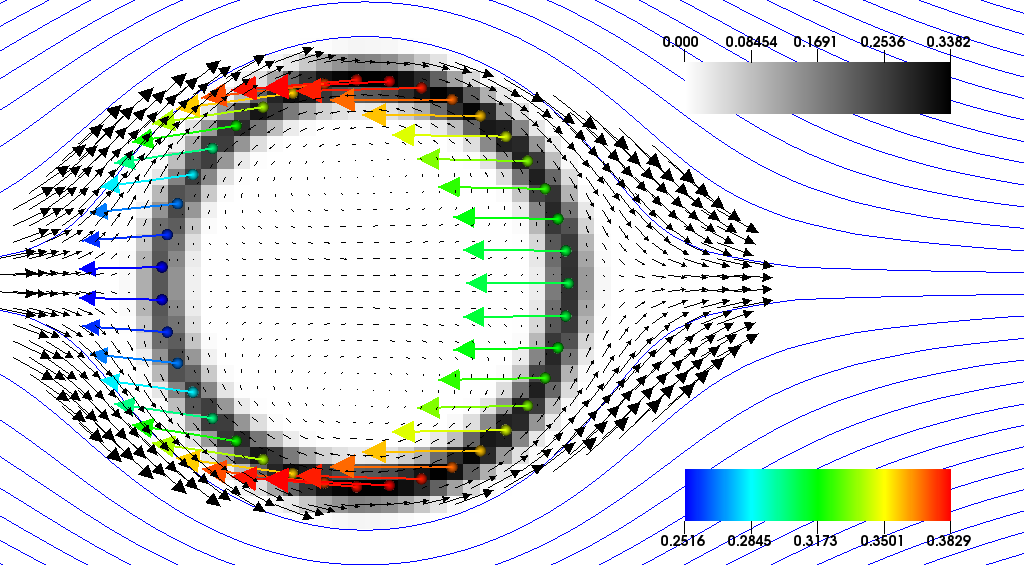}
\par\end{centering}

\centering{}\caption{\label{fig:2DWake_Stokes}Two-dimensional Steady Stokes flow past
a periodic column of circular cylinders (disks) at zero Reynolds number
obtained using our rigid-body IB method (the same setup is also studied
at finite Reynolds numbers in Section \ref{sub:UnsteadyCylinders}).
The markers used to mediate the fluid-body interaction are shown as
small colored circles. The Lagrangian constraint forces $\Lamb$ that
keep the markers at their fixed locations are shown as colored vectors;
the color of the vectors and the corresponding marker $i$ are based
on the magnitude of the constraint force $\lamb_{i}$ (see color bar).
The fluid velocity field is shown as a vector field (black arrows)
in the vicinity and the interior of the body; further from the body,
flow streamlines are shown as solid blue lines. The magnitude of the
Eulerian constraint force $\S\Lamb$ is shown as a gray color plot
(see greyscale bar). (Left panel) A volume marker grid of 121 markers
is used to discretize the disk. The majority of the constraint forces
are seen to act near the surface of the body, but nontrivial constraint
forces are seen also in the interior of the body. (Right panel) A
surface grid of 39 markers is used to discretize the disk, which strictly
localizes the constraint forces to the surface of the body.}
\end{figure}

In the immersed boundary (IB) method \citet{IBM_PeskinReview,IBAMR,IBAMR_HeartValve},
the velocity field $\v(\r,t)$ is extended over the whole domain $\mathcal{D}$,
\emph{including} the body interior. The body is discretized using
a collection of \emph{markers}, which is a set of $N$ points that
cover the interior of the body and at which the interaction between
the body and the fluid is localized. For example, the markers could
be the nodes of a triangular ($d=2$) or tetrahedral ($d=3$) mesh
used to discretize $\Omega$; an illustration of such a \emph{volume}
grid of markers discretizing a rigid disk immersed in steady Stokes
flow is shown in the left panel of Fig. \ref{fig:2DWake_Stokes}.
In the case of Stokes flow, the specification of a no slip condition
on the boundary of a rigid body is sufficient to ensure rigidity of
the fluid inside the body \citet{RegularizedStokeslets}. Therefore,
for Stokes flow, the grid of markers does not need to extend over
the volume of the body and can instead be limited to the surface of
the rigid body, thus substantially reducing the number of markers
required to represent the body. In this case, the markers could be
the nodes of a triangulation ($d=3$) of the surface of the body;
an illustration of such a \emph{surface} grid of markers is shown
in the right panel of Fig. \ref{fig:2DWake_Stokes}. We discuss the
differences between a volume and a surface grid of markers in Section
\ref{sec:Results}.

The traditional IB method is concerned with the motion of elastic
(flexible) bodies in fluid flow, and the collection of markers can
be viewed as a set of quadrature points used to discretize integrals
over the moving body. The elastic body forces are most easily computed
in a Lagrangian coordinate system attached to the deforming body,
and the relative positions of the markers in the fixed Eulerian frame
of reference generally vary in time. For a rigid body, however, the
relative positions of the markers do not change, and it is not necessary
to introduce two distinct coordinate frames. Instead, we use the same
Cartesian coordinate system to describe points in the fluid domain
and in the body; the positions of the $N$ markers in this fixed frame
of reference will be denoted with $\Q=\left\{ \q_{1},\dots,\q_{N}\right\} $,
where $\Q\subset\Omega$ for volume meshes or $\Q\subset\partial\Omega$
for surface meshes.

In the standard IB method for flexible immersed bodies, elastic forces
are computed in the Lagrangian frame and then \emph{spread} to the
fluid in the neighborhood of the markers using a regularized delta
function $\delta_{a}\left(\r\right)$ that integrates to unity and
converges to a Dirac delta function as the regularization width $a\rightarrow0$.
The regularization length scale $a$ is typically chosen to be on
the order of the spacing between the markers (as well as the lattice
spacing of the grid used to discretize the fluid equations), as we
discuss in more detail later. In turn, the motion of the markers is
specified to follow the velocity of the fluid \emph{interpolated}
at the positions of the markers.

The key difference between an elastic and a rigid body is that, for
a rigid object, the motion of the markers is known (e.g, they are
fixed in place or move with a specified velocity) and the body forces
are unknown and must be determined within each time step. To obtain
the fluid-marker interaction forces $\Lamb\left(t\right)=\left\{ \lamb_{i}\left(t\right),\dots,\lamb_{N}\left(t\right)\right\} $
that constrain the motion of the $N$ markers, we solve for the Eulerian
velocity field $\v\left(\r,t\right)$, the Eulerian pressure field
$\p\left(\r,t\right)$, and the Lagrangian constraint forces $\Lamb_{i}\left(t\right)$
the system 
\begin{align}
\rho\left(\partial_{t}\v+\v\cdot\grad\v\right)+\grad\p & =\eta\grad^{2}\v+\sum_{i=1}^{N}\lamb_{i}\delta_{a}\left(\q_{i}-\r\right),\nonumber \\
\div\v & =0,\nonumber \\
\u_{i}=\int\delta_{a}\left(\q_{i}-\r\right)\v\left(\r,t\right)\, d\r & =0,\quad i=1,\dots,N,\label{eq:semi_continuum}
\end{align}
along with suitable boundary conditions. In the case of steady Stokes
flow, we set $\rho=0$. The first two equations are the incompressible
Navier-Stokes equations with an Eulerian constraint force 
\[
\cforce\left(\r,t\right)=\sum_{i=1}^{N}\lamb_{i}\delta_{a}\left(\q_{i}-\r\right).
\]
The last condition is the \emph{rigidity constraint} that requires
that the Eulerian velocity averaged around the position of marker
$i$ must match the known marker velocity $\u_{i}$. This constraint
enforces a regularized no-slip condition at the locations of the IB
markers, which is a numerical approximation of the true no-slip condition
on the surface (or interior) of the body. Observe that flow may still
penetrate the body in-between the markers and this leads to a well-known
small but nonzero ``leak'' in the traditional Peskin IB method.
This leak can be greatly reduced by adopting a staggered-grid formulation
\citet{IBM_Staggered}, as done in the present work. Other more specialized
approaches to reducing spurious fluxes in the IB method have been
developed \citet{ImprovedLeak,IBFE,DivFreeIB}, but will not be considered
in this work.

Notice that for zero Reynolds number, the semi-continuum formulation
\eqref{eq:semi_continuum} is closely related to the popular method
of regularized Stokeslets, which solves a similar system of equations
for $\Lamb$\textcolor{green}{{} }\citet{RegularizedStokeslets,RegularizedStokeslets_2D}.
The key difference %
\footnote{Another important difference is that we follow Peskin and use the
regularized delta function both for spreading and interpolation (this
ensures energy conservation in the formulation \citet{IBM_PeskinReview}),
whereas in the method of Regularized Stokeslets only the spreading
uses a regularized delta function. Our choice ensures that the linear
system we solve is symmetric and positive semi-definite, which is
crucial if one wishes to account for Brownian motion and thermal fluctuations.%
} is that in the method of regularized Stokeslets, the fluid equations
are eliminated using analytic Green's functions; this necessitates
that nontrivial pre-computations of these Green's functions be performed
for each type of boundary condition \citet{RegularizedStokeslets_Walls,RegularizedStokeslets_Periodic}. 

In this work, we treat \eqref{eq:semi_continuum} as the primary continuum
formulation of the problem. This is a semi-continuum formulation in
which the rigid body is represented as a discrete collection of markers
but the fluid description is kept as a continuum, which implies that
different discretizations of the fluid equations are possible. One
can, in principle, try to write a fully continuum formulation in which
the discrete set of rigidity forces $\Lamb$ are replaced by a continuum
force density field $\cforce\left(\q\in\Omega,t\right)$. The well-posedness
and stability of such a fully continuum formulation is mathematically
delicate, however, and there can be subtle differences between weak
and strong interpretations of the equations. To appreciate this, observe
that if each component of the velocity is discretized with $N_{f}$
degrees of freedom, it cannot in general be possible to constrain
the velocity strongly at more than $N_{f}$ points (markers). By contrast,
in our strong formulation \eqref{eq:semi_continuum}, the velocity
is infinite dimensional but it is only constrained in the vicinity
of a finite number of markers. Therefore, the problem \eqref{eq:semi_continuum}
is always well posed and is directly amenable to numerical discretization
and solution, at least when it is well-conditioned. As we show in
this work, the conditioning of the fully discrete problem is controlled
by the relationship between the regularization length $a$ and the
marker spacing. 

The physical interpretation of the constraint forces $\lamb_{i}$
depends on details of the marker grid and the type of the problem
under consideration. For fully continuum formulations, in which the
fluid-body interaction is represented solely as a surface force density,
the force $\lamb_{i}$ can be interpreted as the integral of the traction
(normal component of the fluid stress tensor) over a surface area
associated with marker $i$. Such a formulation is appropriate, for
example, for steady Stokes flow. In particular, for steady Stokes
flow our method can be seen as a discretized and regularized first-kind
integral formulation in which Green's functions are computed by the
fluid solver. This approach is different from the method of regularized
Stokeslets, in which regularized Green's functions must be computed
analytically \citet{RegularizedStokeslets_2D,RegularizedStokeslets}. 

For cases in which markers are placed on both the surface and the
interior of a rigid body, the precise physical interpretation of the
volume force density, and thus of $\Lamb$, is delicate even for steady
Stokes flow. Notably, observe that the splitting between a volume
constraint force density and the gradient of the pressure is not unique
because the pressure inside a rigid body cannot be determined uniquely.
Specifically, only the component of the constraint force density projected
onto the space of divergence-free vector fields is uniquely determined.
In the presence of finite inertia and a density mismatch between the
fluid and the \emph{moving} rigid bodies, the inertial terms in \eqref{eq:semi_continuum}
need to be modified in the interior of the body \citet{ISIBM}. Furthermore,
sufficiently many markers in the interior of the body are required
to prevent spurious angular momentum being generated by motions of
the fluid inside the body \citet{DirectForcing_Uhlmann}. We do not
discuss these physical issues in this work because they do not affect
the numerical algorithm, and because we restrict our numerical studies
to flow past \emph{stationary} rigid bodies, for which the fluid-body
interaction force is localized to the surface of the body in the continuum
limit.

\section{Discrete Formulation}

The spatial discretization of the fluid equation uses a uniform Cartesian
grid with grid spacing $h$ and is based on a second-order accurate
staggered-grid finite-difference discretization, in which vector-valued
quantities, including velocities and forces, are represented on the
faces of the Cartesian grid cells, and scalar-valued quantities, including
the pressure, are represented at the centers of the grid cells \citet{IBAMR,IBAMR_HeartValve,ISIBM,RigidIBAMR}.
Our implicit-explicit temporal discretization of the Navier-Stokes
equation is standard and summarized in prior work; see for example
the work of Griffith \citet{IBAMR_HeartValve}. The key features are
that we treat advection explicitly using a predictor-corrector approach,
and that we treat viscosity implicitly, using either the backward
Euler or the implicit midpoint method. For steady Stokes flow, no
temporal discretization required, although one can also think of this
case as corresponding to a backward Euler discretization of the time-dependent
problem with a very large time step size $\D t$. A key dimensionless
quantity is the viscous CFL number $\beta=\nu\D t/h^{2}$, where the
kinematic viscosity is $\nu=\eta/\rho$. If $\beta$ is small, the
pressure and velocity are weakly coupled, but for large $\beta$,
and in particular for the steady Stokes limit $\beta\rightarrow\infty$,
the coupling between the velocity and pressure equations is strong. 

We do \emph{not} use a fractional time-stepping scheme (i.e., a projection
method) to split the pressure and velocity updates; instead, the pressure
is treated as a Lagrange multiplier that enforces the incompressibility
and must be determined together with the velocity at the end of the
time step \citet{NonProjection_Griffith}; except in special cases,
this is \emph{necessary} for small Reynolds number flows. This approach
also greatly aids with imposing stress boundary conditions \citet{NonProjection_Griffith}.
The constraint force $\cforce\left(\r,t\right)$ is treated analogously
to the pressure, i.e., as a Lagrange multiplier. Whereas the role
of the pressure is to enforce the incompressibility constraint, $\cforce$
enforces the rigidity constraint. Like the pressure, $\cforce$ is
an unknown that must be solved for in this formulation.

\subsection{Force spreading and velocity interpolation}

In the fully discrete formulation of the fluid-body coupling, we replace
spatial integrals by sums over fluid or body grid points in the semi-continuum
formulation \eqref{eq:semi_continuum}. The regularized delta function
is discretized using a tensor product in $d$-dimensional space (see
\citet{ISIBM} for more details), 
\[
\delta_{a}(\r)=h^{-d}\prod_{\alpha=1}^{d}\phi_{a}\left(r_{\alpha}\right),
\]
where $\HV$ is the volume of a grid cell. The one-dimensional kernel
function $\phi_{a}$ is chosen based on numerical considerations of
efficiency and maximized approximate translational invariance \citet{IBM_PeskinReview}.
In this work, for reasons that will become clear in Section \ref{sec:ApproxMob},
we prefer to use a kernel that maximizes translational and rotational
invariance (i.e., improves grid-invariance). We therefore use the
smooth (three-times differentiable) six-point kernel recently described
by Bao \emph{et al.} \citet{New6ptKernel}. This kernel is more expensive
than the traditional four-point kernel \citet{IBM_PeskinReview} because
it increases the support of the kernel to $6^{2}=36$ grid points
in two dimensions and $6^{3}=216$ grid points in three dimensions;
however, this cost is justified because the new six-point kernel improves
the translational invariance by orders of magnitude compared to other
standard IB kernel functions \citet{New6ptKernel}.

The interaction between the fluid and the rigid body is mediated through
two crucial operations. The discrete velocity-interpolation operator
$\J$ averages velocities on the staggered grid in the neighborhood
of marker $i$ via 
\[
\left(\J\v\right)_{i}^{\alpha}=\sum_{k}v_{k}^{\alpha}\;\phi_{a}\left(\q_{i}-\r_{k}^{\alpha}\right),
\]
where the sum is taken over faces $k$ of the grid, $\alpha$ indexes
coordinate directions ($x,y,z$) as a superscript, and $\r_{k}^{\alpha}$
is the position of the center of the grid face $k$ in the direction
$\alpha$. The discrete force-spreading operator $\S$ spreads forces
from the markers to the faces of the staggered grid via
\begin{equation}
\left(\S\Lamb\right)_{k}^{\alpha}=h^{-d}\sum_{i}\Lambda_{i}^{\alpha}\;\phi_{a}\left(\q_{i}-\r_{k}^{\alpha}\right),\label{eq:S_unweighted}
\end{equation}
where now the sum is over the markers that define the configuration
of the rigid body. These operators are adjoint with respect to a suitably-defined
inner product, $\J=\S^{\star}=\HV\,\S^{T}$, which ensures conservation
of energy \citet{IBM_PeskinReview}. Extensions of the basic interpolation
and spreading operators to account for the presence of physical boundary
conditions are described in Appendix \ref{AppendixPhysicalBCs}.

\subsection{Rigidly-constrained Stokes problem}

At every stage of the temporal integrator, we need to solve a linear
system of the form 
\begin{equation}
\left[\begin{array}{ccc}
\A & \Grad & -\S\\
-\Div & \M 0 & \M 0\\
-\J & \M 0 & \V 0
\end{array}\right]\left[\begin{array}{c}
\v\\
\pi\\
\Lamb
\end{array}\right]=\left[\begin{array}{c}
\V g\\
\V h=\V 0\\
\W=\V 0
\end{array}\right],\label{eq:constrained_Stokes}
\end{equation}
which is the focus of this work. The right-hand side $\V g$ includes
all remaining fluid forcing terms, explicit contributions from previous
time steps or stages, boundary conditions, etc. Here, $\Grad$ is
the discrete gradient operator, $\Div=-\Grad^{T}$ is the discrete
divergence operator, and $\A$ is the vector equivalent of the familiar
screened Poisson (or Helmholtz) operator
\[
\A=\frac{\rho}{\D t}\Id-\frac{\kappa\eta}{h^{2}}\Lap,
\]
with $\kappa=1$ for the backward Euler method or for steady Stokes,
and $\kappa=1/2$ for the implicit midpoint rule. Here $\Lap$ is
the dimensionless vector Laplacian operator, which takes into account
boundary conditions for velocity such as no-slip boundaries. Since
the viscosity appears multiplied by the coefficient $\kappa$, we
will henceforth absorb this coefficient into the viscosity, $\eta\leftarrow\kappa\eta$,
which allows us to assume, without loss of generality, that $\mbox{\ensuremath{\kappa}=1}$
and to write the fluid operator in the form
\begin{equation}
\A=\eta h^{-2}\left(\beta^{-1}\Id-\Lap\right).\label{eq:A_def}
\end{equation}

We remark that making the $(3,3)$ block in the matrix in (\ref{eq:constrained_Stokes})
non-zero (i.e., regularizing the saddle-point system) is closely related
to solving the Brinkman equations \citet{Brinman_Original} for flow
through a permeable or porous body suspended in fluid \citet{RigidIBAMR_ZhangGuy}.
In particular, by making the $(3,3)$ block a diagonal matrix with
suitable diagonal elements, one can consistently discretize the Brinkman
equations. Such regularization greatly simplifies the numerical linear
algebra except, of course, when the permeability of the body is so
small that it effectively acts as an impermeable body. In this work,
we focus on developing a solver for (\ref{eq:constrained_Stokes})
that is effective even when there is no regularization (permeability),
and even when the matrix $\A$ is the discretization of an elliptic
operator, as is the case in the steady Stokes regime. This is the
hardest case to consider, and a solver that is robust in this case
will be able to handle the easier cases of finite Reynolds number
or permeable bodies with ease.

It is worth noticing the structure of the linear system (\ref{eq:constrained_Stokes}).
First, observe that the system is symmetric, at least if only simple
boundary conditions such as periodic or no-slip boundaries are present
\citet{NonProjection_Griffith}. In the top $1\times1$ block, $\A\succsim\M 0$
is a symmetric positive-semidefinite (SPD) matrix. The top left $2\times2$
block represents the familiar saddle-point problem arising when solving
the Navier-Stokes or Stokes equations in the absence of a rigid body
\citet{NonProjection_Griffith}. The whole system is a saddle-point
problem for the fluid variables and for $\Lamb$, in which the top-left
block is the Stokes saddle-point matrix.

...

\subsection{Mobility matrix}

We can formally solve (\ref{eq:constrained_Stokes}) through a Schur
complement approach, as described in more detail in Section \ref{sec:Solver}.
For increased generality, which will be useful when discussing preconditioners,
we allow the right hand side to be general and, in particular, do
not assume that $\V h$ and $\V W$ are zero. 

First, we solve the unconstrained fluid equation for pressure and
velocity 
\begin{equation}
\left[\begin{array}{cc}
\A & \Grad\\
-\Div & \M 0
\end{array}\right]\left[\begin{array}{c}
\v\\
\p
\end{array}\right]=\left[\begin{array}{c}
\S\Lamb+\V g\\
\V h
\end{array}\right],\label{eq:Stokes_sub}
\end{equation}
where we recall that $\A=\eta h^{-2}\left(\beta^{-1}\Id-\Lap\right)$.
The solution can be written as $\v=\L^{-1}\left(\S\Lamb+\V g\right)+\L_{p}^{-1}\V h$,
where $\L^{-1}$ is the standard Stokes solution operator for divergence-free
flow ($\V h=0$), given by
\begin{equation}
\L^{-1}=\A^{-1}-\A^{-1}\Grad\left(\Div\A^{-1}\Grad\right)^{-1}\Div\A^{-1},\label{eq:L_inv_general}
\end{equation}
where we have assumed for now that $\A^{-1}$ is invertible. For a
periodic system, the discrete operators commute, and we can write
\begin{equation}
\L^{-1}=\sM P\A^{-1}=\left(\Id-\Grad\left(\Div\Grad\right)^{-1}\Div\right)\A^{-1},\label{eq:L_inv_periodic}
\end{equation}
where $\sM P$ is the Helmholtz projection onto the space of divergence-free
vector fields. We never explicitly compute or form $\L^{-1}$; rather,
we solve the Stokes velocity-pressure subsystems using the projection-method
based preconditioner developed by Griffith \citet{NonProjection_Griffith}.
Let us define $\tilde{\v}=\L^{-1}\V f+\L_{p}^{-1}\V h$ to be the
solution of the \emph{unconstrained} Stokes problem
\begin{equation}
\left[\begin{array}{cc}
\A & \Grad\\
-\Div & \M 0
\end{array}\right]\left[\begin{array}{c}
\tilde{\v}\\
\tilde{\p}
\end{array}\right]=\left[\begin{array}{c}
\V f\\
\V h
\end{array}\right],\label{eq:Stokes_saddle}
\end{equation}
giving $\v=\tilde{\v}+\L^{-1}\S\Lamb$.

Next, we plug the velocity $\v$ into the rigidity constraint, $\J\v=-\W$,
to obtain
\begin{equation}
\Mob\Lamb=-\left(\W+\J\tilde{\v}\right),\label{eq:lambda_sol}
\end{equation}
where the Schur complement or \emph{marker mobility matrix} is
\begin{equation}
\Mob=\J\L^{-1}\S=\S^{\star}\L^{-1}\S.\label{eq:mob_def}
\end{equation}
The mobility matrix $\Mob\succsim\M 0$ is SPD and has dimensions
$dN\times dN$, and the $d\times d$ block $\Mob_{ij}$ relates the
force applied at marker $j$ to the velocity induced at marker $i$.
Our approach to obtaining an efficient algorithm for the constraind
fluid-solid system is to develop a method for approximating the marker
mobility matrix $\Mob$ in a simple and efficient way that leads to
robust preconditioners for solving the mobility subproblem (\ref{eq:lambda_sol});
see Section \ref{sec:ApproxMob}.

Observe that the conditioning of the saddle-point system (\ref{eq:constrained_Stokes})
is controlled by the conditioning of $\Mob$. In particular, if the
(non-negative) eigenvalues of $\Mob$ are bounded away from zero,
then there will be a unique solution to the saddle-point system. If
this bound is uniform as the grid is refined, then the problem is
well-posed and will satisfy a stability criterion similar to the well-known
Ladyzenskaja-Babuska-Brezzi (LBB) condition for the Stokes saddle-point
problem (\ref{eq:Stokes_saddle}). We investigate the spectrum of
the the marker mobility matrix numerically in Section \ref{sec:Solver}.
In practice, there may be some nearly zero eigenvalues of the matrix
$\Mob$ corresponding to physical (rather than numerical) null modes.
An example is a sphere discretized with markers on the surface: we
know that a uniform compression of the sphere will not cause any effect
because of the incompressibility of the fluid filling the sphere.
This compression mode corresponds to a null-vector for the constraint
forces $\Lamb$; it poses no difficulties in principle because the
right-hand side in (\ref{eq:lambda_sol}) is always in the range of
$\Mob$. Of course, when a discrete set of markers is placed on the
sphere, the rotational symmetry will be broken and the corresponding
mode will have a small but nonzero eigenvalue, which can lead to numerical
difficulties if not handled with care.

\subsection{Periodic steady Stokes flow}

In the time-dependent context, $\beta$ is finite, and it is easy
to see that $\A\succ\M 0$ is invertible. The same happens even for
steady Stokes flow if at least one of the boundaries is a no-slip
boundary. In the case of periodic steady Stokes flow, however, $\A=-\eta h^{-2}\Lap$
has in its range vectors that sum to zero, because no nonzero total
force can be applied on a periodic domain. This means that a solvability
condition is
\[
\av{\S\Lamb+\V g}=\text{\ensuremath{\DV}}^{-1}\sum_{i=1}^{N}\lamb_{i}+\av{\V g}=\text{\ensuremath{\DV}}^{-1}\V 1^{T}\lamb+\av{\V g}=0,
\]
where $\av{}$ denotes an average over the whole system, $\V 1$ is
a vector of ones, and $\text{\ensuremath{\DV}}$ is the volume of
the domain. This is an additional constraint that must be added to
the constrained Stokes system (\ref{eq:constrained_Stokes}) for a
periodic domain the steady Stokes case. In this approach, the solution
has an indeterminate mean velocity $\av{\v}$ because momentum is
not conserved. This sort of approach is followed for a scalar (reaction-diffusion)
equivalent of (\ref{eq:constrained_Stokes}) in the Appendix of Ref.
\citet{ReactiveBlobs}, for the traditional Peskin IB method in Ref.
\citet{TetherPoint_IBM}, and for a higher-order IB method in \citet{IBSE_Poisson}.

Here, we instead impose the mean velocity $\av{\v}=\V 0$ and ensure
that the total force applied to the fluid sums to zero, i.e., we enforce
momentum conservation. Specifically, for the special case of periodic
steady Stokes, we solve the system
\begin{equation}
\left[\begin{array}{ccc}
\A & \Grad & -\left(\S-\text{\ensuremath{\DV}}^{-1}\V 1^{T}\right)\\
-\Div & \M 0 & \M 0\\
-\J & \M 0 & \M 0
\end{array}\right]\left[\begin{array}{c}
\v\\
\p\\
\Lamb
\end{array}\right]=\left[\begin{array}{c}
\V g\\
\V h\\
\W
\end{array}\right],\label{eq:constrained_periodic}
\end{equation}
together with the constraint $\av{\v}=\V 0$, where we assume that
$\av{\V g}=\V 0$ for consistency. This change amounts to simply redefining
the spreading operator to subtract the total applied force on the
markers as a uniform force density, $\S\leftarrow\S-\text{\ensuremath{\DV}}^{-1}\V 1^{T}$.
This can be justified by considering the unit cell to be part of an
infinite periodic system in which there is an externally applied constant
pressure gradient, which is balanced by the drag forces on the bodies
so as to ensure that the domain as a whole is in force balance \citet{FiniteRe_3D_Ladd,FiniteRe_2D_Ladd,VACF_Ladd}.

\section{\label{sec:ApproxMob}Approximating the Mobility Matrix}

A key element in the preconditioned Krylov solver for (\ref{eq:constrained_Stokes})
that we describe in Section \ref{sec:Solver} is an approximate solver
for the mobility subproblem (\ref{eq:lambda_sol}). The success of
this approximate solver, i.e., the accuracy with which we can approximate
the Schur complement of the saddle-point problem (\ref{eq:constrained_Stokes}),
is crucial to an effective linear solver and one of the key contributions
of this work.

Because it involves the inverse Stokes operator $\L^{-1}$, the actual
Schur complement $\Mob=\S^{\star}\L^{-1}\S$ cannot be formed efficiently.
Instead of forming the true mobility matrix, we instead approximate
$\Mob\approx\widetilde{\Mob}$ by a \emph{dense} but low-rank \emph{approximate
mobility matrix} $\widetilde{\Mob}$ given by simple \emph{analytical}
approximations. To achieve this, we use two key ideas:
\begin{enumerate}
\item We ignore the specifics of the boundary conditions and assume that
the structure is immersed in an infinite domain at rest at infinity
(in three dimensions) or in a finite periodic domain (in two dimensions).
This implies that the Krylov solver for (\ref{eq:constrained_Stokes})
must handle the boundary conditions.
\item We assume that the IB spatial discretization is translationally and
rotationally invariant; that is, $\Mob$ does not depend on the exact
position and orientation of the body relative to the underlying fluid
grid. This implies that the Krylov solver must handle any grid-dependence
in the solution.
\end{enumerate}
The first idea, to ignore the boundary conditions in the preconditioner,
has worked well in the context of solving the Stokes system (\ref{eq:Stokes_saddle}).
Namely, a simple but effective approximation of the inverse of the
Schur complement for (\ref{eq:Stokes_saddle}), $\left(\Div\A^{-1}\Grad\right)^{-1}$,
can be constructed by assuming that the domain is periodic so that
the finite difference operators commute, and thus the Schur complement
degenerates to a diagonal or nearly-diagonal mass matrix \citet{ApproximateCommutators,NonProjection_Griffith,StokesKrylov}.
The second idea, to make use of the near grid invariance of Peskin's
regularized kernel functions, has previously been used successfully
in implicit immersed-boundary methods by Ceniceros \emph{et al.} \citet{IBM_Implicit_Fisher}.
Note that for certain choices of the kernel function, the assumption
of grid invariance can be a very good approximation to reality; here,
we rely on the recently-developed six-point kernel \citet{New6ptKernel},
which has excellent grid invariance and relatively compact support.

In the remainder of this section, we explain how we compute the entries
in $\widetilde{\Mob}$ in three dimensions, assuming an unbounded
fluid at rest at infinity. The details for two dimensions are given
in Appendix \ref{AppendixMob2D} and are similar in nature, except
for complications for two-dimensional steady Stokes flow resulting
from the well-known Stokes paradox. 

The mobility matrix $\Mob$ is a symmetric block matrix built from
$N\times N$ blocks of size $d\times d$. The block $\Mob_{ij}$ corresponding
to markers $i$ and $j$ relates a force applied at marker $j$ to
the velocity induced at marker $i$. Our basic assumption is that
$\Mob_{ij}$ does not depend on the actual position of the markers
relative to the fluid grid, but rather only depends on the distance
between the two markers and on the viscous CFL number $\beta$ in
the form
\begin{equation}
\widetilde{\Mob}_{ij}=f_{\beta}\left(r_{ij}\right)\Id+g_{\beta}\left(r_{ij}\right)\hat{\r}_{ij}\otimes\hat{\r}_{ij},\label{eq:M_tilde_ij}
\end{equation}
where $\r_{ij}=\q_{i}-\q_{j}$ and $r_{ij}$ is the distance between
the two markers, and hat denotes a unit vector. The functions of distance
$f_{\beta}(r)$ and $g_{\beta}(r)$ depend on the specific kernel
chosen, the specific discretization of the fluid equations (in our
case the staggered-grid scheme), and the the viscous CFL number $\beta$.
To obtain a specific form for these two functions, we empirically
fit numerical data with functions with the proper asymptotic behavior
at short and large distances between the markers. For this purpose,
we first discuss the asymptotic properties of $f_{\beta}(r)$ and
$g_{\beta}(r)$ from a physical perspective.

It is important to note that the true mobility matrix $\Mob$ is guaranteed
to be SPD because of its structure and the adjointness of the spreading
and interpolation operators. This can be ensured for the approximation
$\widetilde{\Mob}$ by placing positivity constraints on suitable
linear combinations of the Fourier transforms of $f_{\beta}\left(r\right)$
and $g_{\beta}(r)$, which ensure that the kernel $\underline{\Mob}\left(\V r_{i},\M r_{j}\right)$
given by (\ref{eq:M_tilde_ij}) is SPD in the sense of integral operators.
It is, however, very difficult to place such constraints on empirical
fits in practice, and in this work, we do not attempt to ensure $\widetilde{\Mob}$
is SPD for all marker configurations.

\subsection{Physical Constraints}

Let us temporarily focus on the semi-continuum formulation \eqref{eq:semi_continuum}
and ignore Eulerian discretization artifacts. The pairwise mobility
between markers $i$ and $j$ for a continuum fluid is
\begin{equation}
\Mob_{ij}=\eta^{-1}\int\delta_{a}(\q_{i}-\V r^{\prime\prime})\Set G(\V r^{\prime\prime},\V r^{\prime})\delta_{a}(\q_{j}-\V r^{\prime})\ d\V r^{\prime\prime}d\V r^{\prime},\label{eq:M_ij_Stokes}
\end{equation}
where $\Set G(\V r,\V r^{\prime})$ is the the Green's function for
the fluid equation, i.e., $\v(\V r)=\left(\L^{-1}\f\right)(\r)=\int\Set G(\V r,\V r^{\prime})\f(\r^{\prime})d\r^{\prime}$,
where
\begin{align}
\frac{\rho}{\D t}\v+\grad\p & -\eta\grad^{2}\v=\V f,\label{eq:Brinkman_Greens}\\
\div\v & =0.\nonumber 
\end{align}
It is well-known that $\Set G$ has the same form as (\ref{eq:M_tilde_ij}),
\[
\Set G(\Q_{1},\Q_{2})=f\left(r_{12}\right)\Id+g\left(r_{12}\right)\hat{\r}_{12}\otimes\hat{\r}_{12}.
\]

For steady Stokes flow ($\beta\rightarrow\infty$), $\Set G\equiv\sM O$
is the well-known Oseen tensor or Stokeslet %
\footnote{Observe that the regularized Stokeslet of Cortez \citet{RegularizedStokeslets}
is similar to (\ref{eq:M_ij_Stokes}) but contains only one regularized
delta function in the integrand; this makes the resulting mobility
matrix asymmetric, which is unphysical.%
}, and corresponds to $f_{S}(r)=g_{S}(r)\approx\left(8\pi\eta r\right)^{-1}$.
For inviscid flow, $\beta=0$, and we have that $\A=\left(\rho/\D t\right)\Id$
and (\ref{eq:L_inv_periodic}) applies, and therefore $\L^{-1}=\left(\Delta t/\rho\right)\sM P$
is a multiple of the projection operator. For finite nonzero values
of $\beta$, we can obtain $\Set G$ from the solution of the screened
Stokes (i.e., Brinkman) equations (\ref{eq:Brinkman_Greens}) \citet{Brinman_Original,RPY_Brinkman,RegularizedBrinkmanlet},
and corresponds to the ``Brinkmanlet'' \citet{RPY_Brinkman,RegularizedBrinkmanlet}
\begin{align}
f_{B}(r) & =\frac{e^{-\alpha r}}{4\pi\eta r}\left(\left(\frac{1}{\alpha r}\right)^{2}+\frac{1}{\alpha r}+1\right)-\frac{1}{4\pi\eta\alpha^{2}r^{3}},\label{eq:Brinkmanlet_3D}\\
g_{B}(r) & =-\frac{e^{-\alpha r}}{4\pi\eta r}\left(3\left(\frac{1}{\alpha r}\right)^{2}+\frac{3}{\alpha r}+1\right)+\frac{3}{4\pi\eta\alpha^{2}r^{3}},\nonumber 
\end{align}
where $\alpha^{2}=\rho/\left(\eta\D t\right)=\left(\beta h^{2}\right)^{-1}$.
Note that in the steady Stokes limit, $\alpha\rightarrow0$ and the
Brinkmanlet becomes the Stokeslet.

We can use (\ref{eq:Brinkmanlet_3D}) to construct $\widetilde{\Mob}_{ij}$
when the markers are far apart. Namely, if $r_{ij}\gg h$, then we
may approximate the IB kernel function by a true delta function, and
thus $f_{\beta}(r)$ and $g_{\beta}(r)$ are well-approximated by
(\ref{eq:Brinkmanlet_3D}). For steady Stokes flow, the interaction
between markers decays like $r^{-1}$. For finite $\beta$, however,
the viscous contribution decays exponentially fast as $\exp\left(-r/\left(h\sqrt{\beta}\right)\right)$,
which is consistent with the fact that markers interact via viscous
forces only if they are at a distance not much larger than $h\sqrt{\beta}=\sqrt{\nu\D t}$,
the typical distance that momentum diffuses during a time step. For
nonzero Reynolds numbers, the leading order asymptotic $r^{-3}$ decay
of $f_{\beta}(r)$ and $g_{\beta}(r)$ is given by the last terms
on the right hand side of (\ref{eq:Brinkmanlet_3D}) and corresponds
to the electric field of an electric dipole; its physical origin is
in the incompressibility constraint, which instantaneously propagates
hydrodynamic information between the markers %
\footnote{In reality, of course, this information is propagated via fast sound
waves and not instantaneously.%
}.

For steady Stokes flow, we can say even more about the approximate
form of $f_{\beta}(r)$ and $g_{\beta}(r)$. As discussed in more
detail by Delong et al. \citet{BrownianBlobs}, for distances between
the markers that are not too small compared to the regularization
length $a$, we can approximate (\ref{eq:M_ij_Stokes}) with (\ref{eq:M_tilde_ij})
using the well-known Rotne-Prager-Yamakawa (RPY) \citet{RotnePrager,RPY_FMM,RPY_Shear_Wall}
tensor for the functions $f_{\beta}(r)$ and $g_{\beta}(r)$,

\begin{align}
f_{RPY}(r) & =\frac{1}{6\pi\eta a}\begin{cases}
\frac{3a}{4r}+\frac{a^{3}}{2r^{3}}, & \quad r>2a,\\
1-\frac{9r}{32a}, & \quad r\leq2a,
\end{cases}\label{eq:RPYTensor}\\
g_{RPY}(r) & =\frac{1}{6\pi\eta a}\begin{cases}
\frac{3a}{4r}-\frac{3a^{3}}{2r^{3}}, & \quad r>2a,\\
\frac{3r}{32a}, & \quad r\leq2a,
\end{cases}\nonumber 
\end{align}
where $a$ is the effective \emph{hydrodynamic radius} of the specific
kernel $\delta_{a}$, defined by $\left(6\pi a\right)^{-1}=\int\delta_{a}(\V r^{\prime\prime})\sM O(\V r^{\prime\prime},\V r^{\prime})\delta_{a}(\V r^{\prime})\ d\V r^{\prime\prime}d\V r^{\prime}.$
Note that for $r\gg a$ the RPY tensor approaches the Oseen tensor
and decays like $r^{-1}$. A key advantage of the RPY tensor is that
it \emph{guarantees} that the mobility matrix (\ref{eq:M_tilde_ij})
is SPD for \emph{all} configurations of the markers, which is a rather
nontrivial requirement \citet{RPY_Shear_Wall}. The actual discrete
pairwise mobility $\Mob_{ij}$ obtained from the spatially-discrete
IB method is well-described by the RPY tensor \citet{BrownianBlobs}
(see Fig. \ref{fig:fg_3d_asympt}). The only fitting parameter in
the RPY approximation is the effective hydrodynamic radius $a$ averaged
over many positions of the marker relative to the underlying grid
\citet{ISIBM,BrownianBlobs}; for the six-point kernel used here %
\footnote{As summarized in Refs. \citet{ISIBM,BrownianBlobs}, $a\approx1.25h$
for the widely used four-point kernel \citet{IBM_PeskinReview}, and
$a\approx0.91h$ for the three-point kernel \citet{StaggeredIBM}.%
}, $a=1.47\, h$. For the Brinkman equation, the equivalent of the
RPY tensor can be computed for $r\geq2a$ by applying a Faxen-like
operator from the left and right on the Brinkmanlet (see Eq. (26)
in Ref. \citet{RPY_Brinkman}); the resulting analytical expressions
are complex and are not used in our empirical fitting.

\subsection{Empirical Fits}

In this work, we use empirical fits to approximate the mobility. This
is because the analytical approximations, such as those offered by
the RPY tensor, are most appropriate for unbounded domains and assume
the markers are far apart compared to the width of the regularized
delta function. In numerical computations, we use a finite periodic
domain, and this requires corrections to the analytic expressions
that are difficult to model. For example, for finite $\beta$, we
find that the periodic corrections to the inviscid (dipole) $r^{-3}$
contribution dominate over the exponentially decaying viscous contribution,
which makes the precise form of the viscous terms in (\ref{eq:Brinkmanlet_3D})
irrelevant in practice. For $r\gg h$, only the asymptotically-dominant
far-field terms survive, and we make an effort to preserve those in
our fitting because the numerical results are obtained using finite
systems and thus not reliable at large marker distances. At shorter
distances, however, the discrete nature of the fluid solver and the
IB kernel functions becomes important, and empirical fitting seems
to be a simple yet flexible alternative to analytical computations.
At the same time, we feel that is important to \emph{constrain} the
empirical fits based on known behavior at short and large distances.

Firstly, for $r\ll h$, the pairwise mobility can be well-approximated
by the self-mobility ($r=0$, corresponding to the diagonal elements
$\widetilde{\Mob}_{ii}$), for which we know the following facts:
\begin{itemize}
\item For the steady Stokes regime ($\beta\rightarrow\infty$), the diagonal
elements are given by Stokes's drag formula, yielding 
\[
f_{\infty}(0)=\left(6\pi\eta a\right)^{-1}\sim\ 1/\eta h\mbox{ and }g_{\infty}(0)=0,
\]
where we recall that $a$ is the effective hydrodynamic radius of
a marker for the particular spatial discretization (kernel and fluid
solver).
\item For the inviscid case ($\beta=0$), it is not hard to show that \citet{ISIBM}
\begin{equation}
f_{0}(0)=\frac{d-1}{d}\frac{\D t}{\rho}\BV^{-1}\sim\beta/\eta h\mbox{ and }g_{0}(0)=0,\label{eq:f_0_inviscid}
\end{equation}
where $d=3$ is the dimensionality, and $\BV=c_{V}h^{3}$ is the ``volume''
of the marker, where the constant $c_{V}$ is straightforward to calculate.
\item The above indicates that $f_{\beta}(0)$ goes from $\sim\beta/\left(\eta h\right)$
for small $\beta$ to $\sim1/\left(\eta h\right)$ for large $\beta$.
At intermediate viscous CFL numbers $\beta$, we can set
\begin{equation}
f_{\beta}\left(0\right)=\frac{C\left(\beta\right)}{\eta h}\mbox{ and }g_{\beta}(0)=0,\label{eq:f_beta_r_0}
\end{equation}
where $C\left(\beta\ll1\right)\approx2\beta/(3c_{V})$ is linear for
small $\beta$ and then becomes $O(1)$ for large $\beta$. We will
obtain the actual form of $C\left(\beta\right)$ from empirical fitting.
\end{itemize}
Secondly, for $r\gg h$, we know the asymptotic decay of the hydrodynamic
interactions from (\ref{eq:Brinkmanlet_3D}):
\begin{itemize}
\item For the steady Stokes regime ($\beta\rightarrow\infty$), we have
the Oseen tensor given by 
\begin{equation}
f_{\infty}(r\gg h)\approx g_{\infty}(r\gg h)\approx\left(8\pi\eta r\right)^{-1}.\label{eq:fg_asympt_3D}
\end{equation}

\item For the inviscid case ($\beta=0$), we get the electric field of an
electric dipole,
\begin{equation}
f_{0}(r\gg h)\approx-\frac{\D t}{4\pi\rho r^{3}}\mbox{ and }g_{0}(r\gg h)\approx\frac{3\D t}{4\pi\rho r^{3}},\label{eq:fg_inv_3D}
\end{equation}
 which is also the asymptotic decay for $\beta>0$ for $r\gg h\sqrt{\beta}$.
\end{itemize}
We obtain the actual form of the functions $f_{\beta}(r)$ and $g_{\beta}(r)$
empirically by fitting numerical data for the parallel and perpendicular
mobilities 
\begin{align*}
\mu_{ij}^{\parallel}=\hat{\r}_{ij}^{T}\widetilde{\Mob}_{ij}\hat{\r}_{ij} & \approx f_{\beta}\left(r_{ij}\right)+g_{\beta}\left(r_{ij}\right),\\
\mu_{ij}^{\perp}=\left(\hat{\r}_{ij}^{\perp}\right)^{T}\widetilde{\Mob}_{ij}\hat{\r}_{ij}^{\perp} & \approx f_{\beta}\left(r_{ij}\right),
\end{align*}
where $\hat{\r}_{ij}^{\perp}\cdot\hat{\r}_{ij}=0$. To do so, we placed
a large number of markers $N$ in a cube of length $l/8$ inside a
periodic domain of length $l$. For each marker $i$, we applied a
unit force $\Lamb_{i}$ with random direction while leaving $\Lamb_{j}=0$
for $j\neq i$, solved (\ref{eq:Stokes_saddle}), and then interpolated
the fluid velocity $\v$ at the position of each of the markers. The
resulting parallel and perpendicular relative velocity for each of
the $N(N-1)/2$ pairs of particles allows us to estimate $f_{\beta}\left(r_{ij}\right)$
and $g_{\beta}\left(r_{ij}\right)$. By making the number of markers
$N$ sufficiently large, we sample the mobility over essentially all
relative positions of the pair of markers. For the self-mobility $\widetilde{\Mob}_{ii}$
($r_{ii}=0$), we take $g_{\beta}(0)=0$ and compute $f_{\beta}(0)$
from the numerical data.

\begin{figure}[tbph]
\begin{centering}
\includegraphics[width=0.49\textwidth]{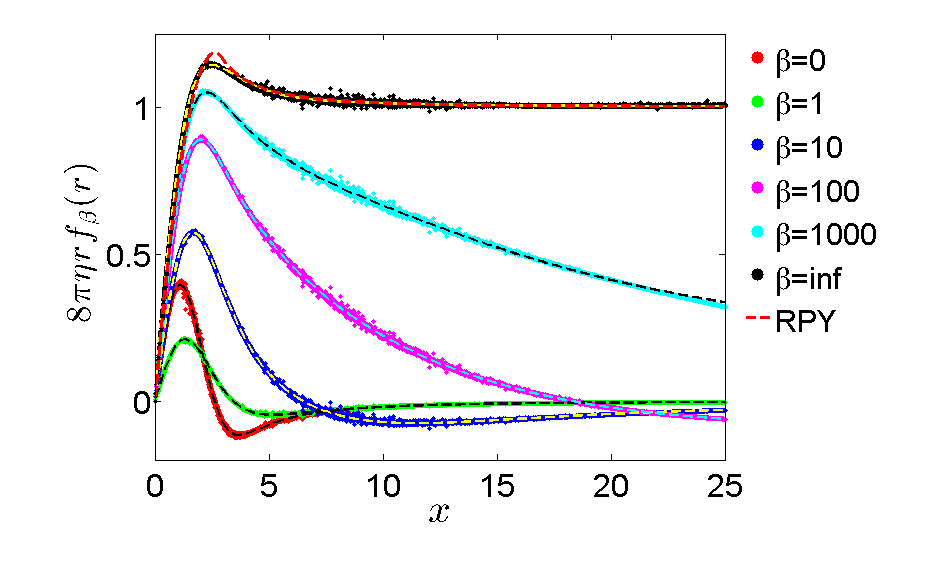}\includegraphics[width=0.49\textwidth]{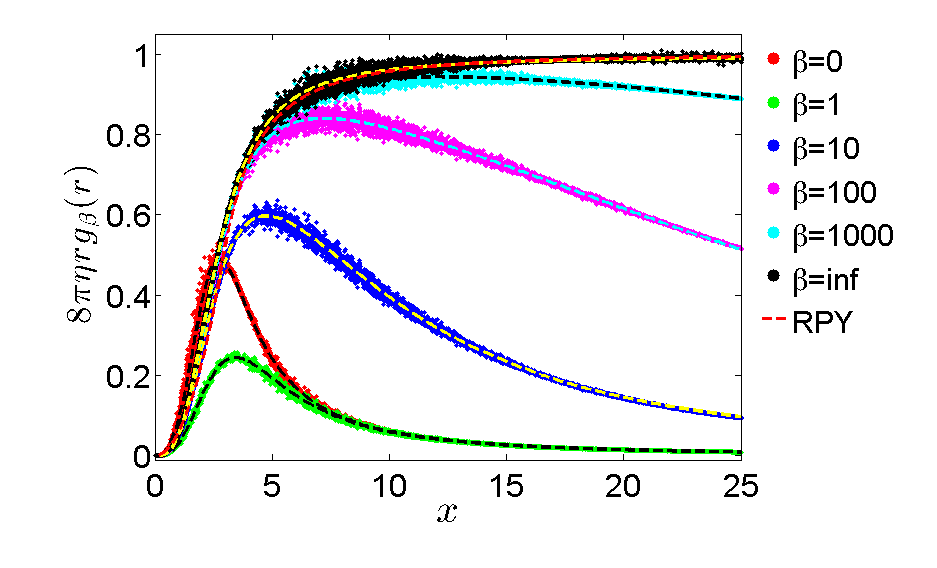}
\par\end{centering}

\caption{\label{fig:fg_3d_asympt}Normalized mobility functions $\tilde{f}(x)$
(left) and $\tilde{g}(x)$ (right) defined similarly to (\ref{eq:fg_tilde_Stokes})
as a function of marker-marker distance $x=r/h$, in three dimensions
for the six-point kernel of Bao et al. \citet{New6ptKernel}, over
a range of viscous CFL numbers (different colors, see legend). Numerical
data is shown with symbols and obtained using a $256^{3}$ periodic
fluid grid, while dashed lines show our empirical fit of the form
(\ref{eq:f_stokes_3D_fit}) for steady Stokes ($\beta\rightarrow\infty$)
and (\ref{eq:f_beta_3D_fit}) for finite $\beta$. For steady Stokes
flow, the numerical data is in reasonable agreement with the RPY tensor
(\ref{eq:RPYTensor}) (dashed red line). }
\end{figure}

If the spatial discretization were perfectly translationally and rotationally
invariant and the domain were infinite, all of the numerical data
points for $f_{\beta}\left(r\right)$ and $g_{\beta}\left(r\right)$
would lie on a smooth curve and would not depend on the actual position
of the pair of markers relative to the underlying grid. In reality,
it is not possible to achieve perfect translational invariance with
a kernel of finite support \citet{IBM_PeskinReview}, and so we expect
some (hopefully small) scatter of the points around a smooth fit.
Normalized numerical data for $f_{\beta}(r)$ and $g_{\beta}(r)$
are shown in Fig. \ref{fig:fg_3d_asympt}, and we indeed see that
the data can be fit well by smooth functions over the whole range
of distances. To maximize the quality of the fit, we perform separate
fits for $\beta\rightarrow\infty$ (steady Stokes flow) and finite
$\beta$. We also make an effort to make the fits change smoothly
as $\beta$ grows towards infinity, as we explain in more detail in
Appendix \ref{AppendixMob3D}. Code to evaluate the empirical fits
described in Appendices \ref{AppendixMob3D} and \ref{AppendixMob2D}
is publicly available to others for a number of kernels constructed
by Peskin and coworkers (three-, four-, and six-point) in both two
and three dimensions at \url{http://cims.nyu.edu/~donev/src/MobilityFunctions.c}.

\section{\label{sec:Solver}Linear Solver}

To solve the constrained Stokes problem (\ref{eq:constrained_Stokes}),
we use the preconditioned flexible GMRES (FGMRES) method, which is
a Krylov solver. We will refer to this as the ``outer'' Krylov solver,
as it must be distinguished from ``inner'' Krylov solvers used in
the preconditioner. Because we use Krylov solvers in our preconditioner
and because Krylov solvers generally cannot be expressed as linear
operators, it is crucial to use a flexible Krylov method such as FGMRES
for the outer solver. The overall method is implemented in the open-source
immersed-boundary adaptive mesh refinement (IBAMR) software infrastructure
\citet{IBAMR}; in this work we focus on uniform grids and do not
use the AMR capabilities of IBAMR (but see \citet{RigidIBAMR,IBAMR_Fish}).
IBAMR uses Krylov solvers that are provided by the PETSc library \citet{PETSc}.

\subsection{\label{sub:Preconditioner}Preconditioner for the constrained Stokes
system}

In the preconditioner used by the outer Krylov solver, we want to
\emph{approximately} solve the nested saddle-point linear system
\[
\left[\begin{array}{ccc}
\A & \Grad & -\S\\
-\Div & \M 0 & \M 0\\
-\J & \M 0 & \V 0
\end{array}\right]\left[\begin{array}{c}
\v\\
\p\\
\Lamb
\end{array}\right]=\left[\begin{array}{c}
\V g\\
\V h\\
\W
\end{array}\right],
\]
where we recall that $\A=\left(\rho/\D t\right)\Id-\eta h^{-2}\Lap$.
Let us set $\alpha=1$ if $\A$ has a null-space, (e.g., for a fully
periodic domain for steady Stokes flow) and we set $\alpha=0$ if
$\A$ is invertible. When $\alpha=1$, let us define the restricted
inverse $\A^{-1}$ to only act on vectors of mean value zero, and
to return a vector of mean zero.

Applying our Schur complement based preconditioner for solving (\ref{eq:constrained_Stokes})
consists of the following steps:
\begin{enumerate}
\item Solve the (unconstrained) fluid sub-problem,
\[
\left[\begin{array}{cc}
\A & \Grad\\
-\Div & \M 0
\end{array}\right]\left[\begin{array}{c}
\v\\
\p
\end{array}\right]=\left[\begin{array}{c}
\V g\\
\V h
\end{array}\right].
\]
To control the accuracy of the solution one can either use a relative
tolerance based stopping criterion or fix the number of iterations
$N_{s}$ in the inner solver.
\item Calculate the slip velocity on the set of markers, $\D{\U}=-\left(\J\v+\W\right)$.
\item Approximately solve the Schur complement system, 
\begin{equation}
\widetilde{\Mob}\Lamb=\D{\U},\label{eq:mob_subproblem}
\end{equation}
where the mobility approximation $\widetilde{\Mob}$ is constructed
as described in Section \ref{sec:ApproxMob}.
\item Optionally, re-solve the corrected fluid sub-problem,
\[
\left[\begin{array}{cc}
\A & \Grad\\
-\Div & \M 0
\end{array}\right]\left[\begin{array}{c}
\v\\
\p
\end{array}\right]=\left[\begin{array}{c}
\V g+\S\Lamb-\alpha\text{\ensuremath{\DV}}^{-1}\V 1^{T}\Lamb\\
\V h
\end{array}\right].
\]

\end{enumerate}
All linear solvers used in the preconditioner can be approximate,
and this is in fact the key to the efficiency of the overall solver
approach. Notably, the inner Krylov solvers used to solve the unconstrained
Stokes sub-problems in steps 1 and 4 above can be done by using a
small number $N_{s}$ of iterations using a method briefly described
in the next section. If the fluid sub-problem is approximately solved
in both steps 1 and 4, which we term the \emph{full} Schur complement
preconditioner, each application of the preconditioner requires $2N_{s}$
applications of the Stokes preconditioner (\ref{P1_practical}). It
is also possible to omit step 4 above to obtain a block \emph{lower
triangular} Schur preconditioner \citet{Elman_FEM_Book}, which requires
only $N_{s}$ applications of the unconstrained Stokes preconditioner
(\ref{P1_practical}). We will numerically compare these two preconditioners
and study the effect of $N_{s}$ on the convergence of the FGMRES
outer solver in Section \ref{sub:GMRES}.

\subsection{Unconstrained Fluid Solver}

A key component we rely on is an approximate solver for the unconstrained
Stokes sub-problem,
\[
\left[\begin{array}{cc}
\A & \Grad\\
-\Div & \M 0
\end{array}\right]\left[\begin{array}{c}
\v\\
\p
\end{array}\right]=\left[\begin{array}{c}
\V g\\
\V h
\end{array}\right],
\]
for which a number of techniques have been developed in the finite-element
context \citet{Elman_FEM_Book}. To solve this system, we use GMRES
with a preconditioner $\sM P_{S}^{-1}$ based on the projection method,
as proposed by Griffith \citet{NonProjection_Griffith} and improved
to some extent by Cai et al. \citet{StokesKrylov}. Specifically,
the preconditioner for the Stokes system that we use in this work
is
\begin{equation}
\sM P_{S}^{-1}=\left(\begin{array}{cc}
\Id & h^{2}\Grad\widetilde{\Lp}^{-1}\\
\V 0 & \widetilde{\M{\mathcal{B}}}^{-1}
\end{array}\right)\left(\begin{array}{cc}
\Id & \V 0\\
-\Div & -\Id
\end{array}\right)\left(\begin{array}{cc}
\widetilde{\A}^{-1} & \V 0\\
\V 0 & \Id
\end{array}\right),\label{P1_practical}
\end{equation}
where $\Lp=h^{2}\left(\Div\Grad\right)$ is the dimensionless pressure
(scalar) Laplacian, and $\widetilde{\A}^{-1}$ and $\widetilde{\Lp}^{-1}$
denote approximate solvers obtained by a \emph{single} V-cycle of
a geometric multigrid solver for the vector Helmholtz and scalar Poisson
problems, respectively. In the time-dependent case, the approximate
Schur complement for the unconstrained Stokes sub-problem is
\[
\widetilde{\M{\mathcal{B}}}^{-1}=-\frac{\rho h^{2}}{\D t}\widetilde{\Lp}^{-1}+\eta\Id,
\]
and for steady Stokes flow, $\widetilde{\M{\mathcal{B}}}^{-1}=\eta\Id$.
Further discussion of the relation of these preconditioners to the
those described in the book \citet{Elman_FEM_Book} can be found in
\citet{NonProjection_Griffith}.

Observe that one application of $\sM P_{S}^{-1}$ is relatively inexpensive
and involves only a few scalar multigrid V-cycles. Indeed, solving
the Stokes system using GMRES with this preconditioner is only a few
times more expensive than solving a scalar Poisson problem, even in
the steady Stokes regime \citet{StokesKrylov}. Note that it is possible
to omit the upper right off-diagonal block in the first matrix on
the right hand side of (\ref{P1_practical}) to obtain a block lower
triangular preconditioner that is also effective, and may in fact
be preferred at zero Reynolds number since it allows one to skip a
sweep of the pressure multigrid solver \citet{StokesKrylov}. We empirically
find that including the Poisson solve (velocity projection) improves
the overall performance of the outer solver.

\subsection{Mobility Solver}

From a computational perspective, one of the most challenging steps
in our preconditioner is solving the mobility sub-problem (\ref{eq:mob_subproblem}).
Since this is done inside a preconditioner, and because $\widetilde{\Mob}$
is itself an approximation of the true mobility matrix $\Mob$, it
is not necessary to solve (\ref{eq:mob_subproblem}) exactly. In the
majority of the examples presented herein, we solve (\ref{eq:mob_subproblem})
using direct solvers provided by LAPACK. This is feasible on present
hardware for up to around $10^{5}$ markers and allows us to focus
on the design of the approximation $\widetilde{\Mob}$ and to study
the accuracy of the overall method.

Let us denote with $s$ the smallest marker-marker spacing. For well-spaced
markers, $s/h\gtrapprox2$, our approximate mobility $\widetilde{\Mob}$
is typically SPD even for large numbers of markers, and in these cases,
we can use the Cholesky factorization to solve (\ref{eq:mob_subproblem}).
In some cases, however, there may be a few small or even negative
eigenvalues of $\widetilde{\Mob}$ that have to be handled with care.
We have found that the most robust (albeit expensive) alternative
is to perform an SVD of $\widetilde{\Mob}$, and to use a pseudoinverse
of $\widetilde{\Mob}$ (keeping only eigenvalues larger than some
tolerance $\epsilon_{SVD}>0$) to solve (\ref{eq:mob_subproblem}).
This effectively filters out the spuriously small or negative eigenvalues.
Note that the factorization of $\widetilde{\Mob}$ needs to be performed
only once per constrained Stokes solve since the body is kept fixed
during a time step. In cases where there is a single body, the factorization
needs to be performed only once per simulation and can be reused;
if the body is translating or rotating, one ought to perform appropriate
rotations of the right hand side and solution of (\ref{eq:mob_subproblem}).
In some cases of practical interest where the number of markers is
not too large, it is possible to precompute the true mobility $\Mob_{0}$
with periodic boundary conditions (for a sufficienly large domain)
and to store its factorization. Even if the structure moves relative
to the underlying grid, such a precomputed (reference) mobility $\Mob_{0}$
is typically a much better approximation to the true mobility than
our empirical approximation $\widetilde{\Mob}$, and can effectively
be used in the preconditioner. Determining effective approaches to
solving the mobility sub-problem in the presence of multiple moving
rigid bodies remains future work, as discussed further in the Conclusions.

\section{\textmd{\textup{\normalsize \label{sec:Conditioning}}}Conditioning
of the mobility matrix}

The conditioning of the constrained Stokes problem (\ref{eq:constrained_Stokes})
is directly related to the conditioning of the Schur complement mobility
matrix $\Mob=\J\L^{-1}\S$, which is intimately connected to the relation
between the fluid solver grid spacing $h$ and the smallest inter-marker
spacing $s$. Firstly, it is obvious that if two markers $i$ and
$j$ are very close to each other, then the fluid solver cannot really
distinguish between $\Lamb_{i}$ and $\Lamb_{j}$ and will instead
effectively see only their sum. We also know that using too many markers
for a fixed fluid grid will ultimately lead to a rank-deficient $\Mob$,
because it is not possible to constrain a finite-dimensional discrete
fluid velocity at too many points. This physical intuition tells us
that the condition number of $\Mob$ should increase as the marker
spacing becomes small compared to the grid spacing. This well-known
intuition, however, does not tell us how closely the markers can or
must be placed in practice. Standard wisdom for the immersed boundary
method, which is based on the behavior of models of elastic bodies,
is to make the marker spacing on the order of half a grid spacing.
As we show, this leads to extremely ill-conditioned mobility matrices
for rigid bodies. We note that the specific results depend on the
dimensionality, the details of the fluid solver, and the specific
kernel used; however, the qualitative features we report appear to
be rather general.

\begin{figure}[tbph]
\begin{centering}
\includegraphics[width=0.9\textwidth]{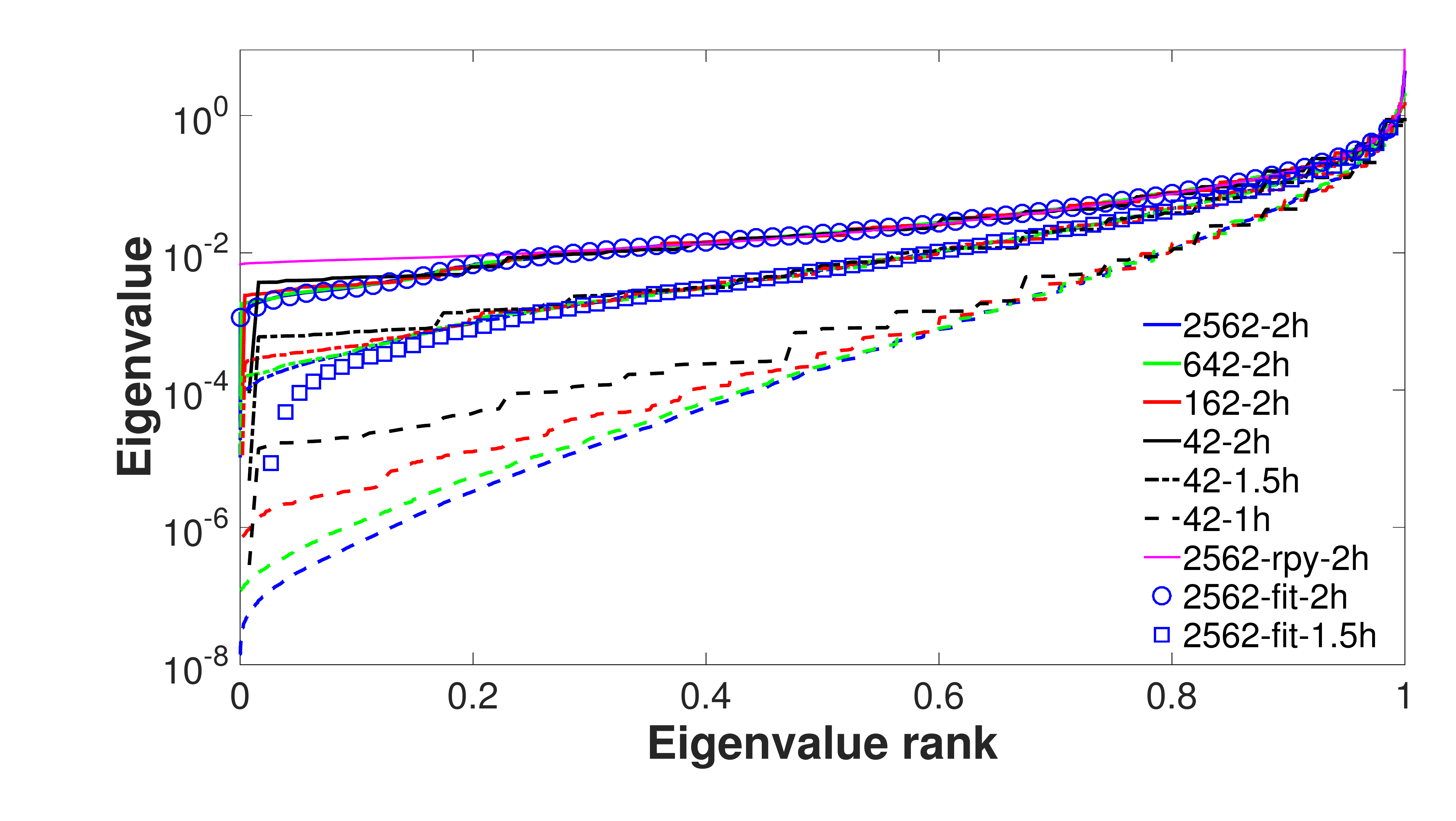}
\par\end{centering}

\begin{centering}

\par\end{centering}

\caption{\label{fig:Eigenvalues}Eigenvalue spectrum of the mobility matrix
for steady Stokes flow around a spherical shell covered with different
numbers of markers (42, 162, 642, or 2562, see legend) embedded in
a periodic domain. Solid lines are for marker spacing of $s\approx2h$,
dashed-dotted lines for spacing $s\approx1.5\, h$, and dashed lines
for spacing of $s\approx1h$. The marker spacing is $s\approx1$ in
all cases; for $s\approx2h$, the fluid grid size is $128^{3}$ for
2562 markers and $64^{3}$ for smaller number of markers, and scaled
accordingly for other spacings. For comparison, we show the spectrum
of $\widetilde{\Mob}_{RPY}$ for the most resolved model ($N=2562$
markers) at $s/h\approx2$. Also shown is the spectrum of the empirical
(fit) approximation to the mobility $\widetilde{\Mob}$ for the two
larger spacings; for $s\approx h$ our empirical approximation is
very poor and includes many spurious negative eigenvalues (not shown).}
\end{figure}

To determine the condition number of the mobility matrix, we consider
``open'' and ``filled'' sphere models. We discretize the surface
of a sphere as a shell of markers constructed by a recursive procedure
suggested to us by Charles Peskin (private communication). We start
with 12 markers placed at the vertices of an icosahedron, which gives
a uniform triangulation of a sphere by 20 triangular faces. Then,
we place a new marker at the center of each edge and recursively subdivide
each triangle into four smaller triangles, projecting the vertices
back to the surface of the sphere along the way. Each subdivision
approximately quadruples the number of vertices, with the $k$-th
subdivision producing a model with $10\cdot4^{k-1}+2$ markers. To
create filled sphere models, we place additional markers at the vertices
of a tetrahedral grid filling the sphere that is constructed using
the TetGen library, starting from the surface triangulation described
above. The constructed tetrahedral grids are close to uniform, but
it is not possible to control the precise marker distances in the
resulting irregular grid of markers. We use models with approximately
equal edges (distances between nearest-neighbor markers) of length
$\approx s$, which we take as a measure of the typical marker spacing.
We numerically computed the mobility matrix $\Mob$ for an isolated
spherical shell in a large periodic domain for various numbers of
markers $N$. Here we keep the ratio $s/h$ fixed and keep the marker
spacing fixed at $s\approx1$; one can alternatively keep the radius
of the sphere fixed %
\footnote{The scaling used here, keeping $s=1$ fixed, is more natural for examining
the small eigenvalues of $\Mob$, which are dominated by discretization
effects, as opposed to the large eigenvalues, which correspond to
physical modes of the Stokes problem posed on a sphere and are insensitive
to the discretization details.%
}. In Fig. \ref{fig:Eigenvalues}, we show the spectrum of $\Mob$
for varying levels of resolution for three different spacings of the
markers, $s/h=1$, $s/h=3/2$, and $s/h=2$. Similar spectra, but
with somewhat improved condition number (i.e., fewer smaller eigenvalues),
are seen for nonzero Reynolds numbers (finite $\beta$).

The results in Fig. \ref{fig:Eigenvalues} strongly suggest that as
the number of markers increases, the low-lying (small eigenvalue)
spectrum of the mobility matrix approaches a limiting shape. Therefore,
the nontrivial eigenvalues remain bounded away from zero even as the
resolution is increased, which implies that for $s/h\gtrsim1$ the
system (\ref{eq:constrained_Stokes}) is uniformly solvable or ``stable''
under grid refinement. Note that in the case of a sphere, there is
a trivial zero eigenvalue in the continuum limit, which corresponds
to uniform compression of the sphere; this is reflected in the existence
of one eigenvalue much smaller than the rest in the discrete models.
Ignoring the trivial eigenvalue, the condition number of $\Mob$ is
$O(N)$ for this example because the largest eigenvalue in this case
increases like the number of markers $N$, in agreement with the fact
that the Stokes drag on a sphere scales linearly with its radius.
This is as close to optimal as possible, because for the continuum
equations for Stokes flow around a sphere, the eigenvalues corresponding
to spherical harmonic modes scale like the index of the spherical
harmonic. However, what we are concerned here is not so much how the
condition number scales with $N$, but with the size of the prefactor,
which is determined by the smallest nontrivial eigenvalues of $\Mob$.

Figure \ref{fig:Eigenvalues} clearly shows that the number of very
small eigenvalues increases as we bring the markers closer to each
other, as expected. The increase in the conditioning number is quite
rapid, and the condition number becomes $O(10^{6}-10^{7})$ for marker
spacings of about one per fluid grid cell. For the conventional choice
$s\approx h/2$, the mobility matrix is so poorly conditioned that
we cannot solve the constrained Stokes problem in double-precision
floating point arithmetic. Of course, if the markers are too far apart
then fluid will leak through the wall of the structure. We have performed
a number of heuristic studies of leak through flat and curved rigid
walls and concluded that $s/h\approx2$ yields both small leak and
a good conditioning of the mobillity, at least for the six-point kernel
used here \citet{New6ptKernel}. Therefore, unless indicated otherwise,
in the remainder of this work, we keep the markers about \emph{two
grid cells apart} in both two and three dimensions. It is important
to emphasize that this is just a heuristic recommendation and not
a precise estimate. We remark that Taira and Colonius, who solve a
different Schur complement ``modified Poisson equation'', recommend
$s/h\approx1$ to ``achieve a reasonable condition number and to
prevent penetration of streamlines caused by a lack of Lagrangian
points.'' 

It is important to observe that putting the markers further than the
traditional wisdom will increase the ``leak'' between the markers.
For rigid structures, the exact positioning of the markers can be
controlled since they do not move relative to one another as they
do for an elastic bodies; this freedom can be used to reduce penetration
of the flow into the body by a careful construction of the marker
grid. In the Conclusions, we discuss alternatives to the traditional
marker-based IB method \citet{IBFE} that can be used to control the
conditioning number of the Schur complement and allow for more tightly-spaced
markers.

It is worthwhile to examine the underlying cause of the ill-conditioning
as the markers are brought close together. One source of ill-conditioning
comes from the \emph{discrete} (finite-dimensional) nature of the
fluid solver, which necessarily limits the rank of the mobility matrix.
But another contributor to the worsening of the conditioning is the
\emph{regularization} of the delta function. Observe that for a true
delta function ($a\rightarrow0$) in Stokes flow, the pairwise mobility
is the length-scale-free Oseen tensor $\sim r^{-1}$, and the shape
of the spectrum of the mobility matrix has to be \emph{independent}
of the spacing among the markers. In the standard immersed boundary
method, $a\sim h$, so the fluid grid scale $h$ and the regularization
scale $a$ are difficult to distinguish.

To try to separate $h$ from $a$, we can take a continuum model of
the fluid, but keep the discrete marker representation of the body;
see (\ref{eq:semi_continuum}). In this case the pairwise mobility
would be given by (\ref{eq:M_ij_Stokes}), which leads to the RPY
tensor (\ref{eq:RPYTensor}) for a kernel that is a surface delta
function over a sphere of radius $a$ (see (4.1) in \citet{RPY_Shear_Wall}).
In Fig. \ref{fig:Eigenvalues} we compare the spectra of the discrete
mobility $\Mob$ with those of the analytical mobility approximation
$\widetilde{\Mob}_{RPY}$ constructed by using (\ref{eq:RPYTensor})
for the pairwise mobility. We observe that the two are very similar
for $s\approx2h$, however, for smaller spacings $\widetilde{\Mob}_{RPY}$
does not have very small eigenvalues and is much better conditioned
than $\Mob$ (data not shown). In Fig. \ref{fig:Eigenvalues} we also
show the spectrum of our approximate mobility $\widetilde{\Mob}$
constructed using the empirical fits described in Section \ref{sec:ApproxMob}.
The resulting spectra show a worsening conditioning for spacing $s\approx1.5h$
consistent with the spectrum of $\Mob$. These observations suggest
that both the regularization of the kernel and the discretization
artifacts contribute to the ill-conditioning, and suggest that it
is worthwhile to explore alternative discrete delta function kernels
in the context of rigid-body IB methods.

We also note that we see a severe worsening of the conditioning of
$\Mob$, independent of $\beta$, when we switch from a spherical
shell to a filled sphere model. Some of this may be due to the fact
that the tetrahedral volume mesh used to construct the marker mesh
is not as uniform as the surface triangular mesh. We suspect, however,
that this ill-conditioning is primarily \emph{physical} rather than
numerical, and comes from the fact that the present marker model cannot
properly distinguish between surface tractions and body (volume) stresses.
Therefore, $\Lamb$ remains physically ill-defined even if one gets
rid of all discretization artifacts. 

Lastly, it is important to emphasize that in the presence of ill-conditioning,
what matters in practice are not only the smallest eigenvalues but
also their associated eigenvectors. Specifically, we expect to see
signatures of these eigenvectors (modes) in $\Lamb$, since they will
appear with large coefficients in the solution of (\ref{eq:lambda_sol})
if the right hand side has a nonzero projection onto the corresponding
mode. As expected, the small-eigenvalue eigenvectors of the mobility
correspond to high-frequency (in the spatial sense) modes for the
forces $\Lamb$. Therefore, if the markers are too closely spaced
the solutions for the forces $\Lamb$ will develop unphysical high-frequency
oscillations or jitter, even for smooth flows, especially in time-dependent
flows, as observed in practice \citet{SmoothingDelta_IBM}. We have
observed that for \emph{smooth} flows (i.e., smooth right hand-side
of (\ref{eq:lambda_sol})), the improved translational invariance
of the 6-point kernel reduces the magnitude of this jitter compared
to the traditional Peskin four-point kernel.

\section{\label{sec:Results}Numerical Tests}

In this section we apply our rigid-body IB method to a number of benchmark
problems. We first present tests of the preconditioned FGMRES solver,
and then demonstrate the advantage of our method over splitting-based
direct forcing methods. We further consider a simple test problem
at zero Reynolds number, involving the flow around a fixed sphere,
and study the accuracy of both the fluid (Eulerian) variables $\v$
and $\p$, as well as of the body (Lagrangian) surface tractions represented
by $\Lamb$, as a function of the grid resolution. We finally study
flows around arrays of cylinders in two dimensions and spheres in
three dimensions over a range of Reynolds numbers, and compare our
results to those obtained by Ladd using the Lattice-Boltzmann method
\citet{VACF_Ladd,FiniteRe_2D_Ladd,SmallRe_3D_Ladd}.

\subsection{\label{sub:GMRES}Empirical convergence of GMRES}

Here we consider the model problem of flow past a sphere in a cubic
domain that is either periodic or with no-slip boundaries. Except
for the largest resolutions studied here, the number of markers is
relatively small, and dense linear algebra can be used to solve the
mobility subproblem (\ref{eq:mob_subproblem}) robustly and efficiently,
so that the cost of the solver is dominated by the fluid solver. We
therefore use the number of total applications of the Stokes preconditioner
(\ref{P1_practical}) as a proxy for the CPU effort, instead of relying
on elapsed time, which is both hardware and software dependent. A
key parameter in our preconditioner is the number of iterations $N_{s}$
used in the iterative unconstrained Stokes solver. We recall that
in the full preconditioner, there are two unconstrained inexact Stokes
solves per iteration, giving a total of $2N_{s}$ applications of
$\sM P_{S}^{-1}$ per outer FGMRES iteration. If the lower triangular
preconditioner is used, then the second inexact Stokes solve is omitted,
and we perform only $N_{s}$ applications of $\sM P_{S}^{-1}$ per
outer FGMRES iteration.

\begin{figure}[h]
\begin{centering}
\includegraphics[width=0.49\textwidth]{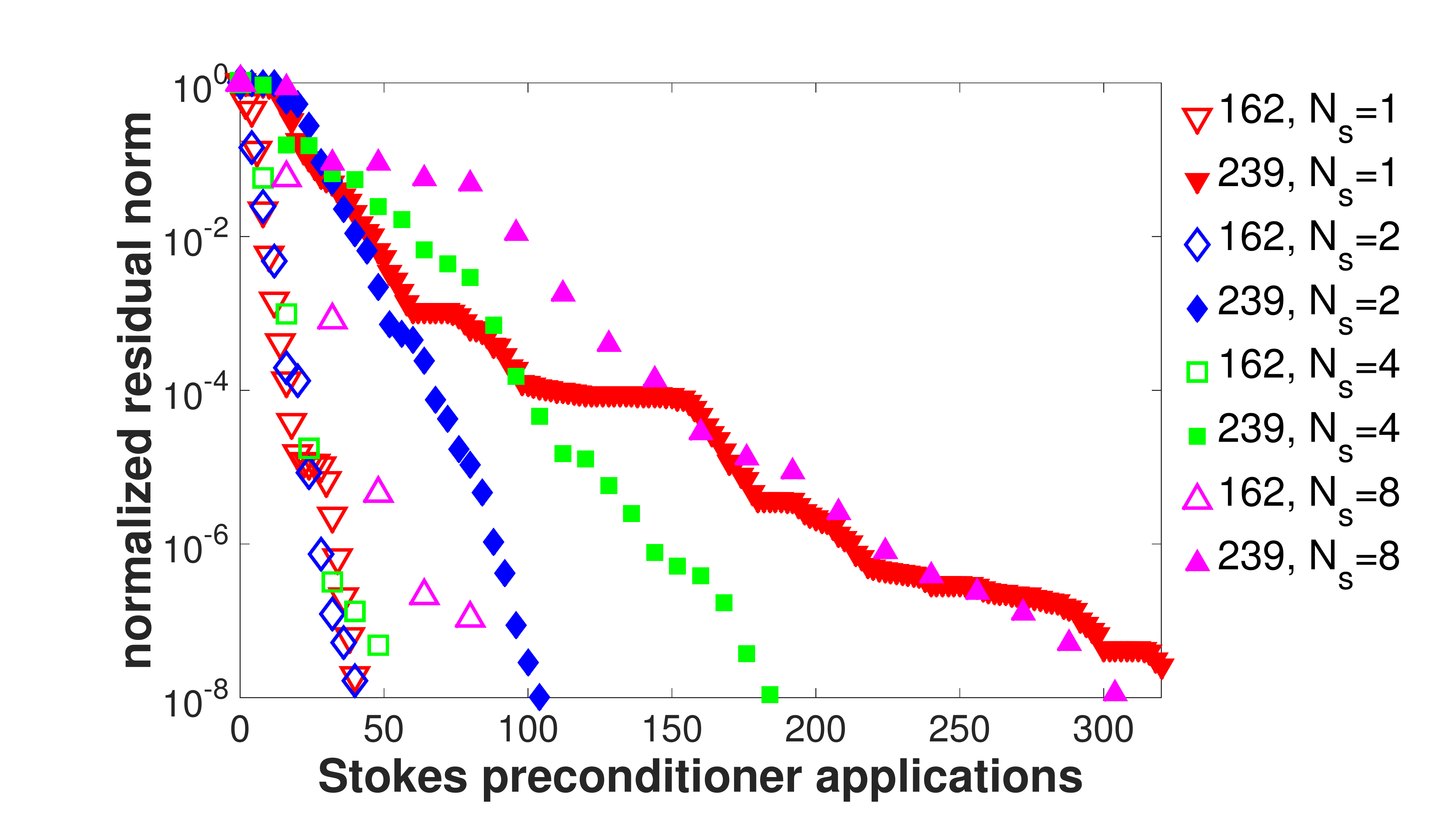}\includegraphics[width=0.49\textwidth]{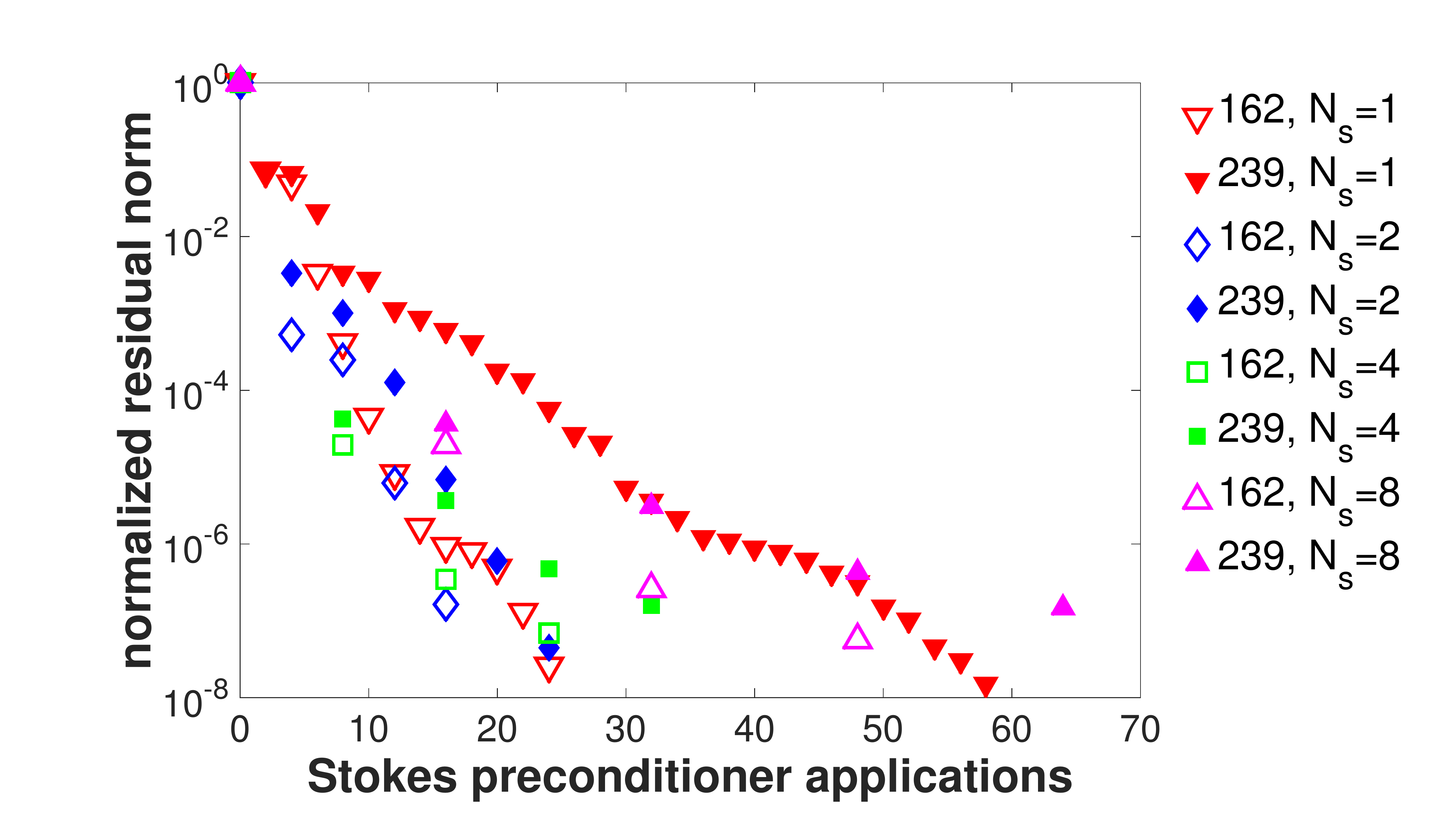}
\par\end{centering}

\begin{centering}
\includegraphics[width=0.49\textwidth]{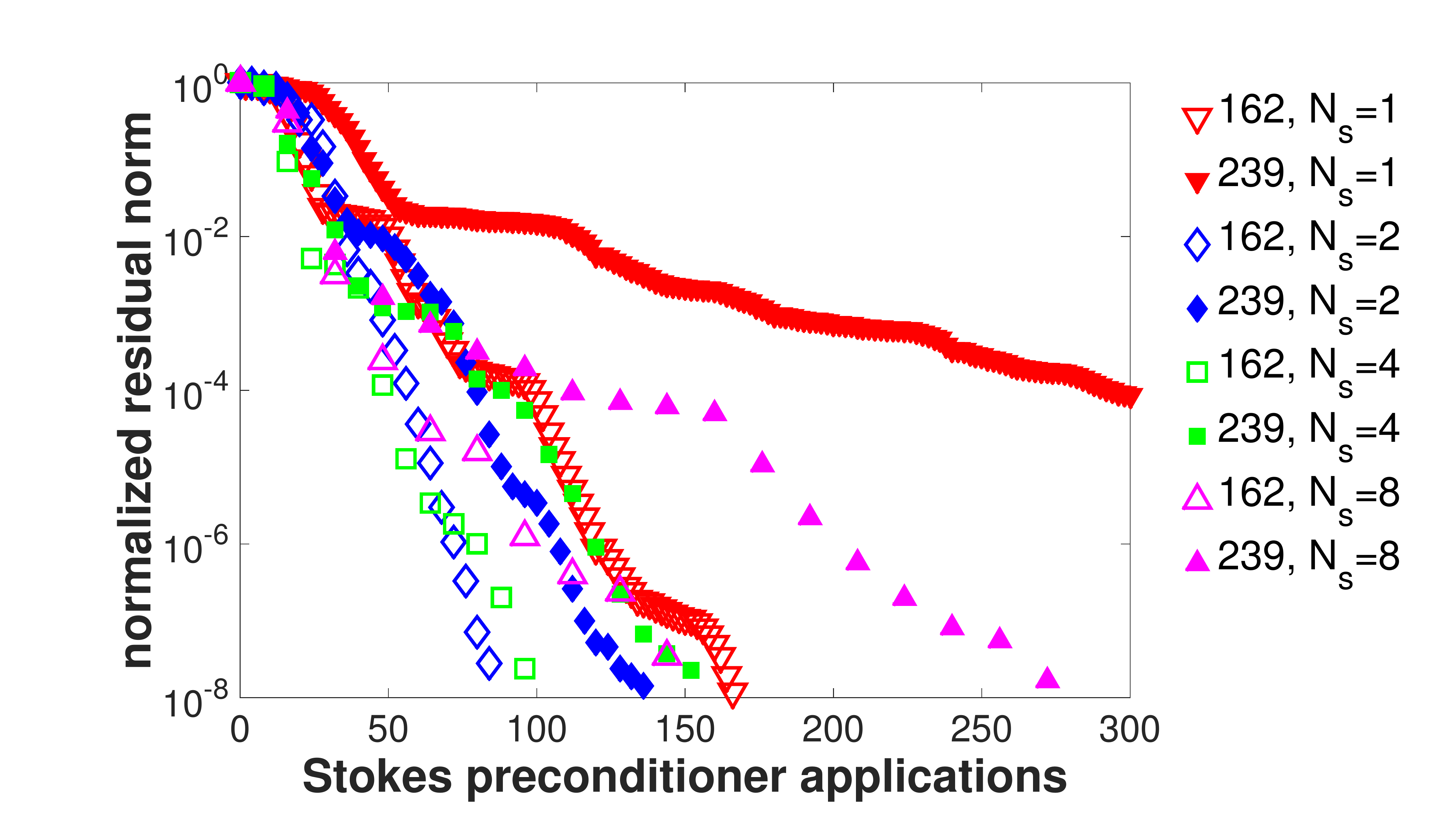}\includegraphics[width=0.49\textwidth]{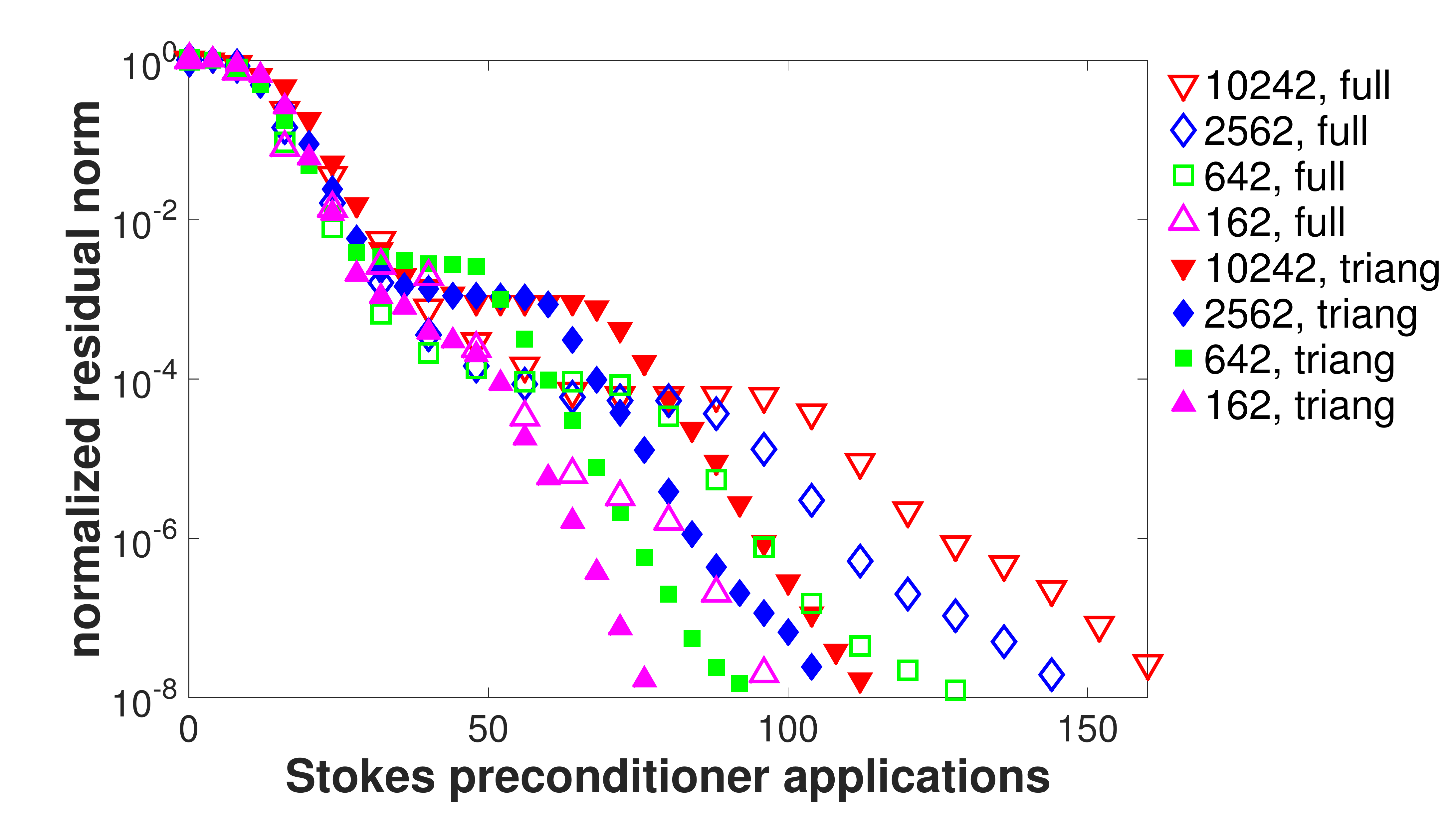}
\par\end{centering}

\centering{}\caption{\label{fig:GMRES_3D}FGMRES convergence for the constrained problem
(\ref{eq:constrained_Stokes}) for different numbers of iterations
$N_{s}$ in the unconstrained Stokes solver used in the preconditioner.
The specific problem is a rigid sphere of hydrodynamic radius $R$
moving through a stationary domain of length $L\approx4.35\, R$,
with the marker spacing fixed at $s/h\approx2$ and the GMRES restart
frequency set to 100 iterations. (Top left panel) Steady Stokes ($Re=0$)
flow for an empty 162-marker shell and a filled 239-marker sphere
moving through a periodic domain of $32^{3}$ fluid grid cells. (Top
right panel) Same as top left but for $Re\approx10$. (Bottom left
panel) As top left panel but now in domain with no slip boundary conditions
applied on all sides of the domain. (Bottom right panel) A spherical
shell moving in a non-periodic domain (as in bottom left panel) for
different resolutions of the shell (162, 642, 2562, and 10242 markers,
respectively) and the fluid solver grid ($32^{3}$, $64^{3}$, $128^{3}$,
and $256^{3}$ grid cells, respectively), fixing $N_{s}=4$, for both
the full Schur complement preconditioner and the lower triangular
approximate Schur complement preconditioner.}
\end{figure}

In the first set of experiments, we use the full preconditioner and
periodic boundary conditions. We represent the sphere by a spherical
shell of markers that is either empty (162 markers) or is filled with
additional markers in the interior (239 markers). The top panels of
Fig. \ref{fig:GMRES_3D} show the relative FGMRES residual as a function
of the total number of applications of $\sM P_{S}^{-1}$ for several
different choices of $N_{s}$, for both steady Stokes flow (left panel)
and a flow at Reynolds number $\Re=10$ (right panel). We see that
for spherical shells with well-conditioned $\Mob$ and $\widetilde{\Mob}$,
the exact value of $N_{s}$ does not have a large effect on solver
performance. However, making $N_{s}$ very large leads to wasted computational
effort by ``over-solving'' the Stokes system. This degrades the
overall performance, especially for tight solver tolerance. For the
ill-conditioned case of a filled sphere model in steady Stokes flow,
the exact value of $N_{s}$ strongly affects the performance, and
the optimal value is empirically determined to be $N_{s}=2$. As expected,
the linear system (\ref{eq:constrained_Stokes}) is substantially
easier to solve at higher Reynolds numbers, especially for the filled-sphere
models.

In the bottom left panel of Fig. \ref{fig:GMRES_3D} we show the FGMRES
convergence for a non-periodic system. In this case, we know that
the Stokes preconditioner $\sM P_{S}^{-1}$ itself does not perform
as well as in the periodic case \citet{NonProjection_Griffith,StokesKrylov},
and we expect slower overall convergence. In this case, we see that
$N_{s}=2$ and $N_{s}=4$ are good choices. Investigations (data not
shown) show that $N_{s}=4$ is more robust for problems with a larger
number of markers. Also, note that increasing $N_{s}$ decreases the
total number of FGMRES iterations for a fixed number of applications
of $\sM P_{S}^{-1}$, and therefore reduces the overall memory usage
and the number of times the mobility subproblem (\ref{eq:mob_subproblem})
needs to be solved; however, note that each of these solves is just
a backward/forward substitution if a direct factorization of $\widetilde{\Mob}$
has been precomputed.

The bottom right panel of Fig. \ref{fig:GMRES_3D} shows the FGMRES
convergence for a non-periodic system as the resolution of the grid
and the spherical shell is refined in unison, keeping $N_{s}=4$.
The results in Fig. \ref{fig:GMRES_3D} demonstrate that our linear
solver is able to cope with the increased number of degrees of freedom
under refinement relatively robustly, although a slow increase of
the total number of FGMRES iterations is observed. Comparing the full
preconditioner with the lower triangular preconditioner, we see that
the latter is computationally more efficient overall; this is in agreement
with experience for the unconstrained Stokes system \citet{StokesKrylov}.
In some sense, what this shows is that it is best to let the FGMRES
solver correct the initial unconstrained solution for the velocity
and pressure in the \emph{next} FGMRES iteration, rather than to re-solve
the fluid problem in the preconditioner itself. However, if very tight
solver tolerance is required, we find that it is necessary to perform
some corrections of the velocity and pressure inside the preconditioner.
In principle, the second unconstrained Stokes solve in the preconditioner
can use a different number of iterations $N_{s}^{\prime}$ from the
first, but we do not explore this option further here. Also note that
if $\widetilde{\Mob}\approx\Mob$ (for example, if it was computed
numerically rather than approximated), then the full Schur complement
preconditioner will converge in one or two iterations and there is
no advantage to using the lower triangular preconditioner.

\subsection{Flow through a nozzle}

In this section we demonstrate the strengths of our method on a test
problem involving steady-state flow through a nozzle in two dimensions.
We compare the steady state flow through the nozzle obtained using
our rigid-body IB method to the flow obtained by using a splitting-based
direct forcing approach \citet{DirectForcing_Uhlmann,RigidIBAMR}.
Specifically, we contrast our monolithic fluid-solid solver to a split
solver based on performing the following operations at time step $n$:Solve
the fluid sub-problem as if the body were not present,
\[
\left[\begin{array}{cc}
\A & \Grad\\
-\Div & \M 0
\end{array}\right]\left[\begin{array}{c}
\tilde{\v}^{n+1}\\
\tilde{\p}^{n+\half}
\end{array}\right]=\left[\begin{array}{c}
\V g^{n+\half}\\
\V 0
\end{array}\right].
\]

\begin{enumerate}
\item Calculate the slip velocity on the set of markers, $\D{\U}=-\left(\J^{n+\frac{1}{2}}\tilde{\v}^{n+1}+\W^{n+\frac{1}{2}}\right)$,
giving the fluid-solid force estimate $\Lamb^{n+\frac{1}{2}}=\left(\rho/\D t\right)\D{\U}$.
\item Correct the fluid velocity to approximately enforce the no-slip condition,
\[
\V v^{n+1}=\tilde{\v}^{n+1}+S^{n+\frac{1}{2}}\D{\U}.
\]

\end{enumerate}
Note that in the original method of \citet{RigidIBAMR} in the last
step the fluid velocity is projected onto the space of divergence-free
vector fields by re-solving the fluid problem with the approximation
$\A\approx\left(\rho/\D t\right)\M I$ (i.e., ignoring viscosity).
We simplify this step here because we have found the projection to
make a small difference in practice for steady state flows, since
the same projection is carried out in the subsequent time step.

\begin{figure}[h]
\begin{centering}
\includegraphics[width=0.9\textwidth]{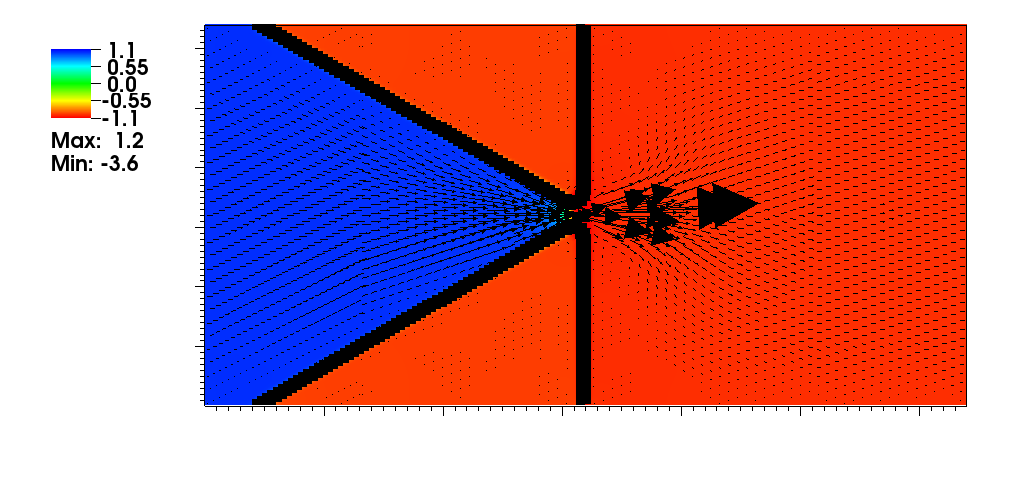}
\par\end{centering}

\begin{centering}
\includegraphics[width=0.49\textwidth]{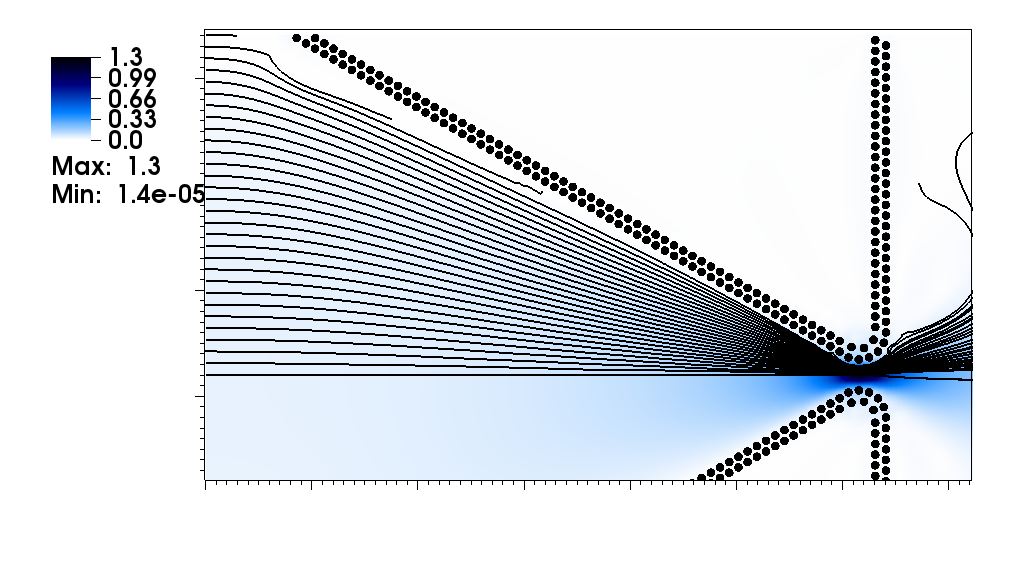}\includegraphics[width=0.49\textwidth]{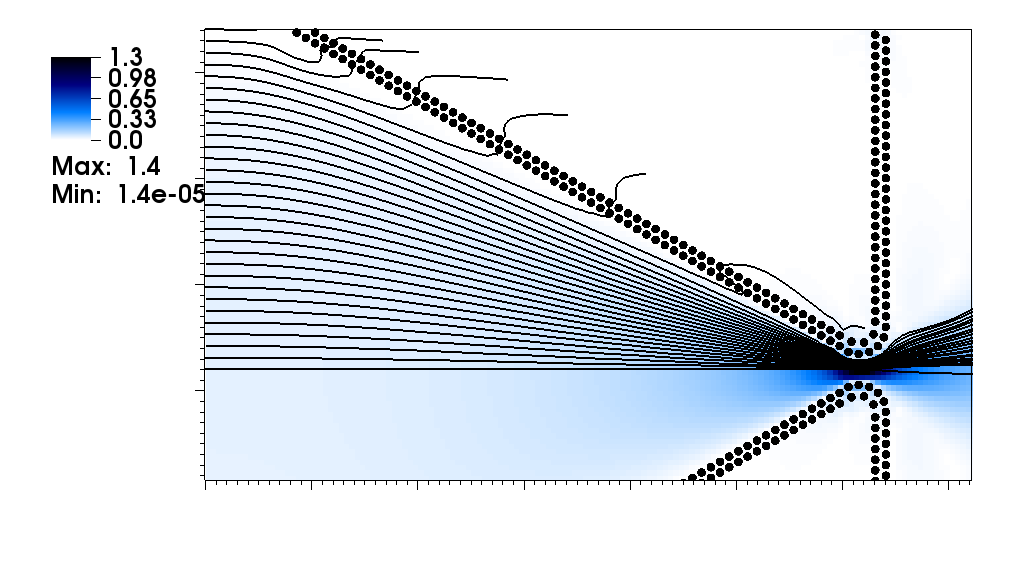}
\par\end{centering}

\centering{}\includegraphics[width=0.49\textwidth]{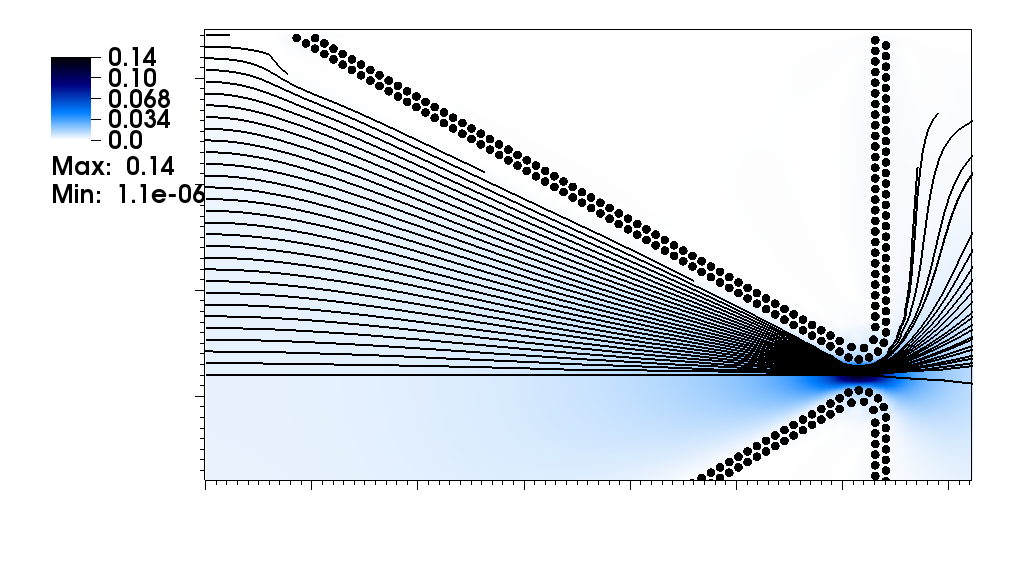}\includegraphics[width=0.49\textwidth]{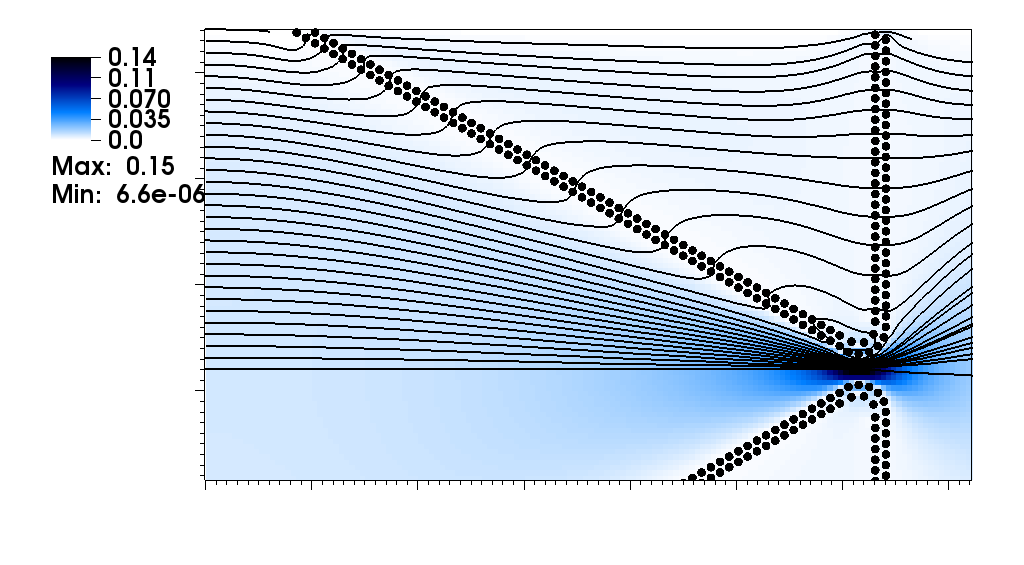}
\caption{\label{fig:nozzle}Two-dimensional flow through a nozzle (top panel)
in a slit channel computed by our rigid IB method (left panels) and
a simplified version of the splitting-based method of Bhalla \emph{et
al.} \citet{RigidIBAMR} (right panels). (Top panel) The geometry
of the channel along with the (approximately) steady state flow at
$\Re\mbox{\ensuremath{\approx}}19$, as obtained using our method.
The color plot shows the pressure and the velocity is shown as a vector
field. (Middle panels) Flow at $\Re\approx19$, computed at time $T=10^{2}$.
For our method (left) we use a time step size of $\D t=5\cdot10^{-2}$
(corresponding to advective CFL number of $U_{max}\D t/\D x\approx0.13$),
while for the splitting method (right) we use $\D t=10^{-3}$. The
streamlines are traced from the entrance to the channel for a time
of $T_{s}=7\cdot10^{3}$ and shown as black lines. (Bottom panels)
Same as the middle row but now for $\Re\approx0.2$, final time $T=10$
and streamlines followed up to $T_{s}=4\cdot10^{4}$, with $\D t=0.125$
for our method (left), and $\D t=10^{-3}$ for the splitting method
(right).}
\end{figure}

We discretize a nozzle constriction in a slit channel using IB marker
points about 2 grid spacings apart. The geometry of the problem is
illustrated in the top panel of Fig. \ref{fig:nozzle}; parameters
are $\rho=1$, grid spacing $\D x=0.5$, nozzle length $l=55.5$,
nozzle opening width $d\approx2.9$, and $\eta$ variable (other parameters
are given in the figure caption). No slip boundary conditions are
specified on the top and bottom channel walls, and on the side walls
the tangential velocity is set to zero and the normal stress is specified
to give a desired pressure jump across the channel of $\D{\p=2}$.
The domain is discretized using a grid of $256\times128$ cells and
the problem evolved for some time until the flow becomes essentially
steady. The Reynolds number is estimated based on the maximum velocity
through the nozzle opening and the width of the opening.

In the bottom four panels in Fig. \ref{fig:nozzle} we compare the
flow computed using our method (left panels) to that obtained using
the splitting-based direct forcing algorithm summarized above (right
panels). Our method is considerably slower (by at least an order of
magnitude) for this specific example because the GMRES convergence
is slow for this challenging choice of boundary conditions at small
Reynolds numbers in two dimensional (recall that steady Stokes flow
in two dimensions has a diverging Green's function). To make the comparison
fairer, we use a considerably smaller time step size for the splitting
method, so that we approximately matched the total execution time
between the two methods. Note that for steady-state problems like
this one with fixed boundaries, it is much more efficient to precompute
the \emph{actual} mobility matrix (Schur complement) once at the beginning,
instead of approximating it with our empirical fits. However, for
a more fair and general comparison we instead use our preconditioner
to solve the constrained fluid problem in each time step anew to a
tight GMRES tolerance of $10^{-9}$. For this test we use $N_{s}=2$
iterations in the fluid solves inside our preconditioner.

The visual results in the right panels of Fig. \ref{fig:nozzle} clearly
show that the splitting errors in the enforcement of the no-slip boundary
condition lead to a notable ``leak'' through the boundary, especially
at small Reynolds numbers. To quantify the amount of leak we compute
the ratio of the total flow through the opening of the nozzle to the
total inflow; if there is no leak this ratio should be unity. Indeed,
this ratio is larger than $0.99$ for our method at all Reynolds numbers,
as seen in the lack of penetration of the flow inside the body in
the left panels in Fig. \ref{fig:nozzle}. For $\Re\approx19$, we
find that even after reducing the time step by a factor of 50, the
splitting method gives a ratio of $0.935$ (i.e., 6.5\% leak), which
can be seen as a mild penetration of the flow into the body in the
middle right panel in Fig. \ref{fig:nozzle}. For $\Re\approx0.2$,
we find that we need to reduce $\D t$ by a factor of 1250 to get
a flow ratio of $0.94$ for the splitting method; for a time step
reduced by a factor of 125 there is a strong penetration of the flow
through the nozzle, as seen in the bottom right panel of Fig. \ref{fig:nozzle}.

\subsection{\label{sub:SphereInShell}Stokes flow between two concentric shells}

Steady Stokes flow around a fixed sphere of radius $R_{1}$ in an
unbounded domain (with fluid at rest at infinity) is one of the fundamental
problems in fluid mechanics, and analytical solutions are well known.
Our numerical method uses a regular grid for the fluid solver, however,
and thus requires a finite truncation of the domain. Inspired by the
work of Balboa-Usabiaga et al. \citet{DirectForcing_Balboa}, we enclose
the sphere inside a rigid spherical shell of radius $R_{2}=4R_{1}$.
This naturally provides a truncation of the domain because the flow
exterior to the outer shell does not affect the flow inside the shell.
Analytical solutions remain simple to compute and are given in Appendix
\ref{AppendixConcentric}.

\begin{figure}[tbph]
\begin{centering}
\includegraphics[width=0.49\textwidth]{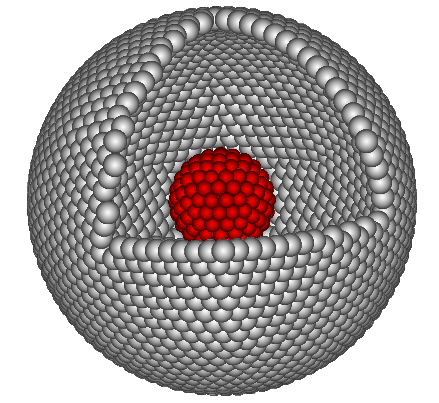}\includegraphics[width=0.49\textwidth]{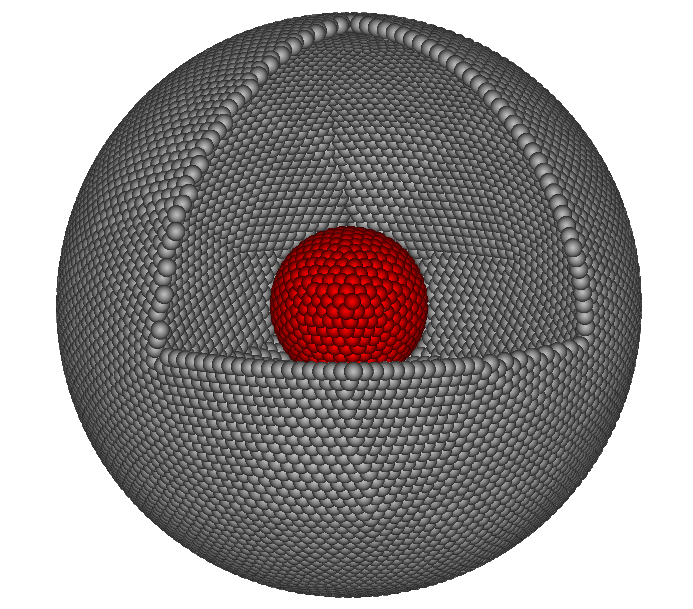}
\par\end{centering}

\caption{\label{fig:Shells}Marker configuration for computing Stokes flow
between two concentric spherical shells. Markers are shown as spheres
with size on the order of their effective hydrodynamic radius. The
inner shell of markers is shown in red, and the outer shell of markers
is shown in gray. (Left) Intermediate resolution, inner shell of 162
markers and outside shell of 2562 markers. (Right) Highest resolution
studied here, inner shell of 642 markers and outside shell of 10242
markers.}
\end{figure}

We discretize the inner sphere using a spherical shell of markers,
since for steady Stokes flow imposing a rigid body motion on the surface
of the inner sphere guarantees a stress-free rigid body motion for
the fluid filling the inner sphere \citet{RegularizedStokeslets}.
We use the same recursive triangulation of the sphere, described in
Section \ref{sec:Solver}, to construct the marker grid for both the
inner and outer shells. The ratio of the number of markers on the
outer and inner spheres is approximately $16$ (i.e., there are two
levels of recursive refinement between the inner and outer shells),
consistent with keeping the marker spacing similar for the two shells
and a fixed ratio $R_{2}/R_{1}=4$. The fluid grid size is set to
keep the markers about two grid cells apart, $s\approx2h$. The rigid-body
velocity is set to $\U=\left(1,0,0\right)$ for all markers on the
outer shell, and to $\U=0$ on all markers on the inner shell. The
outer sphere is placed in a cubic box of length $l=4.15R_{2}$  with
specified velocity $\v=\left(1,0,0\right)$ on all of the boundaries;
this choice ensures that the flow outside of the outer shell is nearly
uniform and equal to $\v=\left(1,0,0\right)$. In the continuum setting,
this exterior flow does not affect the flow of interest (which is
the flow in-between the two shells), but this is not the case for
the IB discretization since the regularized delta function extends
a few grid cells on \emph{both} sides of the spherical shell.

A spherical shell of \emph{geometric }radius $R_{g}$ covered by markers
acts hydrodynamically as a rigid sphere of effective \emph{hydrodynamic}
radius $R_{h}\approx R_{g}+a$ \citet{MultiblobSprings}, where $a$
is the hydrodynamic radius of a single marker \citet{IBM_Sphere,ISIBM,BrownianBlobs}
(we recall that for the six-point kernel used here, $a\approx1.47\, h$).
A similar effect appears in the Lattice-Boltzmann simulations of Ladd,
with $a$ being related to the lattice spacing \citet{VACF_Ladd,MultiblobSprings,FiniteRe_3D_Ladd}.
When comparing to theoretical expressions, we use the effective hydrodynamic
radii of the spherical shells (computed as explained below) and \emph{not}
the geometric radii. Of course, the enhancement of the effective hydrodynamic
radius over the geometric one is a numerical discretization artifact,
and one could choose not to correct the geometric radius. However,
this comparison makes immersed boundary models of steady Stokes flow
appear much less accurate than they actually are in practice. For
example, one should not treat a line of markers as a zero-thickness
object of zero geometric radius; rather, such a line of rigidly-connected
markers should be considered to model a rigid cylinder with finite
thickness proportional to $a$ \citet{IBM_Sphere}.

We can measure the effective hydrodynamic radius $R_{h}$ of a spherical
shell of markers from the drag force on a periodic cubic lattice of
such objects moving with velocity $V$. Specifically, we place a single
shell of $N$ markers in a triply-periodic domain with cubic unit
cell of length $l$, set $\U=\left(1,0,0\right)$ on all markers,
solve (\ref{eq:constrained_Stokes}), and measure the total drag force
as $F=\sum_{i=1}^{N}\Lamb_{i}$. The periodic correction to the Stokes
drag formula is well-known \citet{Mobility2D_Hasimoto},
\begin{equation}
\frac{F}{\eta V}=\frac{6\pi R_{h}}{1-2.8373(R_{h}/l)+4.19(R_{h}/l)^{3}-27.4(R_{h}/l)^{6}+\text{h.o.t.}},\label{eq:drag_periodic}
\end{equation}
and allows us to obtain a very accurate estimate of $R_{h}$ from
the drag for $l\gg R_{h}$. The results are given in the left half
of table \ref{tab:R_h_3D} in the form of the dimensionless ratio
$R_{h}/R_{g}$; we see that as the resolution is increased $R_{h}\rightarrow R_{g}$
with an approximately linear rate of convergence, as expected. Since
this computation refers to flow outside of the shell of markers, we
can call the computed $R_{h}$ the effective \emph{outer} hydrodynamic
radius and use it to set the value of $R_{1}$ in the theory. We use
a similar procedure to measure an effective \emph{inner} hydrodynamic
radius $R_{2}$ for the outer spherical shell. Specifically, we obtain
$R_{2}$ from the drag on the inner sphere based on the theoretical
formula (\ref{eq:drag_inner}), where we use the previously-determined
value of $R_{1}$ for the effective radius of the inner sphere. The
results are given in the right half of Table \ref{tab:R_h_3D} and
again show that as the grid is refined the hydrodynamic radii converge
to the geometric ones.

\begin{table}[tbph]
\caption{\label{tab:R_h_3D}Ratio of the effective hydrodynamic and geometric
radii of the inner (left half) and outer (right half) spherical shells
for simulations of steady Stokes flow around a fixed sphere embedded
within a moving spherical cavity, at different resolutions. }

\centering{}%
\begin{tabular}{|l|l|c|c|c|c|}
\hline 
Resolution  &
Number markers  &
$R_{h}/R_{g}$  &
 &
Number markers  &
 $R_{h}/R_{g}$ \tabularnewline
\hline 
grid size &
Inner shell  &
 &
 &
Outer shell &
\tabularnewline
\hline 
$30^{3}$  &
12 &
1.48  &
 &
162  &
0.93 \tabularnewline
\hline 
$60^{3}$  &
42  &
1.22  &
 &
642  &
0.96 \tabularnewline
\hline 
$120^{3}$  &
162  &
1.09  &
 &
2562  &
0.98 \tabularnewline
\hline 
$240^{3}$  &
642  &
1.04  &
 &
10242  &
0.99 \tabularnewline
\hline 
\end{tabular} 
\end{table}

\subsubsection{Convergence of fluid flow (pressure and velocity)}

\begin{figure}[h]
\begin{centering}
\includegraphics[width=0.98\textwidth]{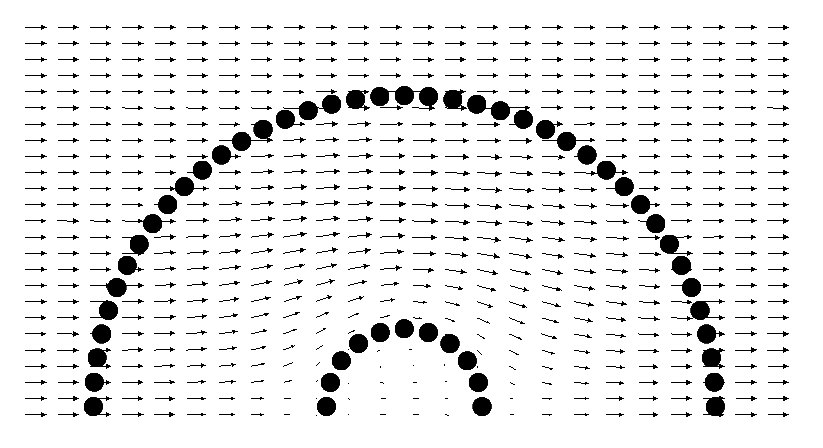}
\par\end{centering}

\begin{centering}
\includegraphics[width=0.49\textwidth]{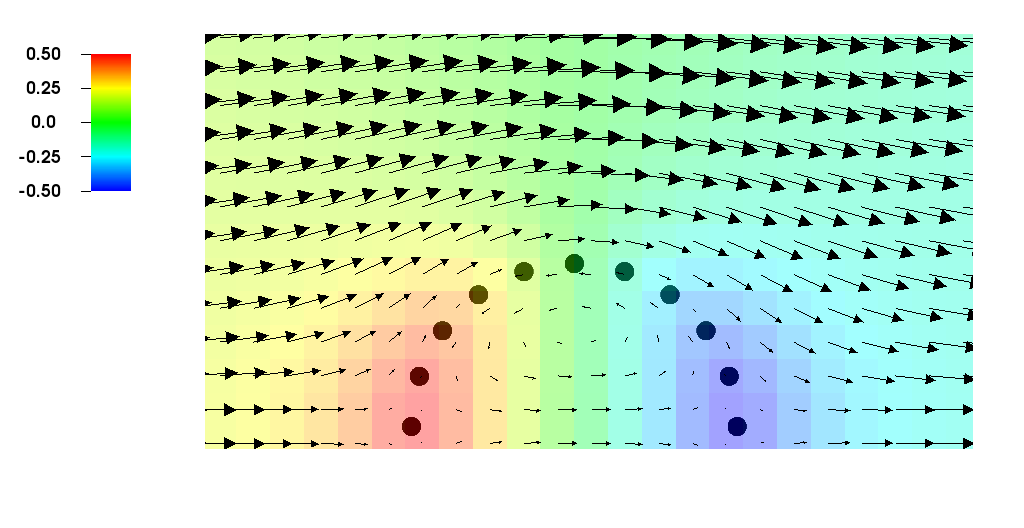}\includegraphics[width=0.49\textwidth]{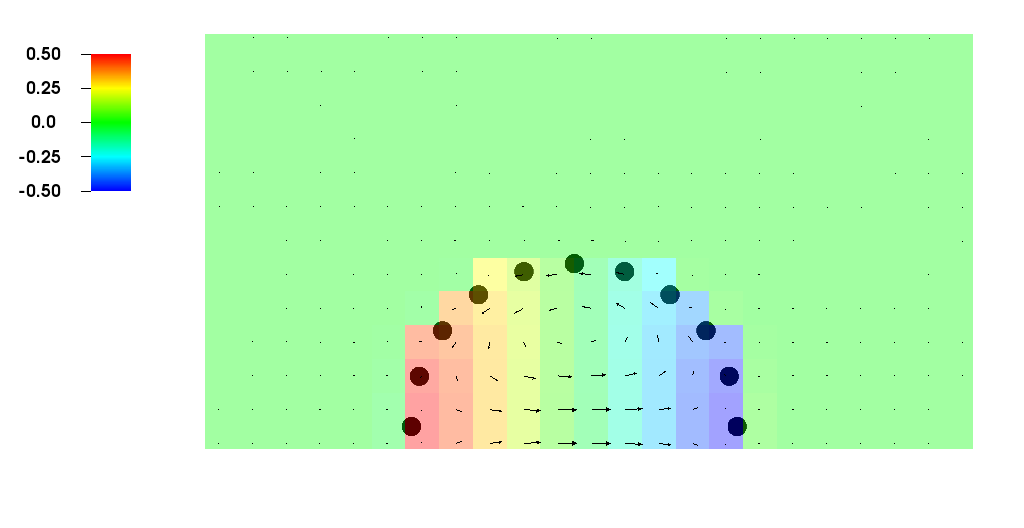}
\par\end{centering}

\caption{\label{fig:ShellFlow}Flow field around a fixed sphere inside a moving
spherical cavity. The outer shell is discretized using 2562 markers,
while the inner one has 162 markers, shown as black circles. (Top)
Velocity field. (Bottom Left) Zoom of the velocity (vector field)
and pressure (color plot) around the inner shell. (Bottom Right) Same
as bottom left but now showing the \emph{error} in the velocity and
pressure compared to the theoretical expressions.}
\end{figure}

The top panel of Fig. \ref{fig:ShellFlow} shows a slice through the
middle of the nested spherical shells along with the fluid velocity
$\v$. Recall that the flow inside the inner sphere should vanish,
implying that the pressure inside the inner shell should be constant
(set to zero here), and the flow outside of the outer sphere should
be uniform. The bottom right panel of the figure zooms in around the
inner sphere to reveal that there is some spurious pressure gradient
and an associated counter-rotating vortex flow generated inside the
inner sphere. The bottom right panel shows the error in the computed
fluid flow $(\v,\p)$, that is, the difference between the computed
flow and the theoretical solution given in Appendix \ref{AppendixConcentric}.
It is clear that the majority of the error is localized in the vicinity
of the inner shell and in the interior of the inner sphere. Note that
these errors would be much larger if the theory had used the geometric
radii instead of the hydrodynamic radii for the shells.

Tables \ref{tab:v_errors} and \ref{tab:p_errors} show the norms
of the error in the computed flow field as a function of resolution.
Asymptotically first-order convergence is observed in the $L_{1}$
and $L_{2}$ norms for both the velocity and the pressure. In the
$L_{\infty}$ norm, we expect the velocity to also converge linearly,
but we do not expect to see convergence in the pressure, since the
velocity is continuous across the interface but the pressure has a
jump; this is consistent with the numerical data.

\begin{table}[h]
\caption{\label{tab:v_errors}Normalized norms of the error in the computed
velocity ($\Delta\v$) for steady Stokes flow around a fixed sphere
embedded within a moving spherical cavity, at different resolutions
(see two left-most columns). An estimated order of convergence based
on successive refinements is indicated in the column to the right
of the corresponding error norm.}

\centering{}%
\begin{tabular}{|l|c|c|c|c|c|c|c|}
\hline 
Markers  &
Resolution  &
$||\Delta\v||_{1}/||\v||_{1}$ &
Rate &
$||\Delta\v||_{2}/||\v||_{2}$ &
Rate &
$||\Delta\v||_{\infty}/||\v||_{\infty}$ &
Rate\tabularnewline
\hline 
162-12  &
$30^{3}$  &
$4.08\cdot10^{-2}$  &
 &
$6.39\cdot10^{-2}$  &
 &
0.558  &
\tabularnewline
\hline 
642-42  &
$60^{3}$  &
$1.14\cdot10^{-2}$  &
1.83 &
$2.08\cdot10^{-2}$ &
1.62 &
0.322  &
0.79\tabularnewline
\hline 
2562-162  &
$120^{3}$  &
$4.61\cdot10^{-3}$  &
1.30 &
$8.74\cdot10^{-3}$  &
1.24 &
0.160  &
1.01\tabularnewline
\hline 
10242-642 &
$240^{3}$  &
$2.16\cdot10^{-3}$  &
1.09 &
$4.26\cdot10^{-3}$  &
1.04 &
0.091 &
0.82\tabularnewline
\hline 
\end{tabular}
\end{table}

\begin{table}[h]
\caption{\label{tab:p_errors}Normalized norms of the error in the pressure
($\Delta\pi$) for steady Stokes flow around a fixed sphere embedded
within a moving spherical cavity, at different resolutions (see two
left-most columns). An estimated order of convergence based on successive
refinements is indicated in the column to the right of the corresponding
error norm.\textbf{ }}

\centering{}%
\begin{tabular}{|l|c|c|c|c|c|c|c|}
\hline 
Markers  &
Resolution  &
$||\Delta\pi||_{1}/||\pi||_{1}$ &
Rate &
$||\Delta\pi||_{2}/||\pi||_{2}$ &
Rate &
$||\Delta\pi||_{\infty}/||\pi||_{\infty}$ &
Rate\tabularnewline
\hline 
162-12  &
$30^{3}$  &
$0.849$  &
 &
0.788 &
 &
1.0 &
\tabularnewline
\hline 
642-42  &
$60^{3}$  &
0.567 &
0.58 &
0.486  &
0.70 &
1.0  &
0.0\tabularnewline
\hline 
2562-162  &
$120^{3}$  &
$0.344$  &
0.72 &
0.275  &
0.82 &
0.860 &
0.22\tabularnewline
\hline 
10242-642 &
$240^{3}$  &
$0.196$  &
0.81 &
0.164  &
0.75 &
0.704  &
0.29\tabularnewline
\hline 
\end{tabular}
\end{table}

\subsubsection{Convergence of Lagrangian forces (surface stresses)}

The first-order convergence of the pressure and velocity is expected
and well-known in the immersed boundary community. The convergence
of the \emph{tractions} $\left(\M{\sigma}\cdot\V n\right)$ on the
fluid-body interface is much less well studied, however. This is in
part because in penalty-based or splitting methods, it is difficult
to estimate tractions precisely (e.g., for penalty methods using stiff
springs, the spring tensions oscillate with time), and in part because
a large number of other studies have placed the markers too closely
to obtain a well-conditioned mobility matrix and thus to obtain accurate
forces. Furthermore, there are at least two ways to estimate surface
tractions in IB methods, as discussed in extensive detail in Ref.
\citet{Tractions_Fauci}. One method is to estimate fluid stress from
the fluid flow and extrapolate toward the boundary. Another method,
which we use here, is to use the computed surface forces $\Lamb$
to estimate the tractions. 

\begin{figure}[h]
\begin{centering}
\includegraphics[width=0.33\textwidth]{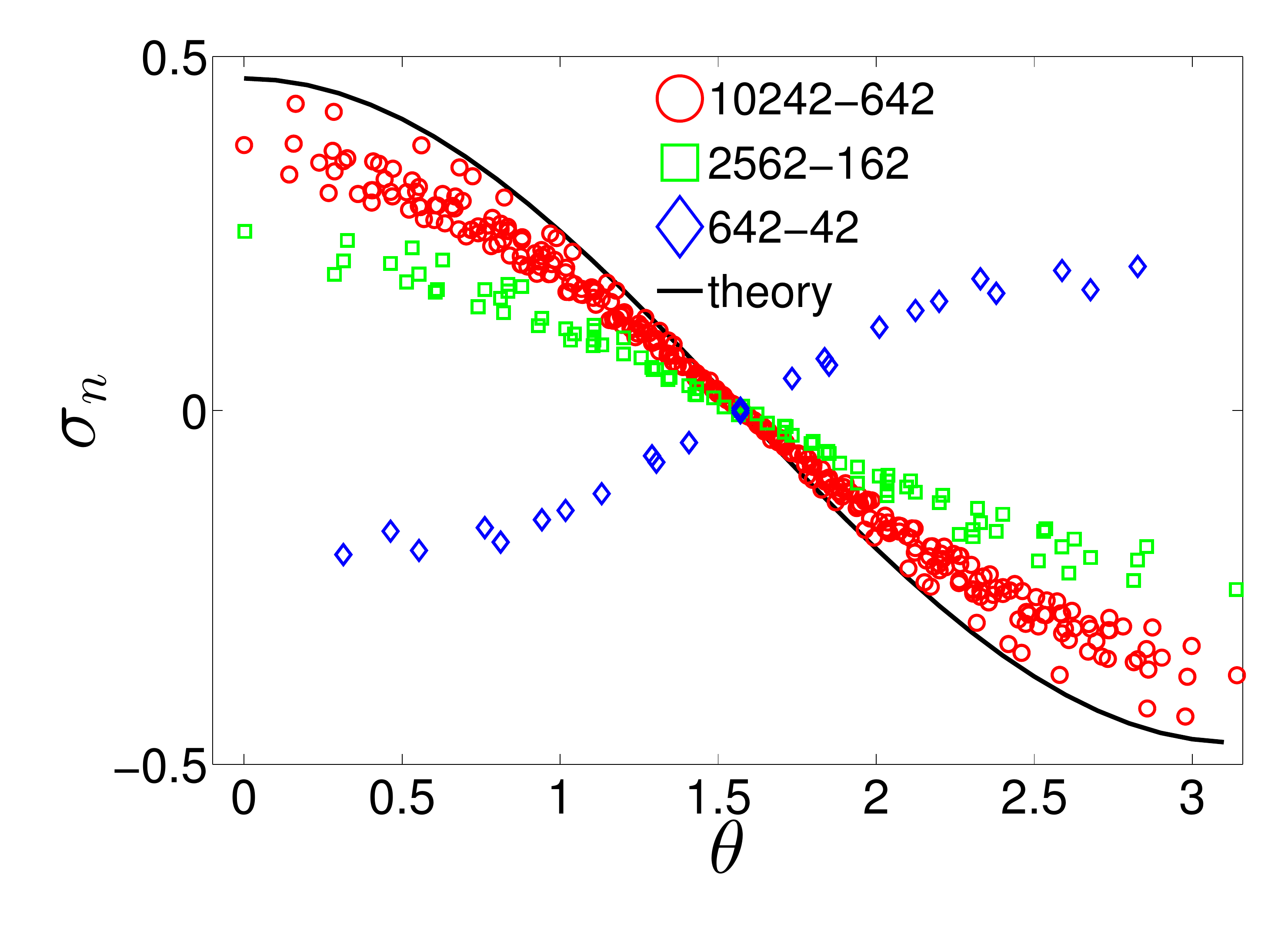}\includegraphics[width=0.33\textwidth]{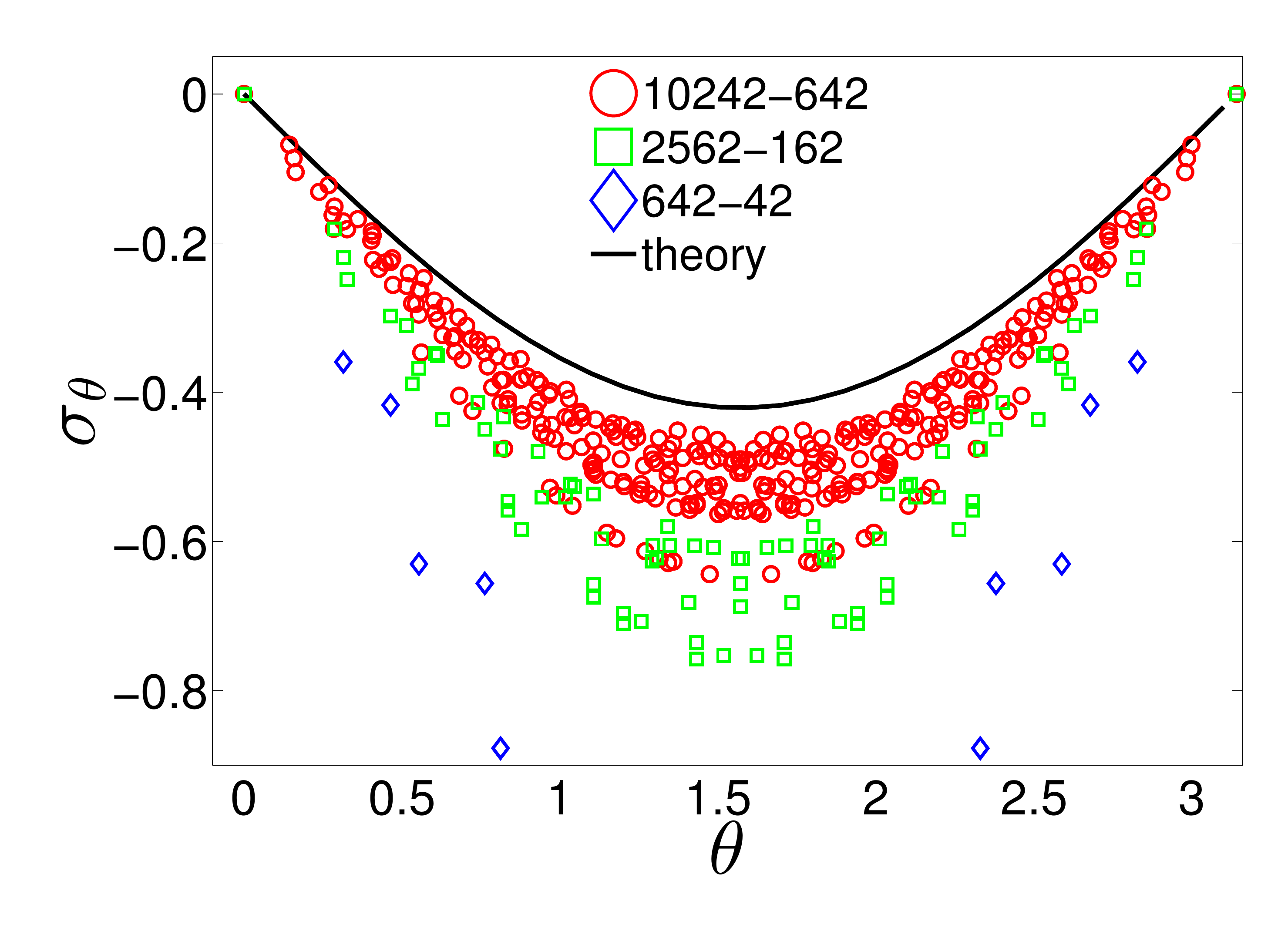}\includegraphics[width=0.33\textwidth]{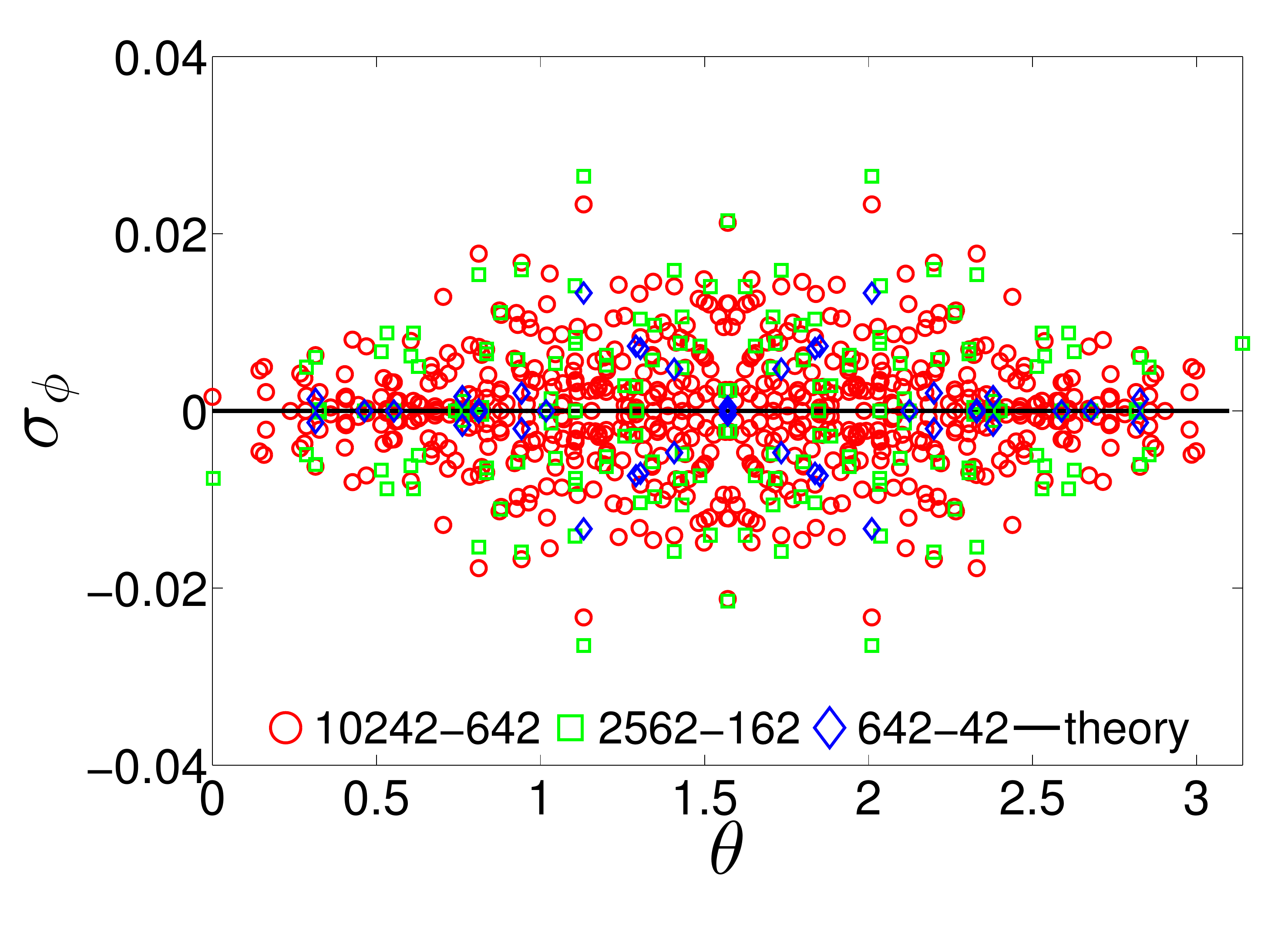}
\par\end{centering}

\begin{centering}
\includegraphics[width=0.33\textwidth]{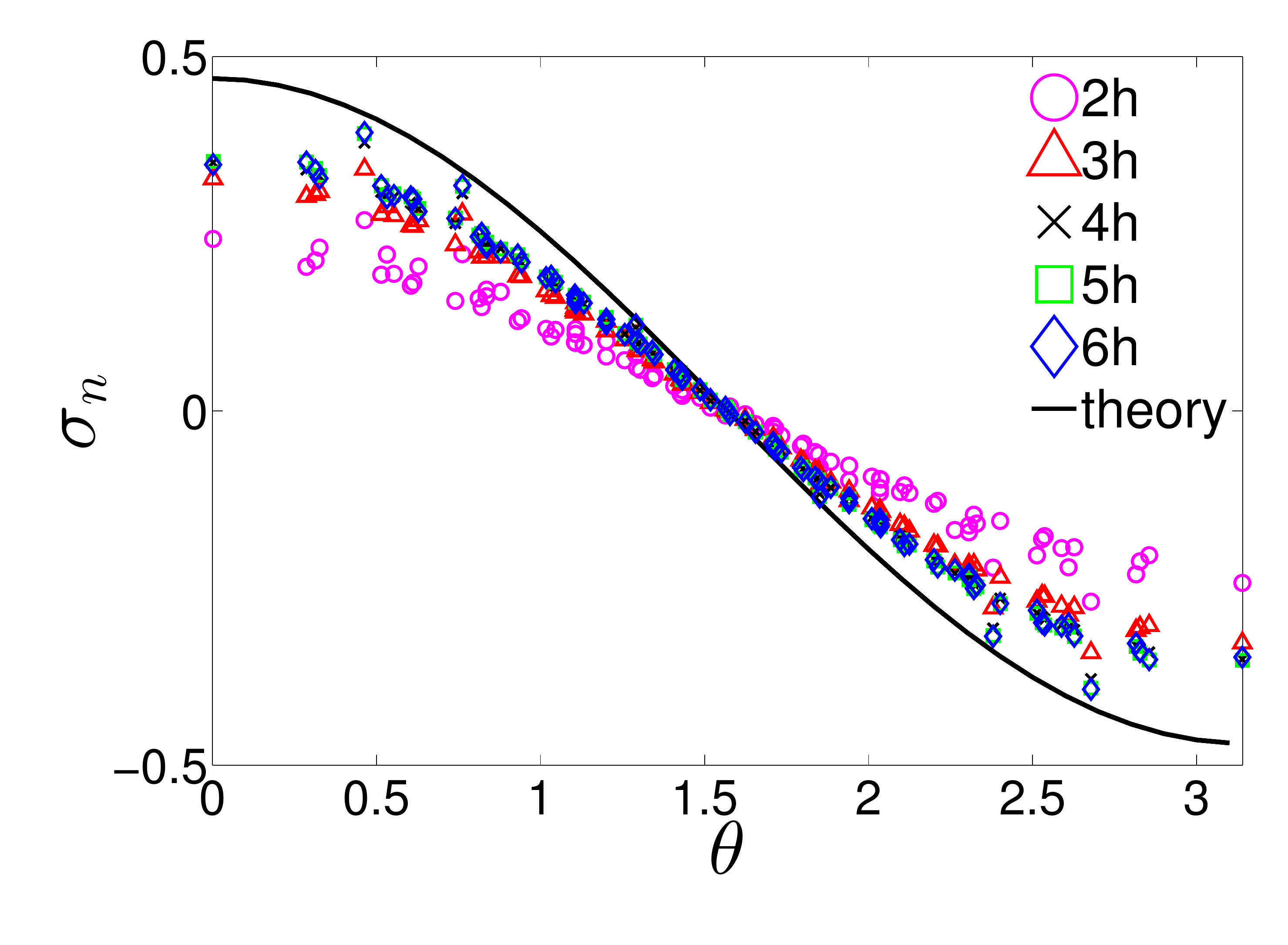}\includegraphics[width=0.33\textwidth]{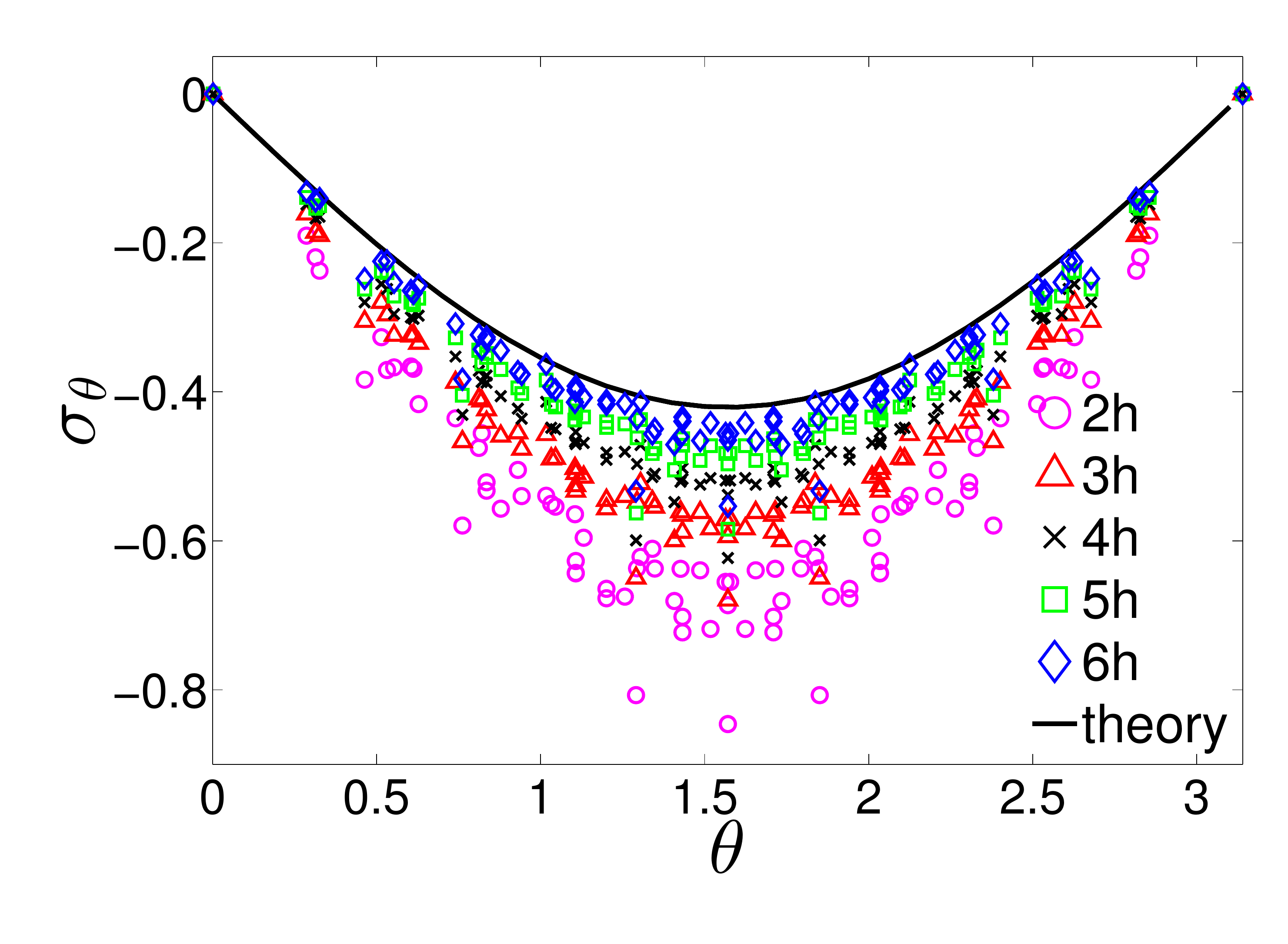}\includegraphics[width=0.33\textwidth]{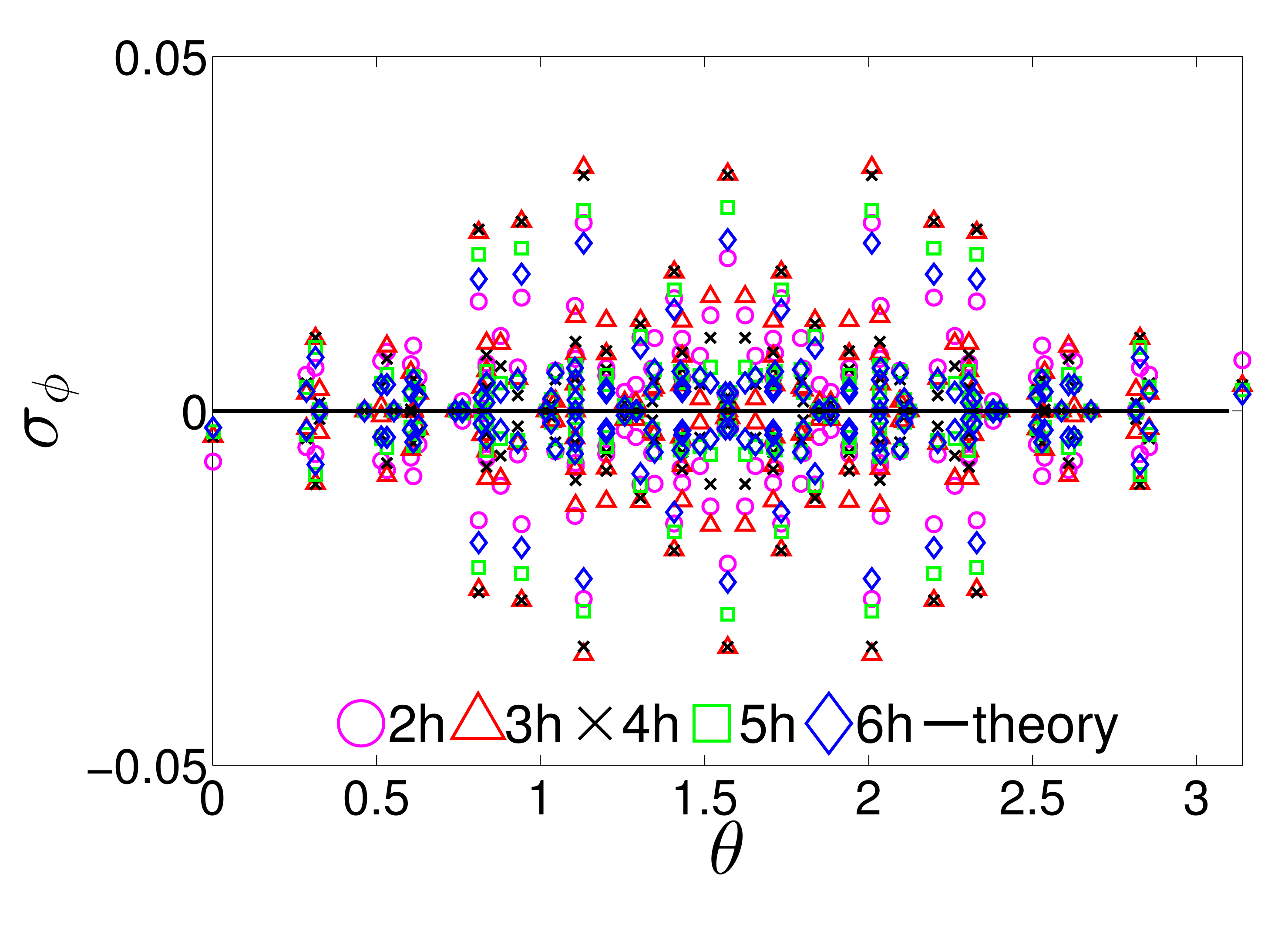}
\par\end{centering}

\caption{\label{fig:ConvTractions}Convergence of surface stresses to their
theoretical values for the three different resolutions. Pointwise
traction estimates are shown with symbols as a function of the angle
$\theta$ relative to the direction of the flow, while the theory
is shown with a solid black line. (\emph{Left columns}) Normal component
of the traction $\sigma_{n}=\hat{\V r}\cdot\M{\sigma}\cdot\V n$.
(\emph{Center columns}) Tangential component of traction in direction
of flow, $\sigma_{\theta}=\hat{\V{\theta}}\cdot\M{\sigma}\cdot\V n$.
(\emph{Right columns}) Tangential component in the direction perpendicular
to the flow, $\sigma_{\phi}=\hat{\V{\phi}}\cdot\M{\sigma}\cdot\V n$,
which should vanish by symmetry. (\emph{Top row}) Different resolutions
(see legend) for a fixed spacing $s\approx2h$. Note that for the
coarsest resolution of only 12 markers on the inner sphere, the computed
tractions have values off the scale of this plot and are thus not
shown. (\emph{Bottom row}) The most resolved case of $2562-162$ markers
for different spacing between markers, as indicated in the legend.
Note that using $s\approx h$ leads to severe ill-conditioning and
the computed tractions show random scatter well beyond the scale of
the plot and are thus not shown.}
\end{figure}

We obtain pointwise estimates of the tractions at the positions of
the markers from the relation $\left(\M{\sigma}\cdot\V n\right)\left(\q_{i}\right)\approx\lamb_{i}/\D A_{i}$,
where $\D A_{i}$ is the surface area associated with marker $i$.
We obtain $\D A_{i}$ from the surface triangulation used to construct
the marker of grids by assigning one third of the area of each triangle
to each of its nodes. In Fig. \ref{fig:ConvTractions} we show the
computed normal and tangential components of the traction in polar
coordinates, with the $z$ symmetry axes along the direction of the
flow. The theoretical prediction given in Appendix \ref{AppendixConcentric}
is shown with a black line and is based on the \emph{geometric} radii. 

The top row of Fig. \ref{fig:ConvTractions} shows the computed tractions
for several resolutions with marker spacing $s\approx2h$. It is seen
that as the grid is refined, the computed tractions appear to converge
\emph{pointwise} to the correct values. However, the convergence is
very slow, and even for the large resolutions reported here, it is
evident that the asymptotic convergence regime has not been reached.
Consequently, no precise statement about the order of convergence
can be made from these data. At lower resolution, some of the results
even show qualitatively wrong behavior. For example, the normal traction
$\sigma_{n}=\V n\cdot\M{\sigma}\cdot\V n$ for a resolution of 42
inner and 642 outer markers grows with $\theta$, but the theoretical
result decreases with $\theta$. We also see scatter in the values
among individual markers, indicating that the geometrical and topological
non-uniformity of the marker grid affect the pointwise values.

Nonetheless, we remark that low-order \emph{moments} of the surface
tractions are much more accurate than the pointwise tractions. For
example, the total drag on the inner sphere is much more accurate,
as seen in Table \ref{tab:R_h_3D}. Other test problems not reported
here indicate that stresslets are also computed quite accurately,
especially if one accounts for the distinction between geometric and
hydrodynamic radii. These findings suggest that weak convergence of
the tractions is more robust than strong convergence. In fact, lower
order moments can show reasonable behavior even if the marker spacing
is small and the pointwise forces are numerically unstable to compute.
Somewhat unsurprisingly, we find that the pointwise traction estimates
are improved as the spacing among the markers is increased; see the
bottom row in Fig. \ref{fig:ConvTractions}. The improvement is not
only due to the reduction of the scatter, as expected from the improvement
in conditioning number of the mobility matrix, but also due to global
reduction of the error in the tractions; the observed global reduction
may, however, be specific to steady Stokes flow. For widely spaced
markers, however, the error in the computed flow field field will
increase because the flow will penetrate the shell boundary. This
once again demonstrates the delicate balance that is required in choosing
the marker spacing for rigid bodies, as we discuss further in the
Conclusions.

\subsection{Steady Stokes flow around sphere in a slit channel}

In this section we study a problem at zero Reynolds numbers with nontrivial
boundary conditions, namely, steady Stokes flow around a sphere in
a slit channel (flow between two parallel walls). It is well-known
that computing flows in such geometries using Green's function based
methods such as boundary-integral methods is highly-nontrivial \citet{BD_LB_Ladd,BoundaryIntegral_Wall,StokesianDynamics_Slit}.
Specific methods for spheres in a channel have been developed \citet{HE_Spheres_TwoWalls}
but these are not general, in particular, flow in a square channel
requires a different method, and incorporating the periodicity in
some of the dimensions is nontrivial \citet{BoundaryIntegral_Periodic3D,BoundaryIntegral_Wall}.
At the same time, we wish to point out that boundary-integral methods
have some advantages over our IB method as well. Notably, they are
\emph{considerably} more accurate, and handling domains unbounded
in one or more directions is possible by using the appropriately-decaying
Green's function.

Unlike the case of a single no-slip boundary, writing down an analytical
solution for slit channels is complex and requires numerically-evaluating
the coefficients in certain series expansions \citet{StokesianDynamics_Slit}.
For the component of the mobility $\mu=F/V$ of a sphere in an infinite
slit channel, Faxen has obtained exact series expansions for the mobility
at the half and quarter channel locations,
\begin{eqnarray}
\mu_{\parallel}\left(H=\frac{d}{2}\right) & = & \frac{1}{6\pi\eta R_{h}}\Bigg[1-1.004\frac{R_{h}}{H}+0.418\frac{R_{h}^{3}}{H^{3}}+0.21\frac{R_{h}^{4}}{H^{4}}-0.169\frac{R_{h}^{5}}{H^{5}}+\dots\Biggr]\nonumber \\
\mu_{\parallel}\left(H=\frac{d}{4}\right) & = & \frac{1}{6\pi\eta R_{h}}\Biggl[1-0.6526\frac{R_{h}}{H}+0.1475\frac{R_{h}^{3}}{H^{3}}-0.131\frac{R_{h}^{4}}{H^{4}}-0.0644\frac{R_{h}^{5}}{H^{5}}+\dots\Biggr]\label{eq:Faxen_slit}
\end{eqnarray}
where $R_{h}$ is the (hydrodynamic) radius of the sphere, $H$ is
the distance from the center of the particle to the nearest wall,
and $d$ is the distance between the walls.

\begin{figure}[tbph]
\begin{centering}
\includegraphics[width=0.75\textwidth]{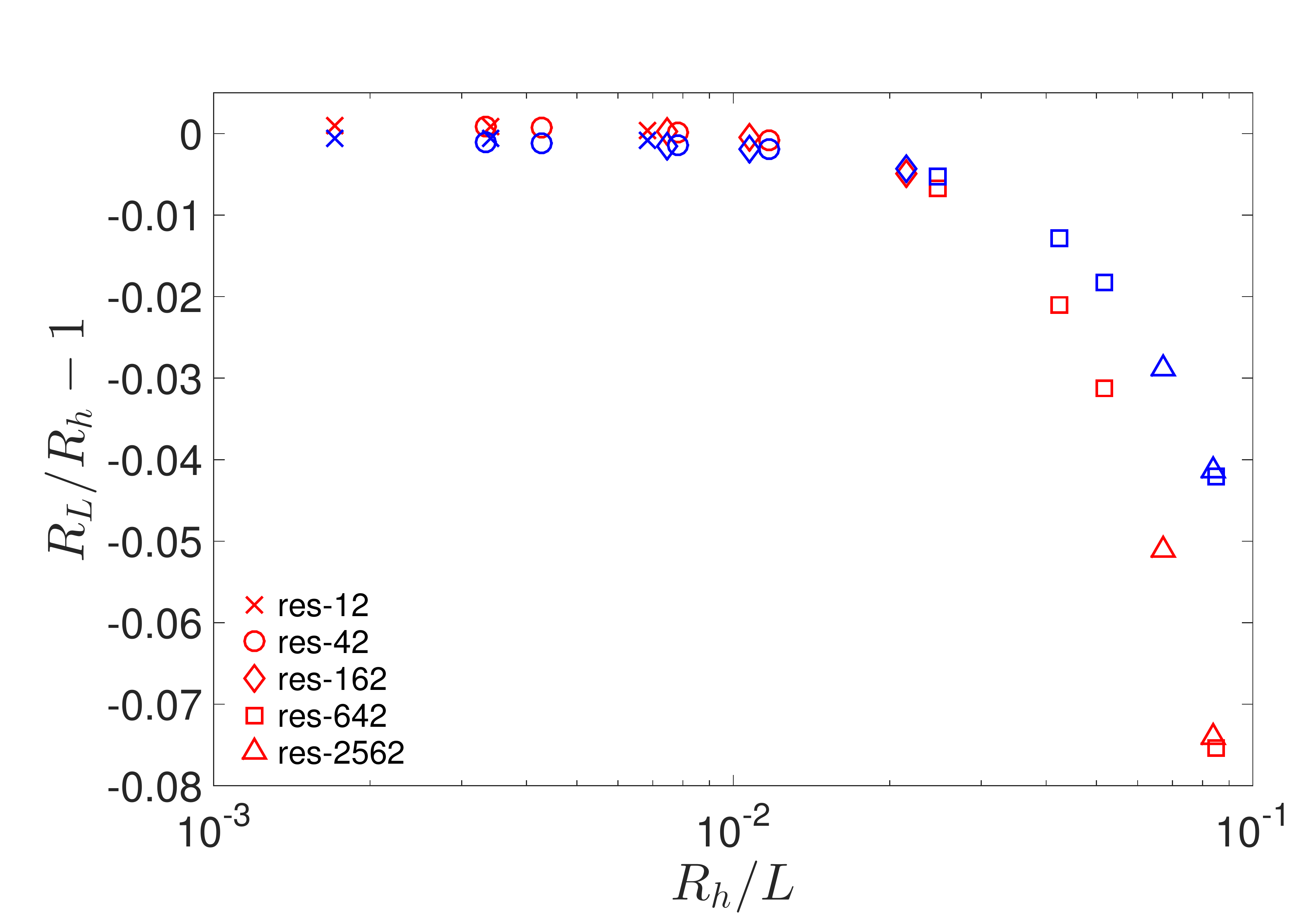}
\par\end{centering}

\caption{\label{fig:Slit-study}Effective hydrodynamic radius $R_{L}(L)$ of
a sphere of hydrodynamic radius $R_{h}$, translating parallel to
the walls of a slit channel of dimensions $L\times L\times d$. Red
symbols are for the sphere at the midplane of the channel, $H/d=0.5$,
and blue symbols are for $H/d=0.25$; these two give different dependence
$R_{L}(L)$ as expected. The channel width is taken $d\approx8R_{H}$
and different numbers of markers are used for the sphere (see legend),
and the grid spacing is set to give a marker spacing $s$ as close
as possible to $a/s\approx0.5$. Note that similar to the example
of flow between concentric spheres the correct value of the drag is
determined by the larger hydrodynamic and \emph{not} by the geometric
radius of the shell.}
\end{figure}

To simulate a spherical particle in a slit channel we place a single
spherical shell with different number of IB markers in a domain of
size $L\times L\times d$, at either a quarter or half distance from
the channel wall. No slip walls are placed at $z=0$ and $z=L$, and
periodic boundary conditions are applied in the $x$ and $y$ directions.
For each $L$, we compute an effective hydrodynamic radius $R_{L}$
by assuming (\ref{eq:Faxen_slit}) holds with $R_{h}$ replaced by
$R_{L}$. We know that as $L\rightarrow\infty$ we have $R_{L}\rightarrow R_{h}$,
however, we are not aware of theoretical results for the dependence
$R_{L}(L)$ at finite $L$. In Fig. \ref{fig:Slit-study}, we plot
$R_{L}/R_{h}-1$ versus $R_{h}/L$ for $d\approx8R_{h}$. Here the
effective hydrodynamic radius of the shell $R_{h}$ is estimated by
using (\ref{eq:drag_periodic}), as shown in Table \ref{tab:R_h_3D}
(see inner radius). We see that we have consistent data for $R_{L}(L)$
among different resolutions, and we obtain consistency in the limit
$L\rightarrow\infty$. This indicates that even a low-resolution model
with as few as 12 markers offers a reasonably-accurate model of a
sphere of effective radius $R_{h}$, \emph{independent} of the boundary
conditions.

\subsection{\label{sub:CylindersSteady}Steady Stokes flow around cylinders}

Here we study the drag force on a periodic square array of cylinders
(i.e., disks in two dimensions) with lattice spacing $l$. The corresponding
study in three dimensions is presented in Section \ref{sub:CubicSpheres}.
The equivalent of (\ref{eq:drag_periodic}) in two dimensions for
\emph{dilute} systems is \citet{Mobility2D_Hasimoto,FiniteRe_2D_Ladd}
\begin{equation}
\frac{F}{\eta V}=\frac{4\pi}{-\ln(\sqrt{\phi})-0.738+\phi-0.887\phi^{2}+2.038\phi^{3}+O\left(\phi^{4}\right)},\label{eq:drag_periodic_2D}
\end{equation}
where $\phi=\pi R_{h}^{2}/l^{2}$ is the packing fraction of the disks
and $R_{h}$ is the hydrodynamic radius of the cylinder, which is
defined from (\ref{eq:drag_periodic_2D}). Observe that in two dimensions,
there is no limit as $\phi\rightarrow0$, in agreement with Stokes's
paradox for flow around a single cylinder; one must account for inertial
effects for very small volume fractions in order to obtain physically-relevant
results. Table \ref{tab:R_h_2D} reports $R_{h}$ for several different
marker models of a cylinder, as estimated by computing the drag for
a range of packing fractions and extrapolating to $\phi\ll1$ using
(\ref{eq:drag_periodic_2D}). As expected, the more resolved the cylinder,
the closer $R_{h}$ is to $R_{g}$. Filling the cylinder with markers
both substantially enlarges the effective hydrodynamic radius and
also degrades the conditioning of the mobility matrix, and is therefore
not advised at zero Reynolds number.

\begin{table}[h]
\caption{\label{tab:R_h_2D}Hydrodynamic radii of several discretizations of
a cylinder with different numbers of markers on the surface and the
interior of the body, keeping $s/h\approx2$. Two models have markers
only on the surface of the cylinder (see the right panel of Fig. \ref{fig:2DWake}
for a 39-marker shell). The rest of the models are constructed from
a regular polar grid of markers filling the interior of the cylinder
(see, for example, left panel of Fig. \ref{fig:2DWake} for a 121-marker
cylinder).}

\centering{}%
\begin{tabular}{|l|c|c|c|}
\hline 
Number markers  &
Surface markers &
Interior markers &
$R_{h}/R_{g}$\tabularnewline
\hline 
39 shell  &
39 &
0 &
1.04 \tabularnewline
\hline 
121 cylinder &
37 &
84 &
1.15 \tabularnewline
\hline 
100 shell &
100 &
0 &
1.02 \tabularnewline
\hline 
834 cylinder &
100 &
734 &
1.03 \tabularnewline
\hline 
\end{tabular}
\end{table}

Another interesting limit for which there are theoretical results
is the \emph{dense} limit, in which the disks/cylinders almost touch,
so that there is a lubrication flow between them. In this limit \citet{FiniteRe_2D_Ladd},
\begin{equation}
\frac{F}{\eta V}\approx\frac{9\pi}{2^{\frac{3}{2}}}\varepsilon^{-\frac{5}{2}},\label{eq:Lubrication2D}
\end{equation}
where $\varepsilon=1-\sqrt{4\phi/\pi}=\left(l-2R_{h}\right)/l$ is
the relative gap between the particles. Note that because the number
of hydrodynamic cells must be an integer, we cannot get an arbitrary
gap between the cylinders for a given cylinder model and fixed $s/h$
(i.e., a fixed $R_{h}/R_{g}$). Also note that when the gap between
the cylinders is too small, the kernels from markers on two cylinders
start to overlap, and the problem becomes ill-conditioned; we have
been able to compute reliable results down to a relative gap of $\varepsilon\gtrsim10^{-2}$
for the resolutions studied here.

\begin{figure}[h]
\centering{}\includegraphics[width=0.49\textwidth]{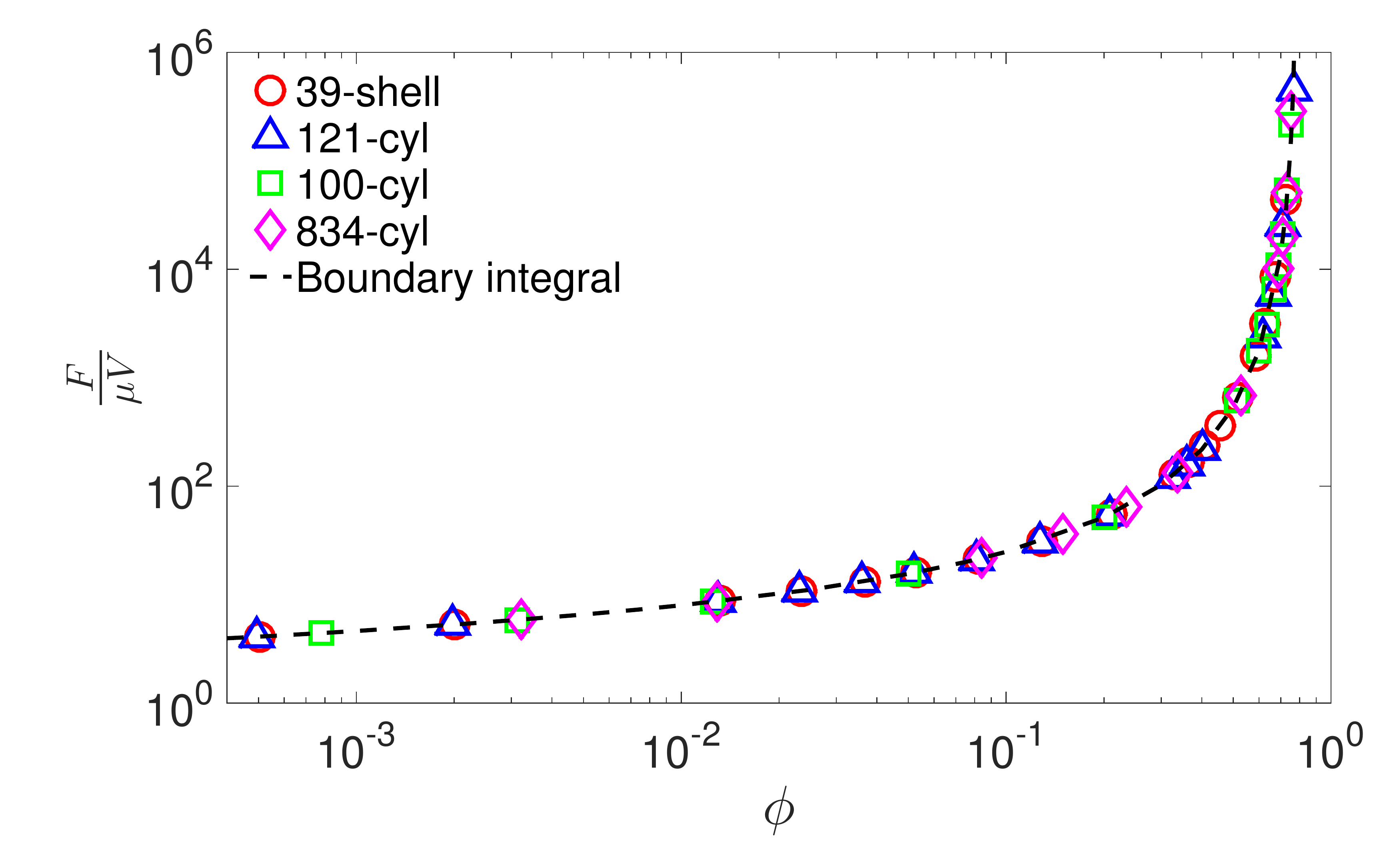}\includegraphics[width=0.49\textwidth]{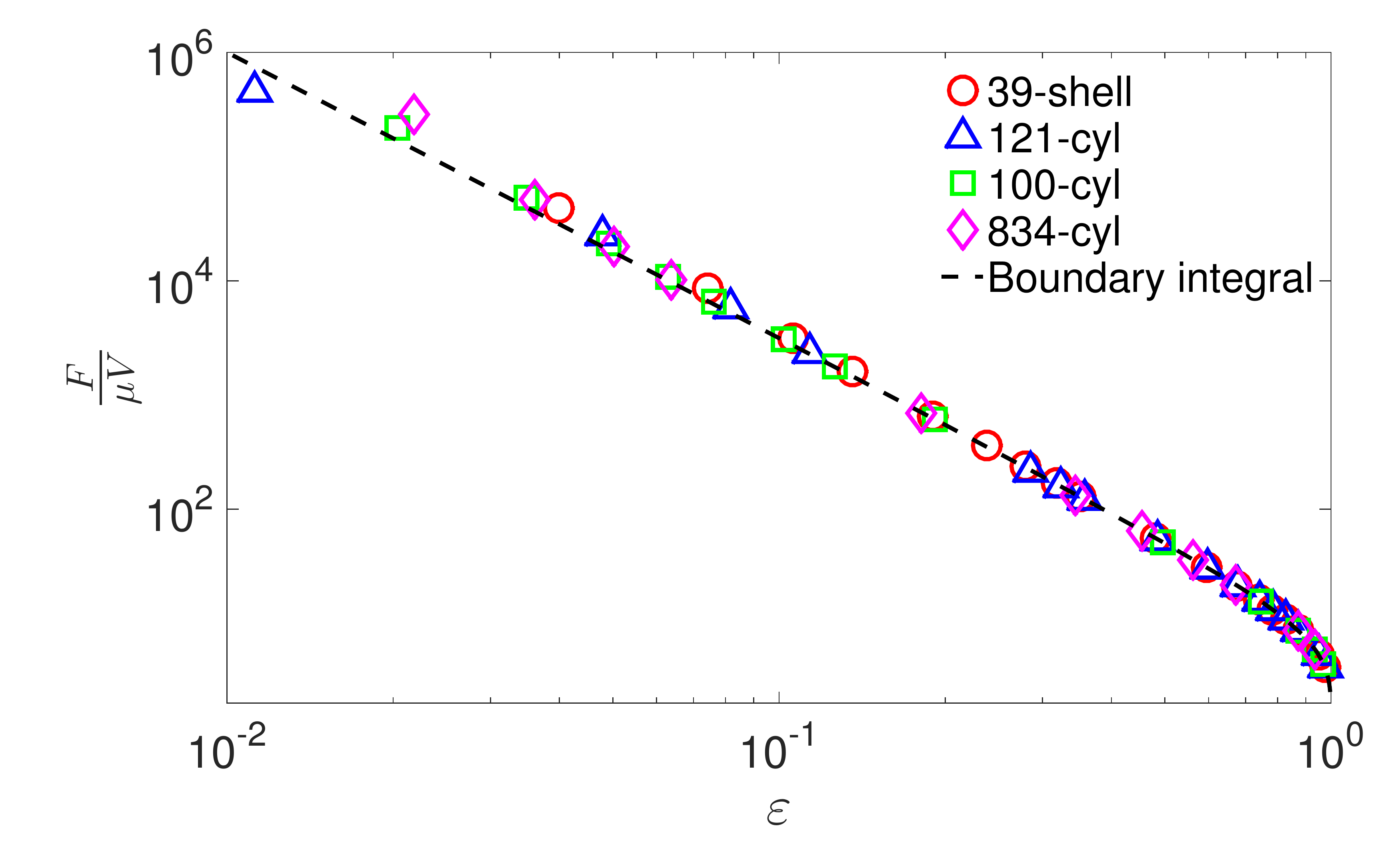}
\caption{\label{fig:Drag2D_phi_Re_0}The drag coefficient for a periodic array
of cylinders in steady Stokes flow for different resolutions (see
Table \ref{tab:R_h_2D}). (Left) As function of volume fraction, compared
to the results of a highly-accurate boundary-integral method. (Right)
Zoom in for close-packed arrays with inter-particle gap plotted on
a log scale to show the asymptotic $\varepsilon^{-\frac{5}{2}}$ divergence
of the lubrication force.}
\end{figure}

Numerical results for the normalized drag over a broad range of volume
fractions are shown in Fig. \ref{fig:Drag2D_phi_Re_0} and compared
to results obtained using an in-house two-dimensional version of the
spectrally-accurate boundary integral method proposed in Ref. \citet{BoundaryIntegral_Periodic3D}.
We obtain very good agreement, similar to that observed using the
Lattice Boltzmann method \citet{FiniteRe_2D_Ladd}, indicating that
even moderately-resolved cylinders are good representations so long
as one uses their hydrodynamic rather than their geometric radius
when computing the effective volume fraction. In particular, in the
right panel, we obtain excellent agreement with the lubrication result
(\ref{eq:Lubrication2D}), seeing an increase in the drag of over
six orders of magnitude consistent with theory. Of course, the IB
method results for the drag do not have a true divergence as $\epsilon\rightarrow0$
because of the regularization of the singular kernel; one must use
adaptively-refined non-regularized boundary integral methods to truly
resolve the divergence. In practice, however, effects not included
in the theoretical model, such as surface roughness or partial slip,
will mollify the unphysical divergence.

\subsection{\label{sub:UnsteadyCylinders}Unsteady flow around cylinders}

Next, we examine the ability of our rigid-body IB method to model
unsteady two-dimensional flow around cylinders (disks). We define
Reynolds number by 
\[
\Re=\frac{\rho VR_{h}}{\eta}=\frac{VR_{h}}{\nu},
\]
where $R_{h}$ is the hydrodynamic radius of the cylinder measured
at $\Re=0$ (see table \ref{tab:R_h_2D}), and $V$ is the velocity
of the incident flow. For small Reynolds numbers, the mean drag per
unit length $F$ is given by \citet{FiniteRe_2D_Ladd} 
\[
\frac{F}{\eta V}=k_{0}+k_{2}\Re^{2},
\]
where $k_{0}\left(\phi\right)$ and $k_{2}\left(\phi\right)$ are
constants that depend on the packing fraction $\phi$ (defined using
the hydrodynamic radius). In the range $\Re\sim2-5$, the drag becomes
quadratic in the flow rate \citet{FiniteRe_2D_Ladd}, and for moderate
Reynolds numbers, a drag coefficient is defined from the empirical
relation
\[
C_{D}=\frac{F}{\rho V^{2}R_{h}}.
\]
As the Reynolds number is increased, the flow becomes unsteady and
vortex shedding occurs, and eventually there is a transition to three-dimensional
flow. Here we focus on steady flow at $\Re\leq100$. 

A staggered-grid variant of the piecewise-parabolic Godunov method
is used for spatial discretization of the advective terms, as explained
in detail by Griffith \citet{NonProjection_Griffith}. In our tests,
the time step size is determined by fixing the advective Courant number
$V\D t/h=0.1$; this value is well bellow the stability limit and
ensures that the discretization errors coming from the (unconstrained)
fluid solver are small. The Adams-Bashforth method is used to handle
advection explicitly. The viscous terms are handled implicitly using
the backward Euler method rather than the implicit midpoint rule because
we are interested in steady states and not transient phenomena. We
initialize the simulations with the fluid moving at a uniform velocity
but allow enough time for a steady-state to be reached.

\subsubsection{Drag on periodic array of disks}

The permeability of a periodic array of aligned cylinders is a well-studied
problem and can be computed by placing a single cylinder in a periodic
domain. To create flow through the periodic system, we follow Ladd
\emph{et al.} \citet{FiniteRe_3D_Ladd,FiniteRe_2D_Ladd} and apply
a constant body force $\f$ throughout the domain (including in the
interior of the body). We solve the constrained time-dependent problem
to a steady state, keeping the cylinder at rest, and measure the average
velocity in the domain, $\bar{\v}=\DV^{-1}\,\int_{V}\v\, d\r.$ In
two spatial dimensions, the dimensionless drag coefficient is defined
by 
\[
k=\frac{F_{x}}{\eta\bar{v}_{x}},
\]
where the force $\F=\DV\,\f=-\V 1^{T}\Lamb$ is the total force applied
to the fluid, which must also equal the negative of the total force
exerted on the rigid body. 

Theory suggests that the correction to the drag scales as $\Re^{2}$
for small Reyonds numbers due to the anti-symmetry of the correction
to the flow (relative to steady Stokes) of order $\Re$ \citet{FiniteRe_2D_Ladd},
so that 
\begin{equation}
k=k_{0}+k_{2}\Re^{2},\label{eq:k_theory}
\end{equation}
where the values $k_{0}(\phi)$ and $k_{2}(\phi)$ depend on the packing
fraction $\phi$. To obtain $k_{0}$, we move the body at a constant
velocity and obtain the drag force $\V 1^{T}\Lamb$ from the solution
of the constrained steady Stokes problem (\ref{eq:constrained_periodic}).
Because marker-based models of rigid bodies do not have perfect symmetry,
the force $\f_{0}=-\DV^{-1}\,\left(\V 1^{T}\Lamb\right)$ has small
nonzero components in the direction perpendicular to the flow. To
ensure that in the limit $\Re\rightarrow0^{+}$ we have perfect consistency
between the finite $\Re$ and zero $\Re$ computations, we use the
force $\V f=\left(k/k_{0}\right)\V f_{0}$ to drive the flow at finite
$\Re$ numbers. Note that it can take thousands of time steps for
the steady state to be established for $\Re\gtrsim1$; to accelerate
convergence, we initialize the computation for a given $\Re$ from
the steady state for the closest smaller $\Re$. Also note that the
exact mobility matrix $\Mob$ and its factorization can be precomputed
once at the beginning and used repeatedly for these steady-state calculations.

\begin{figure}[h]
\centering{}\includegraphics[width=0.85\textwidth]{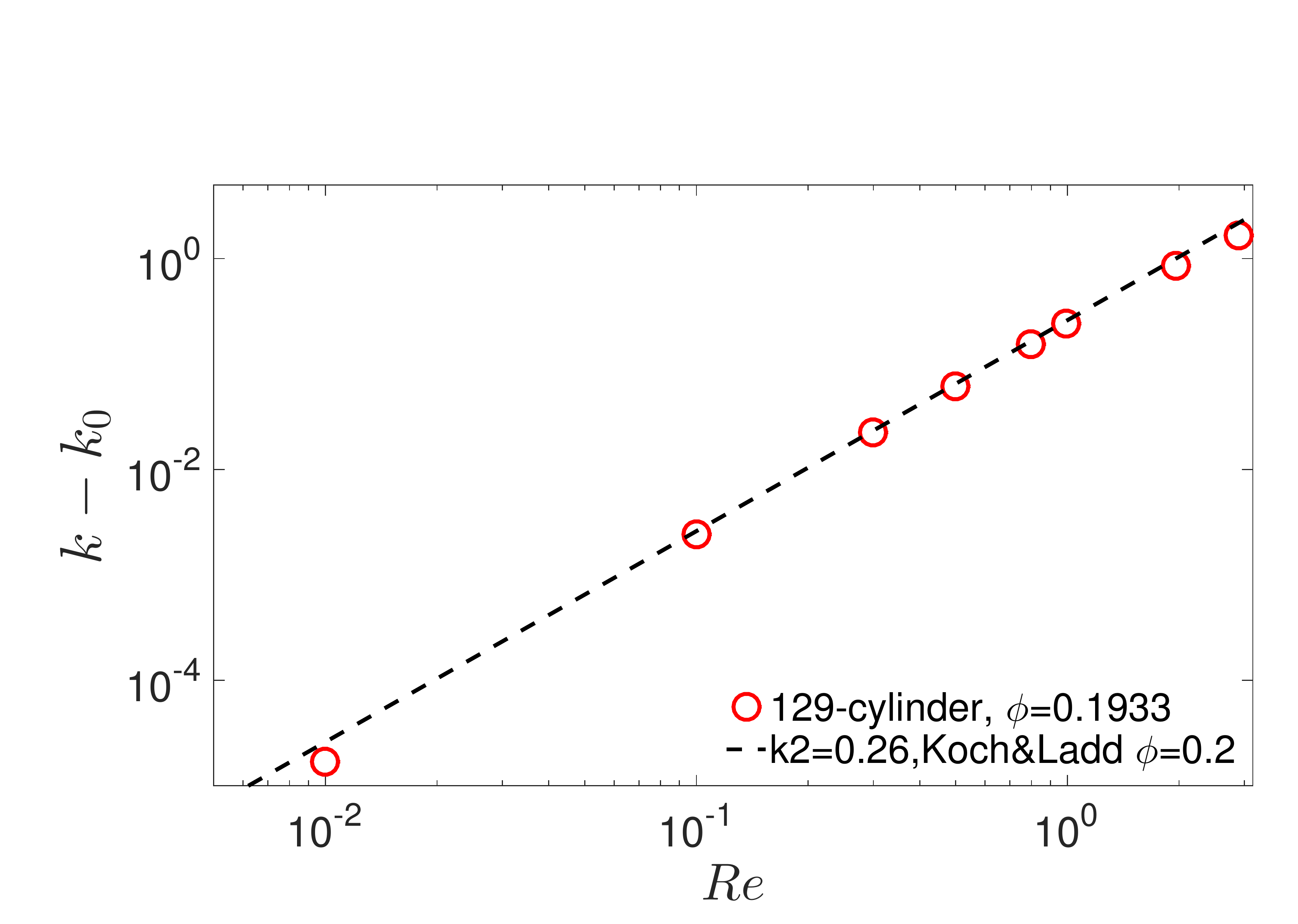}
\caption{\label{fig:Drag2D_array}The dimensionless excess (over Stokes flow)
drag coefficient for a square array of disks with packing fraction
$\phi\approx0.193$ (129-marker filled cylinder model, fluid grid
of $64^{2}$ cells). Comparison is made to known small-$\Re$ dependence
of the form $k_{0}+k_{2}\text{Re}^{2}$, with the coefficients $k_{0}$
and $k_{2}$ taken from the work of Koch and Ladd \citet{FiniteRe_2D_Ladd}
at $\phi=0.2$.}
\end{figure}

Fig. \ref{fig:Drag2D_array} shows the dimensionless excess drag $k-k_{0}$
as a function of $\Re$ at packing fraction $\phi=0.193$, which is
close to the packing fraction $\phi=0.2$ studied using the Lattice-Boltzmann
method in Ref. \citet{FiniteRe_2D_Ladd}. We see very good agreement
of the theoretical formula (\ref{eq:k_theory}) with our results using
the values of $k_{0}=49.2$ and $k_{2}=0.24$, which are in good agreement
with the values of $k_{0}=51.2$ and $k_{2}=0.26$ given in the caption
of Fig. 1 in Ref. \citet{FiniteRe_2D_Ladd}.

\subsubsection{Flow past a periodic column of cylinders}

Here we compute several solutions for flow past a column of cylinders
at somewhat larger Reynolds numbers, mimicking the setup of Ladd \citet{VACF_Ladd}.
The domain is a long narrow channel of $2048\times128$ grid cells
\footnote{For very elongated domains, our multigrid-based preconditioner converges
much faster for grid sizes that are powers of two.%
} with grid spacing $h=0.5$, keeping the markers at a distance $2h$.
Periodic boundary conditions are used in the direction of the short
side of the channel ($y$). The flow is driven by ``uniform'' inflow
and outflow boundary conditions in the long direction ($x$). Specifically,
we impose a specified normal velocity $V$ and zero tangential velocity
at both ends of the channel %
\footnote{An alternative is to use zero tangential stress on both boundaries,
or zero normal and tangential stress on the outflow; such stress boundary
conditions are supported in the fluid solver in the IBAMR library
\citet{NonProjection_Griffith}. %
}. The center of the cylinder is fixed at a quarter channel length
from the inlet. The cylinders in the periodic column are separated
by approximately 10 hydrodynamic radii in the $y$ direction (the
separation is $9.958R_{h}$ for the 121-marker cylinder, and $9.875R_{h}$
for the 39-marker shell).

\begin{figure}[h]
\begin{centering}
\includegraphics[width=0.49\textwidth]{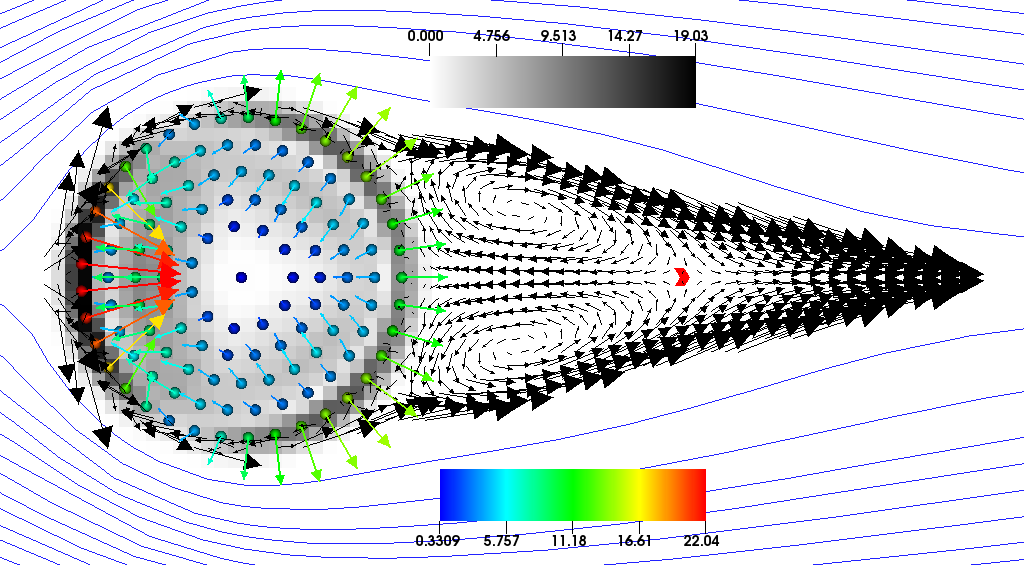}\includegraphics[width=0.49\textwidth]{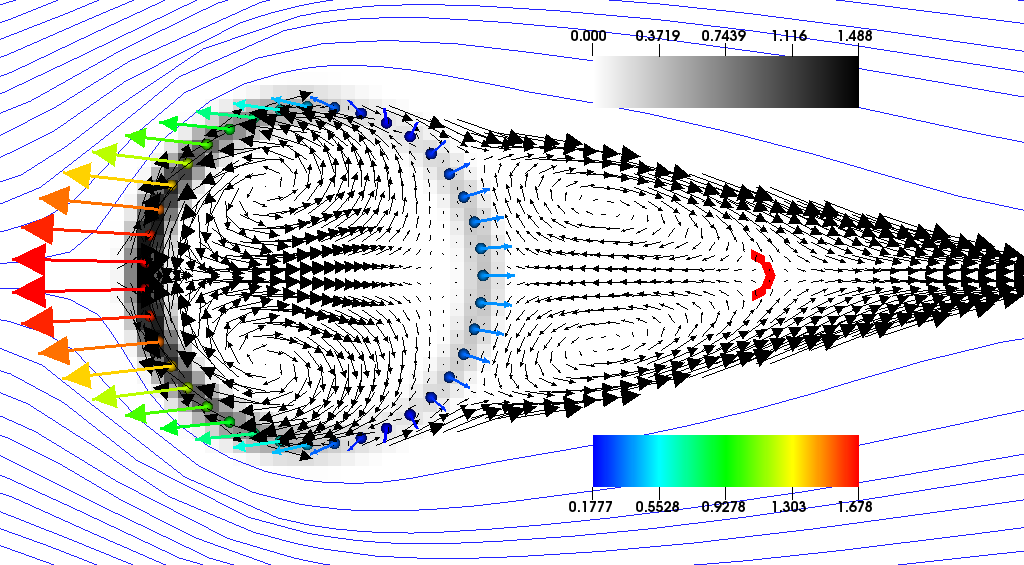}
\par\end{centering}

\centering{}\caption{\label{fig:2DWake}Steady incompressible flow at $\Re=10$ past a
periodic column of cylinders represented as either filled disks of
$121$ markers (left) or a shell of $39$ markers (right). We show
the magnitude of the Eulerian constraint force $\S\V{\Lambda}$ (gray
color map), the streamlines outside of the wake (solid blue lines),
the wake velocity field (black arrows), and the Lagrangian constraint
forces associated with each marker (color arrows). The red arrow marks
the stagnation point where $v_{x}=0$, as used to determine the wake
length.}
\end{figure}

Representative flow fields are shown in Fig. \ref{fig:2DWake} for
$\Re=10$ for a filled cylinder model (left) and an empty shell model
(right). Note that for the computation of total drag on a fixed cylinder
either model can be used since the spurious flow seen inside the empty
shell does not generate any overall acceleration of the fluid inside
the body. Also note that the spurious counter-rotating vortex pair
inside the shell diminishes under refinement, at an approximately
linear convergence rate, just as for steady Stokes flow. Computed
drag coefficients $k$ and wake lengths are shown in Table \ref{tab:Drag2D},
and good agreement is seen with the results of Lattice-Boltzmann and
finite difference schemes \citet{VACF_Ladd}. The wake length measures
the distance from the cylinder center to the stagnation point, which
is obtained by finding the largest $x$ coordinate on the contour
of zero horizontal velocity, $v_{x}=0$.

\begin{table}[h]
\caption{\label{tab:Drag2D}Numerical results for steady flow past a periodic
column of cylinders at different Reynolds number, for two different
models of the body (see Fig. \ref{fig:2DWake}), either a filled cylinder
or an empty shell of markers. For comparison we reproduce the results
in Table 5 in \citet{VACF_Ladd}, which are computed either using
either a Lattice-Boltzmann (LB) or a finite difference (FD) method.
(Left) Mean drag coefficient. (Right) Wake length in units of $R_{h}$. }

\centering{}%
\begin{tabular}{|l|c|c|c|c|}
\hline 
$\Re$  &
121 cyl  &
39 shell  &
LB  &
FD \tabularnewline
\hline 
5  &
4.31  &
4.35  &
4.21  &
4.32 \tabularnewline
\hline 
10  &
2.96  &
2.99  &
2.91  &
2.98 \tabularnewline
\hline 
20  &
2.16  &
2.19  &
2.17  &
2.19 \tabularnewline
\hline 
50  &
1.55  &
1.58  &
1.67  &
1.61 \tabularnewline
\hline 
\end{tabular}\hspace{1cm}%
\begin{tabular}{|l|c|c|c|c|}
\hline 
$\Re$  &
121 cyl  &
39 shell  &
LB  &
FD \tabularnewline
\hline 
5  &
1.52  &
1.40  &
1.5  &
1.49 \tabularnewline
\hline 
10  &
2.55  &
2.59  &
2.6  &
2.65 \tabularnewline
\hline 
20  &
4.50  &
4.61  &
4.7  &
4.74 \tabularnewline
\hline 
50  &
9.96  &
9.91  &
10.7  &
10.3 \tabularnewline
\hline 
\end{tabular}
\end{table}

\subsection{\label{sub:CubicSpheres}Flow past periodic arrays of spheres }

Finally, we study the drag on a cubic arrays of spheres of radius
$a$ at zero and finite Reynolds numbers, and compare our results
to those of Hill et al. \citet{FiniteRe_3D_Ladd}. At small packing
(volume) fractions $\phi$ and Reynolds numbers, according to Eqs.
(1-2) in \citet{SmallRe_3D_Ladd}, $F-F_{0}=3\Re/8+\mbox{h.o.t.}$
if $\sqrt{\phi}\ll\Re\ll1$, or, more relevant to our study, $F-F_{0}\sim\Re^{2}/\sqrt{\phi}$
if $\Re\ll\sqrt{\phi}\ll1$. For small $\text{Re}$ and at larger
densities, the theoretical arguments in \citet{FiniteRe_3D_Ladd,FiniteRe_2D_Ladd}
predict that the dimensionless drag is quadratic in $\Re$ because
the linear term vanishes by symmetry, so that 
\[
k=\frac{F}{6\pi\eta aV}\approx k_{0}+k_{2}\Re^{2}.
\]
For larger $\Re$, the dependence is expected to switch to linear
in $\Re$.

Here we focus on close-packed cubic lattices of spheres with packing
fraction $\phi=\pi/6\approx0.5236$. Note that unlike the case of
two spatial dimensions, in three dimensions the flow does not need
to squeeze in-between the (nearly) touching bodies, so the drag does
not diverge even at close packing. The value of the steady Stokes
drag $k_{0}$ is tabulated in Table \ref{tab:k_0_3D} for several
resolutions. Different resolutions are examined: an empty shell (see
Table \ref{tab:R_h_3D}) of 162 (grid size is $16^{3}$) or 642 markers
($30^{3}$ grid), as well as a filled sphere of 56 (42 on the surface,
$10^{3}$ grid) or 239 (162 on the surface, $16^{3}$ grid) markers;
the actual value of the packing fraction based on the effective hydrodynamic
radius of the model is indicated in the table. A large difference
is seen between the filled and empty shell models at this high packing
fractions because the spheres are very close to each other and discretization
artifacts become pronounced. We have also performed simulations at
a lower (but still high) packing fraction of $\phi=0.44$, and there
we see much better agreement between the filled and empty sphere models;
note that at small $\phi\ll1$ the value of $k_{0}$ must match among
resolutions since we \emph{define} the packing fraction from $R_{h}$,
which is itself determined from the value of $k_{0}$ at small $\phi$
using (\ref{eq:drag_periodic_2D}).

\begin{table}[h]
\caption{\label{tab:k_0_3D}Dimensionless drag force $k_{0}$ for steady Stokes
flow ($Re=0$) past a simple-cubic array of spheres at volume fraction
$\phi\approx\pi/6$ (close packing). For the highest-resolution LB
simulations in \citet{VACF_Ladd} the reported value is $k_{0}=42.8$. }

\centering{}%
\begin{tabular}{|l|c|c|}
\hline 
Number of markers  &
$\phi$  &
$k_{0}$ \tabularnewline
\hline 
56 filled &
0.5236  &
40.08 \tabularnewline
\hline 
239 filled &
0.5238  &
40.73 \tabularnewline
\hline 
\hline 
162 shell &
0.5213  &
44.49 \tabularnewline
\hline 
642 shell &
0.5236  &
43.29 \tabularnewline
\hline 
\end{tabular}
\end{table}

\begin{figure}[h]
\centering{}\includegraphics[width=0.49\textwidth]{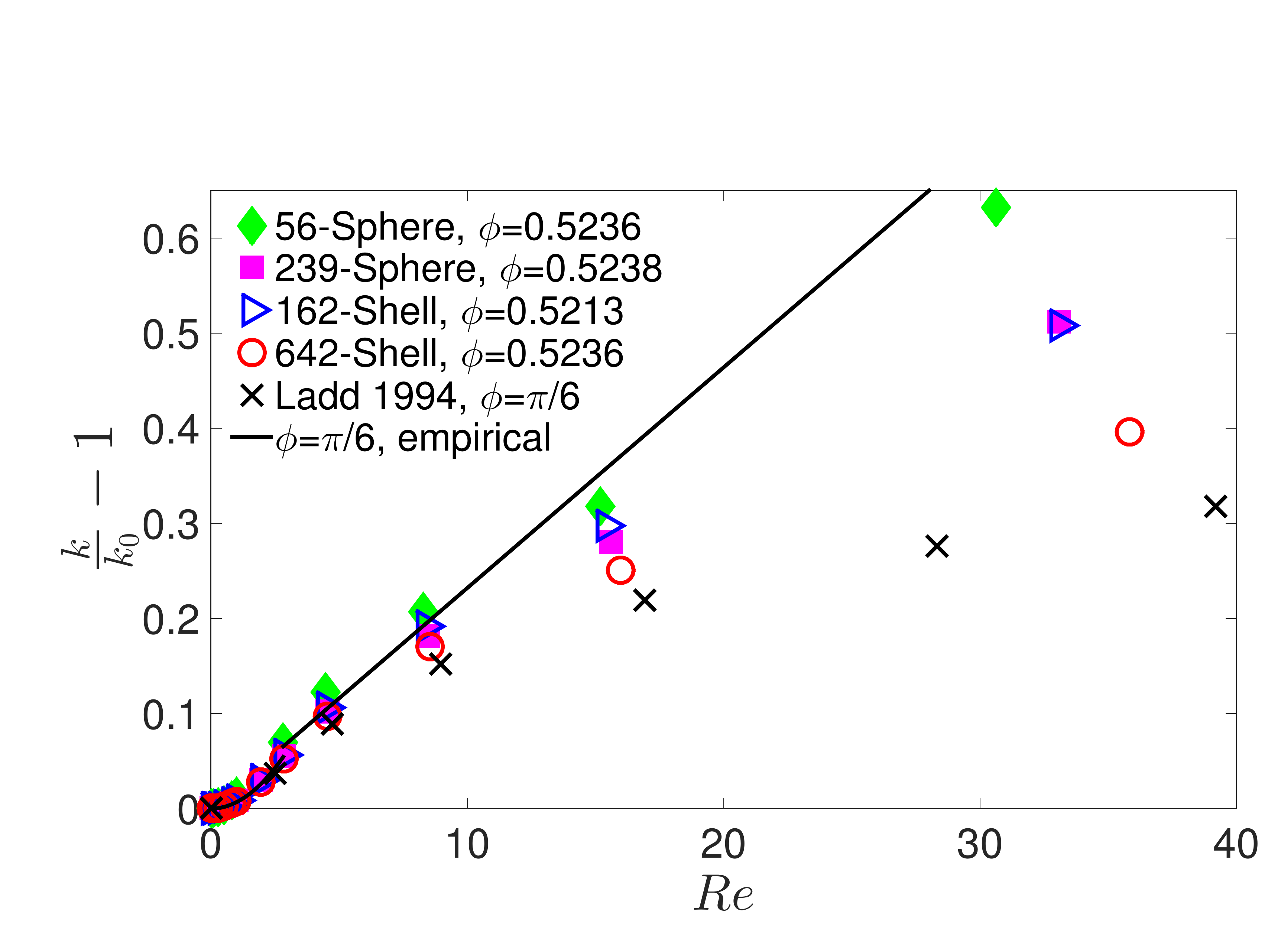}\includegraphics[width=0.49\textwidth]{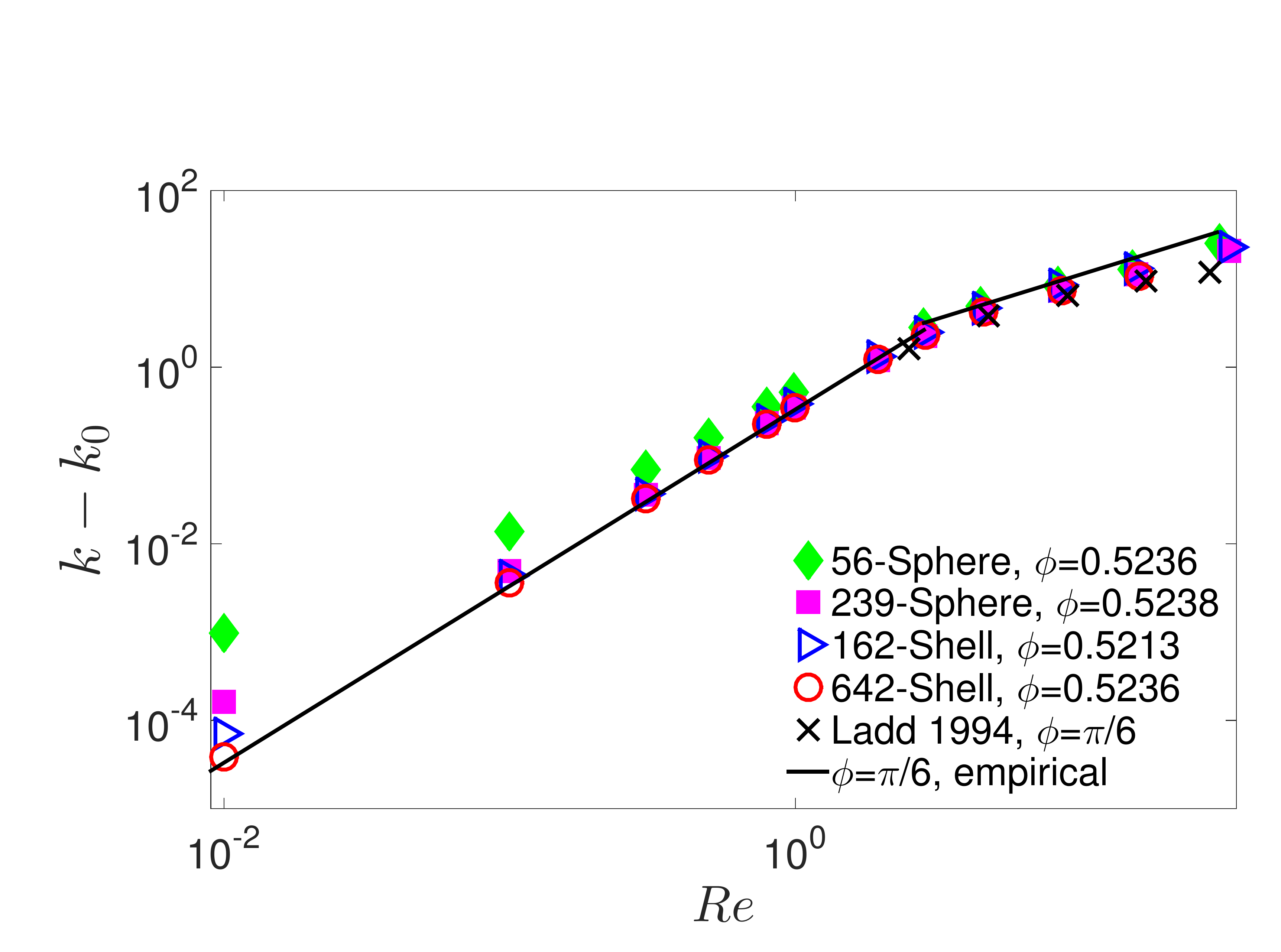}
\caption{\label{fig:Drag3D_array}Numerical values (symbols) for the drag coefficient
of a periodic array of spheres close packed in a cubic lattice of
volume fraction $\phi=\pi/6\approx0.52$, for several resolutions
(see legend), using a linear (left) or log scaling (right). Comparison
is made to empirical formulas given in \citet{EmpiricalDrag_3D} (lines),
as well as Lattice Boltzmann results for close-packed cubic arrays
given in Table 6 in \citet{VACF_Ladd} (crosses).}
\end{figure}

Numerical results for the dimensionless drag coefficient $k$ near
the close-packed density $\phi\approx0.52$ are shown in Fig. \ref{fig:Drag3D_array}.
Because our discrete models of spheres do not have the same symmetry
as a perfect sphere, we numerically observe a small $O(\Re)$ correction
that can dominate the true correction $k_{2}\Re^{2}$ for $\Re\ll1$;
this is especially evident in the right panel of Fig. \ref{fig:Drag3D_array}
for coarsely-resolved models (e.g., a 56-marker sphere) %
\footnote{In two dimensions, we can more easily make the discrete models symmetric
and this is why Fig. \ref{fig:Drag2D_array} does not show deviations
from the expected quadratic behavior even at rather small $\Re$.%
}. Empirical fits to literature data for $k_{0}$, $k_{1}$, $k_{2}$,
and the range of $\Re$ values over which the various fits are valid
are tabulated in Ref. \citet{EmpiricalDrag_3D}. Fig. \ref{fig:Drag3D_array}
compares our results to these fits, as well as to reference results
obtained using the Lattice Boltzmann method \citet{VACF_Ladd} at
close packing. We observe the expected switch from linear to quadratic
dependence on $\Re$ and also a reasonable agreement with the literature
data, and the agreement appears to improve with increasing resolution.

\section{Conclusions}

This paper develops an immersed boundary method that enforces strict
rigidity of immersed bodies at both zero and finite Reynolds numbers.
Unlike existing approaches, we do not rely on penalty or splitting
approaches, and we instead directly solve a saddle-point system that
couples the fluid velocity and pressure to the unknown rigidity forces.
We developed a physics-inspired approximation $\widetilde{\Mob}$
to the Schur complement (mobility matrix) $\Mob$ of the constrained
system, based on analytical considerations for a continuum fluid model,
and demonstrated that this leads to a robust preconditioner so long
as the immersed boundary markers are kept sufficiently far to ensure
a well-conditioned mobility matrix. Contrary to common practice, we
found that the markers should be kept approximately two fluid grid
cells apart in rigid-body models in order to obtain accurate and stable
pointwise estimates for the traction. We tested our method on a number
of standard test problems in both two and three spatial dimensions,
and at both zero and finite Reynolds number, and we observed good
agreement with theory and literature values. Although in this work
we focused on rigid bodies, our method can directly be applied to
study fluid flow around bodies with specified kinematics. For example,
it can be used to model the flow around a swimming body deforming
with a specified gait. We have implemented the method described here
in the open-source IBAMR software infrastructure \citet{IBAMR} in
the hope it will be useful to other users of the IB method.

Another challenge that we did not explore here is the efficient computation
of the action of $\widetilde{\Mob}^{-1}$ when there are many markers
present; there are many approximate solvers and emerging fast solvers
we plan to explore in the future. Of course, using dense linear algebra
to solve (\ref{eq:mob_subproblem}) is likely to be suboptimal, as
these solves have $O(N^{2})$ memory complexity and $O(N^{3})$ time
complexity. The problem of solving a linear systems similar in structure
to (\ref{eq:mob_subproblem}) appears in many other methods for hydrodynamics
of suspensions, including Brownian \citet{BrownianDynamics_OrderN,BrownianDynamics_FMM}
and Stokesian \citet{BrownianDynamics_OrderNlogN} dynamics, the method
of regularized Stokeslets \citet{RegularizedStokeslets_2D,RegularizedStokeslets},
computations based on bead models of rigid bodies \citet{HYDROLIB,SphereConglomerate,HYDROPRO,HYDROPRO_Globular},
and first-kind boundary integral formulations of Stokes flow \citet{BoundaryIntegral_Pozrikidis}.
Similar matrices appear in static Poisson problems such as electrostatics
or reaction-diffusion models \citet{ReactiveBlobs}, and there is
a substantial ongoing work that can be applied to our problem. Notably,
the approximate mobility matrix $\widetilde{\Mob}$ is dense but has
a well-understood low-rank structure that can be exploited. Specifically,
matrix-vector products $\widetilde{\Mob}\Lamb$ can be performed in
almost linear time using the Fast Multipole Method \citet{RPY_FMM}.
If the condition number of $\widetilde{\Mob}$ is not too large, one
can solve linear systems involving $\widetilde{\Mob}$ efficiently
using an unpreconditioned Krylov solver. For poorly conditioned cases,
however, a good preconditioner based on an approximate factorization
of $\widetilde{\Mob}$ is required. In recent years, several approximate
low-rank factorizations of matrices of this type have been developed
\citet{FastHierarchicalSolver,FastSolver_Ho,HODLR_BDLR}, and can
be used as preconditioners in Krylov methods. We have had reasonable
success using a fast hierarchically off-diagonal low-rank (HODLR)
factorization code developed by Ambikasaran and Darve \citet{FastHierarchicalSolver},
with significant improvement offered by a recently-developed boundary
distance low-rank approximation \citet{HODLR_BDLR}. Preliminary results
indicate great promise for the inverse fast multipole (iFMM) method
\citet{iFMM}; we have been able to use iFMM to solve the system (\ref{eq:mob_subproblem})
for as many as $5\cdot10^{5}$ markers to a relative tolerance of
$10^{-8}$. These methods are, however, still under active development,
and a significant amount of investigation is necessary to integrate
them into the method described here. Notably, we only require an approximate
solver for (\ref{eq:mob_subproblem}) and the impact of the innacuracy
in solving (\ref{eq:mob_subproblem}) on the overall convergence of
the outer Krylov solver needs to be assessed.

The type of linear system we solve here is closely connected to those
appearing in implicit immersed boundary methods \citet{IBM_Implicit_Comparison,IBMultigrid_Guy,ImplicitIB_Projection}.
It is in fact possible to recast the saddle-point problem we consider
here into a form closely-related to that appearing in implicit IB
methods; the Schur complement for this system is in Eulerian rather
than Lagrangian variables as it was for this work, and involves the
matrix
\begin{equation}
\Lap+\kappa\S\left(\J\S\right)^{-1}\J,\label{eq:Schur_Eulerian}
\end{equation}
for some constant $\kappa$ that does \emph{not} need to go to infinity.
It may be that geometric multigrid methods \citet{IBMultigrid_Guy}
developed for implicit IB methods can be applied to the Eulerian Schur
complement (\ref{eq:Schur_Eulerian}). At the same time, techniques
developed herein may be useful in the development of more efficient
implicit IB methods for nearly-rigid bodies.

Our work is only the first step toward the ultimate goal of developing
methods able to handle large numbers of rigid bodies in flow. Several
computational challenges need to be tackled to realize this goal.
Firstly, and most importantly, it is crucial to develop a preconditioner
for the enlarged linear system (\ref{eq:free_kinematics_Stokes})
that appears in the context of freely-moving rigid bodies. An additional
Schur complement appears when solving this saddle-point problem, and
the challenge for future work is approximating the body mobility matrix
$\M{\mathcal{N}}=\left(\M{\mathcal{K}}^{\star}\M{\mathcal{M}}^{-1}\M{\mathcal{K}}\right)^{-1}.$
Initial investigations have shown great promise in block-diagonal
preconditioners with one block per body. In this approach, we neglect
the hydrodynamic interactions between bodies, but use the mobility
approximation developed in this work together with dense linear algebra
for each body.

In the marker-based method described in this work, one must adjust
the marker spacing to be ``neither too small nor too large''. The
sensitivity of the solver performance and the numerical results to
the exact spacing of the markers, which comes from the ill-conditioning
of the mobility matrix, is one of the key deficiencies of the marker-based
representation inherent to the traditional IB method. Recently, Griffith
and Luo have proposed an alternative IB approach that models the deformations
and stresses of immersed elastic body using a finite element (FE)
representation \citet{IBFE}. In their IB/FE approach, the degrees
of freedom associated with $\Lamb$ are represented on an FE mesh
that may be coarser than the fluid grid, and the interaction between
the fluid grid and body mesh is handled by placing IB markers at the
numerical quadrature points of the FE mesh. When such an approach
is generalized to rigid bodies, the conditioning of the mobility becomes
much less sensitive to the marker spacing. Using a finite-element
basis to represent the unknown fluid-body interaction force amounts
to applying a \emph{filter $\M{\Psi}$} to the marker-based mobility
matrix, which is a well-known and robust technique to regularize ill-conditioned
systems. Specifically, in the context of the IB/FE approach, the mobility
operator becomes 
\[
\M{\mathcal{M}}_{FE}=\M{\Psi}\left(\M{\mathcal{J}}\M{\mathcal{L}}^{-1}\M{\mathcal{S}}\right)\M{\Psi}^{T}=\M{\Psi}\Mob\M{\Psi}^{T},
\]
where $\M{\Psi}$ is a matrix that contains quadrature weights as
well as geometric information about the relation between the nodes
and quadrature points of the FE mesh. The FE mobility matix $\M{\mathcal{M}}_{FE}$
is still symmetric, but now can be much smaller because the number
of unknowns is equal to the number of FE degrees of freedom rather
than the number of markers. Even if markers are closely spaced, the
filtering of the high-frequency modes performed by representing forces
in a smooth FE basis makes the mobility much better conditioned than
for marker-based schemes. Furthermore, the mobility matrix, or approximations
of it used for preconditioning, will be smaller and thus easier to
fit in memory. We also expect the resulting method to be more accurate
because the tractions are represented in a smoother basis. We will
explore this promising extension of our rigid-body IB methods in future
work.
\begin{acknowledgments}
We thank John Brady, Bob Guy, Anthony Ladd, Neelesh Patankar, and
Charles Peskin for numerous stimulating and informative discussions.
We thank Sivaram Ambikasaran, Pieter Coulier, AmirHossein Aminfar,
and Eric Darve for their help with low-rank approximate factorizations
of mobility matrices. A. Donev and B. Kallemov were supported in part
by the Air Force Office of Scientific Research under grant number
FA9550-12-1-0356, as well as by the National Science Foundation under
award DMS-1418706. B. E. Griffith and A. P. S. Bhalla were supported
in part by the National Science Foundation under awards DMS-1016554
and ACI-1047734 (to New York University School of Medicine) and awards
DMS-1460368 and ACI-1460334 (to the University of North Carolina at
Chapel Hill).
\end{acknowledgments}
\appendix

\section*{Appendix}

\section{\label{AppendixMob3D}Approximating the mobility in three dimensions}

In this appendix we give the details of our empirical fits for the
approximations to the functions $f_{\beta}(r)$ and $g_{\beta}(r)$
in (\ref{eq:M_tilde_ij}) in three spatial dimensions, following the
physics-based constraints discussed in Section \ref{sec:ApproxMob}.
To maximize the quality of the fit, we perform separate fits for $\beta\rightarrow\infty$
(steady Stokes flow) and finite $\beta$. We also make an effort to
make the fits change smoothly as $\beta$ grows towards infinity.

\subsubsection{Steady Stokes flow}

Because our numerical computations are done in a periodic domain of
length $l$ rather than an unbounded domain, we need to apply a well-known
correction to the Oseen tensor \citet{Mobility2D_Hasimoto,ISIBM},
\[
f_{\infty}(h\ll r\ll l)\approx\left(8\pi\eta r\right)^{-1}-2.84\left(6\pi\eta l\right)^{-1}.
\]
From the numerical data, we calculated the normalized functions
\begin{align}
\tilde{f}(x) & =\left(8\pi\eta r\right)\left(f_{\infty}(r)+2.84/\left(6\pi\eta l\right)^{-1}\right),\label{eq:fg_tilde_Stokes}\\
\tilde{g}(x) & =\left(8\pi\eta r\right)g_{\infty}(r),\nonumber 
\end{align}
where $x=r/h$ is the normalized distance between the markers. As
explained previously, we know that $\tilde{f}\approx\left(8\pi\eta r\right)/\left(6\pi\eta a\right)=4r/(3a)$
for $x\ll1$ (in practice, markers are never too close to each other
so we only need the self-mobility, i.e., $x=0$), and that $\tilde{g}$
grows at least quadratically for small $x$ (since $g(0)=0$). We
also know that $\tilde{f}\approx1$ and $\tilde{g}\approx1$ for large
$r\gg h$. The numerical data for the normalized functions $\tilde{f}(x)$
and $\tilde{g}(x)$ are shown in Fig. \ref{fig:fg_3d_asympt} along
with fits to the following semi-empirical rational functions,
\begin{align}
\tilde{f}(x) & =\begin{cases}
\frac{x}{\left(3a\right)/\left(4h\right)+b_{0}x^{2}} & \quad\text{if }x<0.8,\\
b_{1}xe^{-b_{2}x}+\frac{b_{3}x^{2}+x^{4}}{1+b_{4}x^{2}+x^{4}} & \quad\text{if }x\geq0.8,
\end{cases}\label{eq:f_stokes_3D_fit}\\
\tilde{g}(x) & =\frac{x^{3}}{b_{5}+b_{6}x^{2}+x^{3}}.\nonumber 
\end{align}
As the figure shows, the numerical data are well described by these
formulas, and there is only small scatter of the numerical data around
the fit, indicating approximate discrete translational and rotational
invariance %
\footnote{Most of the scatter comes from the finite size of the periodic box
and can be explained using a known periodic correction to the RPY
tensor \citet{RotnePrager_Periodic}.%
}. We also obtain a reasonable agreement with the RPY tensor (\ref{eq:RPYTensor})
approximation; however, as expected, the empirical fits yield a better
match to the data.

\subsubsection{Nonzero Reynolds numbers}

For finite $\beta$, we consider separately the case $r=0$ (giving
the diagonal elements $\widetilde{\Mob}_{ii}$) and $r>0.1h$ (giving
the off-diagonal elements). For $r=0$ we use an empirical fit designed
to conform to (\ref{eq:f_beta_r_0}), 
\begin{eqnarray}
\varphi_{0}\left(\beta\right)=\frac{\eta hf_{\beta}(0)}{\beta} & = & \frac{1+z_{1}\sqrt{\beta}+z_{2}\beta}{z_{0}+z_{3}\beta+6\pi(a/h)z_{2}\beta^{2}},\\
g_{\beta}(0) & = & 0,\nonumber 
\end{eqnarray}
where $z_{1}-z_{3}$ are coefficients obtained by fitting the numerical
data for the self mobility for different $\beta$. Note that $z_{0}=2h^{3}/\left(3V_{m}\right)$
is fixed by the inviscid condition (\ref{eq:f_0_inviscid}). Also
note that as $\beta\rightarrow\infty$, our fit obeys the correct
Stokes limit,
\[
\varphi_{0}\left(\beta\gg1\right)\rightarrow\frac{1}{6\pi\left(a/h\right)}\cdot\frac{1}{\beta}.
\]
We show the empirical fit for $\varphi_{0}\left(\beta\right)$ in
Fig. \ref{fig:f_0_3D} in Appendix \ref{AppendixMob2D}.

For nonzero $r$, we introduce normalized functions $\tilde{f}_{\beta}$
and $\tilde{g}_{\beta}$ via 
\begin{eqnarray}
f_{\beta}(r) & = & -\frac{\beta}{\eta h}\cdot\frac{1}{4\pi x^{3}}\cdot\tilde{f}_{\beta}(x),\\
g_{\beta}(r) & = & \frac{\beta}{\eta h}\cdot\frac{3}{4\pi x^{3}}\cdot\tilde{g}_{\beta}(x),\nonumber 
\end{eqnarray}
where $x=r/h$ is the normalized distance, and $\beta/\eta=\D t/\rho h^{2}$.
For finite $\beta$, we know that $\tilde{f}_{\beta}(x\gg\sqrt{\beta})\approx\tilde{g}_{\beta}(x\gg\sqrt{\beta})\approx1$
according to (\ref{eq:fg_inv_3D}). As $\beta\rightarrow\infty$,
we want to reach the Stokes limit
\begin{align}
\tilde{f}_{\infty}\left(x\gg1\right) & \rightarrow-\frac{x^{2}}{2\beta},\\
\tilde{g}_{\infty}\left(x\gg1\right) & \rightarrow\frac{x^{2}}{6\beta},\nonumber 
\end{align}
and for finite $\beta$, we want the viscous contribution to decay
as $\exp\left(-x/\left(C\sqrt{\beta}\right)\right)$ for some constant
$C$ that should be close to unity. Furthermore, we would like to
ensure continuity near the origin with the fit for $r=0$,
\[
\tilde{f}_{\beta}(x\rightarrow0)\rightarrow-4\pi x^{3}\varphi_{0}\left(\beta\right).
\]

A fitting formula that obeys these conditions that we find to work
well for $r>0.1h$ is 
\begin{eqnarray}
\tilde{f}_{\beta}(x) & = & \varphi_{0}(\beta)\frac{-4\pi x^{3}+a_{4}\left[x^{5}-x^{7}e^{(-a_{3}x/\sqrt{\beta})}/(2\beta)\right]}{1+a_{0}x+a_{1}x^{2}+a_{2}x^{3}+a_{4}x^{5}\varphi_{0}(\beta)}+\nonumber \\
 &  & +\frac{a_{5}x^{4}e^{-a_{6}x}+a_{7}x^{4}}{1+a_{8}x^{3}+a_{9}x^{5}},\label{eq:f_beta_3D_fit}\\
\tilde{g}_{\beta}(x) & = & \varphi_{0}(\beta)\frac{b_{5}\left[x^{5}+x^{7}e^{(-b_{0}x/\sqrt{\beta})}/(6\beta)\right]}{1+b_{1}x+b_{2}x^{2}+b_{3}x^{3}+b_{4}x^{4}+b_{5}\varphi_{0}(\beta)x^{5}},\nonumber 
\end{eqnarray}
where $a_{0}$-$a_{9}$ and $b_{0}$-$b_{5}$ are empirical coefficients.
It is important to emphasize that (\ref{eq:f_beta_3D_fit}) was chosen
in large part based on empirical trial and error. Many other alternatives
exist. For example, one could use the analytical Brinkmanlet (\ref{eq:Brinkmanlet_3D})
for sufficiently large distances and then add short-ranged corrections
for nearby markers. Alternatively, one could first subtract the inviscid
part $f_{0}(r)$ and $g_{0}(r)$ and then fit the viscous contribution
only. As discussed above, ideally the fits would be constrained to
guarantee an SPD approximate mobility matrix, but this seems difficult
to accomplish in practice.

\begin{figure}[tbph]
\begin{centering}
\includegraphics[width=0.49\textwidth]{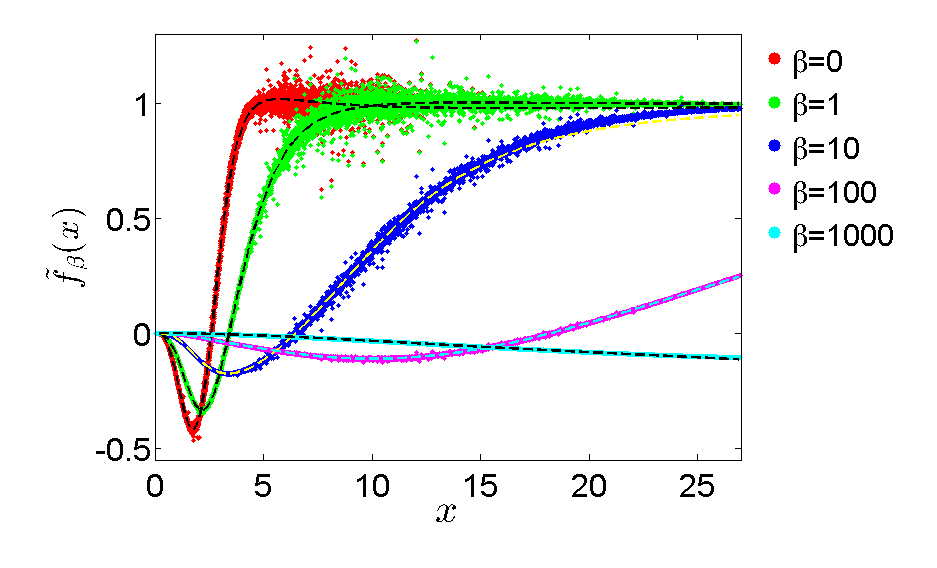}\includegraphics[width=0.49\textwidth]{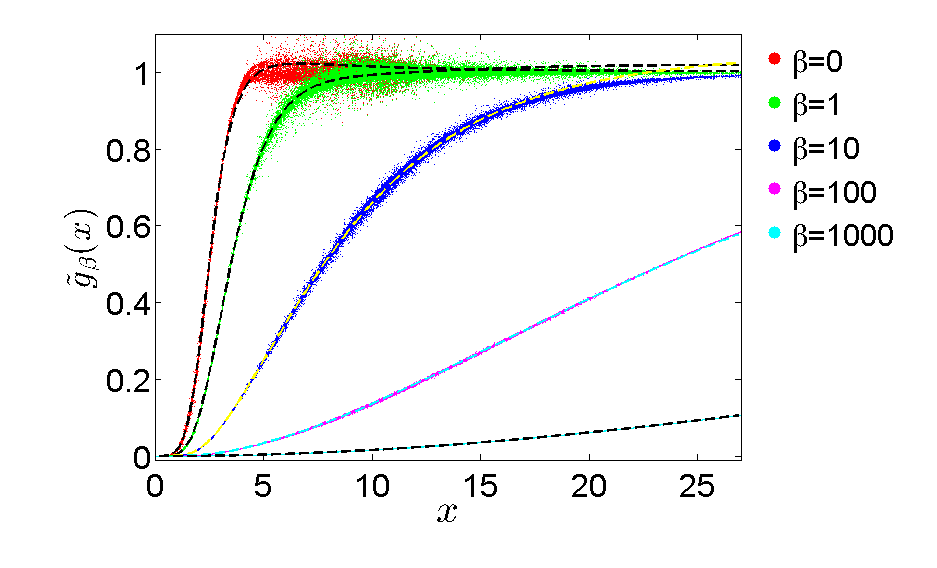}
\par\end{centering}

\caption{\label{fig:fg_3d_fit}Normalized fitting functions $\tilde{f}_{\beta}(x)$
(left) and $\tilde{g}_{\beta}(x)$ (right) at finite $\beta$ in three
dimensions for the 6-point kernel, for different values of the viscous
CFL number (see legend). Symbols are numerical data obtained by using
a $256^{3}$ periodic fluid grid, and dashed lines show the best fit
of the form (\ref{eq:f_beta_3D_fit}).}
\end{figure}

We computed the fitting coefficients in (\ref{eq:f_beta_3D_fit})
for $\beta\in\{0,\ 0.1,\ 0.25,\ 0.5,\ 1,\ 10,\ 100,\ 1000\}$; the
coefficients for other values in the range $0<\beta<1000$ are interpolated
using linear interpolation, and $\beta>1000$ is treated using the
steady Stokes fitting. We see a good match between the numerical data
and our empirical fits in Fig. \ref{fig:fg_3d_fit}, with good translational
and rotational invariance (i.e., relatively small scatter of the numerical
points around the fits).

\section{\label{AppendixMob2D}Approximating the mobility in two dimensions}

To construct empirical approximations to the functions $f_{\beta}(r)$
and $g_{\beta}(r)$ in (\ref{eq:M_tilde_ij}) in two spatial dimensions,
we follow the same approach as we did for three dimensions in Appendix
\ref{AppendixMob3D}. Specifically, we first discuss the known asymptotic
behavior of these functions at short and large distances, and use
this to guide the construction of empirical fitting formulas.

\subsection{Physical constraints}

In two dimensions, we need to modify (\ref{eq:f_beta_r_0}) to agree
with (\ref{eq:f_0_inviscid}) for small $\beta$. For $d=2$, $V_{m}=c_{V}^{'}h^{2}$
and $f_{0}(0)\sim\beta/\eta$ so that we use the fit
\begin{equation}
f_{\beta}\left(0\right)=\frac{C\left(\beta\right)}{\eta}\mbox{ and }g_{\beta}(0)=0,\label{eq:f_beta_r_0_2D}
\end{equation}
where $C(\beta)$ has the same asymptotic scaling as in three dimensions
and is obtained from empirical fits (see Fig. \ref{fig:f_0_2D}).
A key difference exists between two and three spatial dimensions in
the limit $\Re\rightarrow0$. For steady Stokes flow in a square two
dimensional periodic domain, the Green's function diverges logarithmically
with the system size $l$. Therefore, it is not possible to write
a formula for the asymptotic behavior at large distances for an infinite
system. Instead, we must subtract the divergent piece to get a well-defined
answer. The standard Green's function for Stokes flow in two dimensions
has logarithmic growth at infinity, which suggests that (\ref{eq:fg_asympt_3D})
should be replaced by 
\begin{equation}
f_{\infty}(r\gg h)-f_{\infty}(r=0)\approx-\frac{\ln\,\left(r/h\right)}{4\pi\eta}\mbox{ and }g_{\infty}(r\gg h)\approx\frac{1}{4\pi\eta}.\label{eq:fg_asympt_2D}
\end{equation}
For inviscid flow we should replace (\ref{eq:fg_inv_3D}) by the field
of a dipole in two dimensions,
\begin{equation}
f_{0}(r\gg h)\approx-\frac{\D t}{2\pi\rho r^{2}}\mbox{ and }g_{0}(r\gg h)\approx\frac{\D t}{\pi\rho r^{2}}.\label{eq:fg_inv_2D}
\end{equation}
In two dimensions the solutions of the Brinkmann equation (\ref{eq:Brinkman_Greens})
are analytically complicated and involve special functions. Even without
solving these equations, however, physical scaling suggests that the
same physical length scale $h\sqrt{\beta}$ should enter, in particular,
the viscous corrections should decay to zero exponentially fast with
$h\sqrt{\beta}$.

\subsection{Empirical fits}

We have used the analytical results above to construct empirical fitting
formulas that have the correct asymptotic behavior, as we now explain
in more detail.

\subsubsection{Steady Stokes flow}

In two dimensions, steady Stokes flow ($\beta\rightarrow\infty$)
is not well behaved because the Green's function does not decay sufficiently
rapidly (Stokes paradox). This makes the mobility an essentially dense
matrix that is sensitive to boundary conditions and difficult to approximate.
Nevertheless, we have used a periodic system to fit empirical data
based on the theory (\ref{eq:fg_asympt_2D}). The diagonal value $f_{\infty}(0)$
diverges logarithmically with the system size $L$ for periodic boundaries.
Specifically, for a square unit cell of length $l\gg h$, it is known
that \citet{Mobility2D_Hasimoto}
\[
f_{\infty}(0)=\left(4\pi\eta\right)^{-1}\ln\left(\frac{l}{3.708\, a}\right),
\]
and this relation defines the effective hydrodynamic radius of a marker
$a$ (note that $a/h$ is a universal value for a given spatial discretization,
as it is in three dimensions). Since the precise form depends on boundary
conditions and is not known in general, we treat $f_{\infty}(0)$
as an \emph{input} parameter. 

\begin{figure}[tbph]
\begin{centering}
\includegraphics[width=0.49\textwidth]{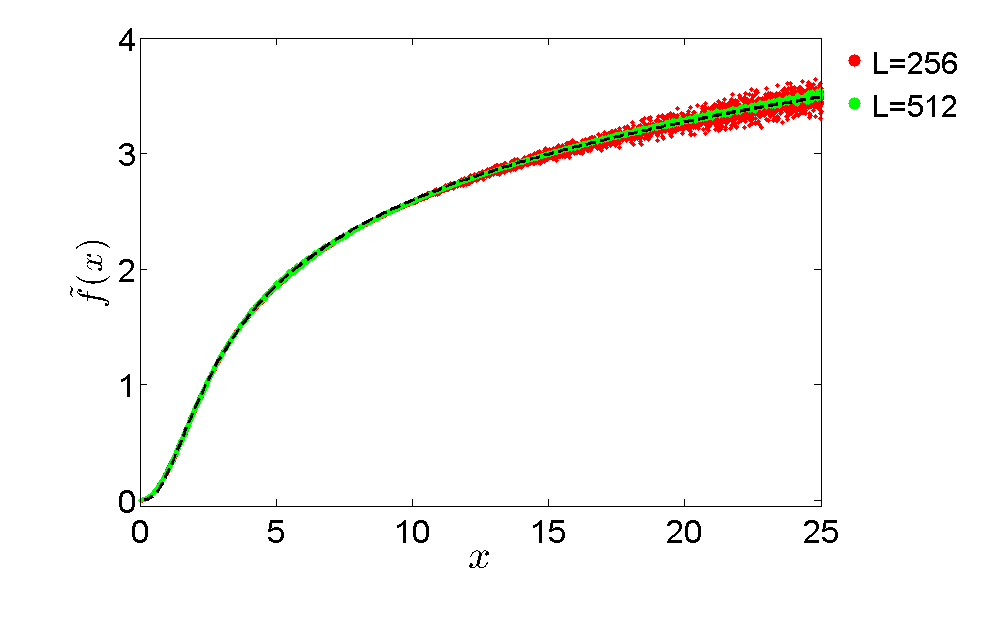}\includegraphics[width=0.49\textwidth]{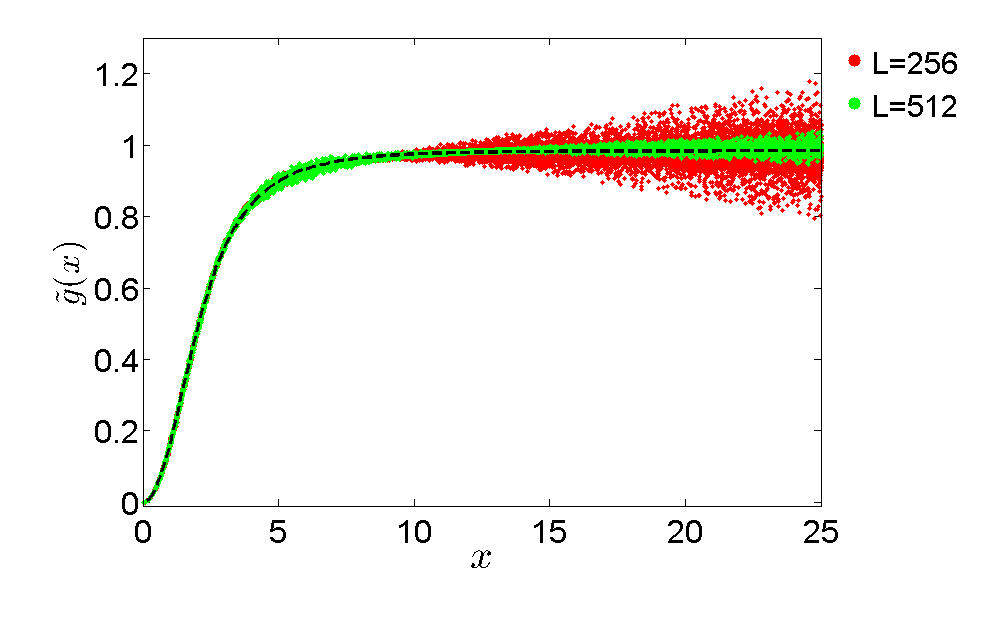}
\par\end{centering}

\caption{\label{fig:fg_steady_2d}Empirical fits (lines) to numerical data
(symbols) for $\tilde{f}(x)$ (left) and $\tilde{g}(x)$ (right),
for the 6-point kernel in two dimensions, obtained using a periodic
system of either $256^{2}$ or $512^{2}$ grid cells. Observe that
both follow the correct asymptotic behavior at large distances, with
scatter dominated by boundary effects.}
\end{figure}

From the numerical data, we calculated the normalized functions 
\begin{align}
\tilde{f}(x) & =-\left(4\pi\eta\right)\left(f_{\infty}(x)-f_{\infty}(0)\right)\label{eq:fg_tilde_Stokes_2d}\\
\tilde{g}(x) & =\left(4\pi\eta\right)g_{\infty}(x),\nonumber 
\end{align}
where $x=r/h$ is the normalized distance between the markers. Observe
that from (\ref{eq:fg_asympt_2D}) we know that $\tilde{f}(x\gg1)\approx\ln x$
and $\tilde{g}(x\gg1)\approx1$. For the normalized functions, we
use the fits 
\begin{eqnarray}
\tilde{f}(x) & = & \frac{a_{0}x^{2}+a_{1}x^{3}+a_{2}x^{3}\ln x}{1+a_{3}x+a_{4}x^{2}+a_{2}x^{3}},\label{eq:fg_steady_2d}\\
\tilde{g}(x) & = & \frac{b_{0}x^{2}+b_{1}x^{3}}{1+b_{2}x+b_{3}x^{2}+b_{1}x^{3}}.\nonumber 
\end{eqnarray}
Numerical results and empirical fits for $\tilde{f}(x)$ and $\tilde{g}(x)$
are shown in Fig. \ref{fig:fg_steady_2d}. While the numerical data
do conform to the theoretical asymptotic behavior, there is substantial
scatter for larger distances because of the strong sensitivity to
the boundaries.

\subsubsection{Nonzero Reynolds numbers}

\begin{figure}[h]
\centering{}\includegraphics[width=0.49\textwidth]{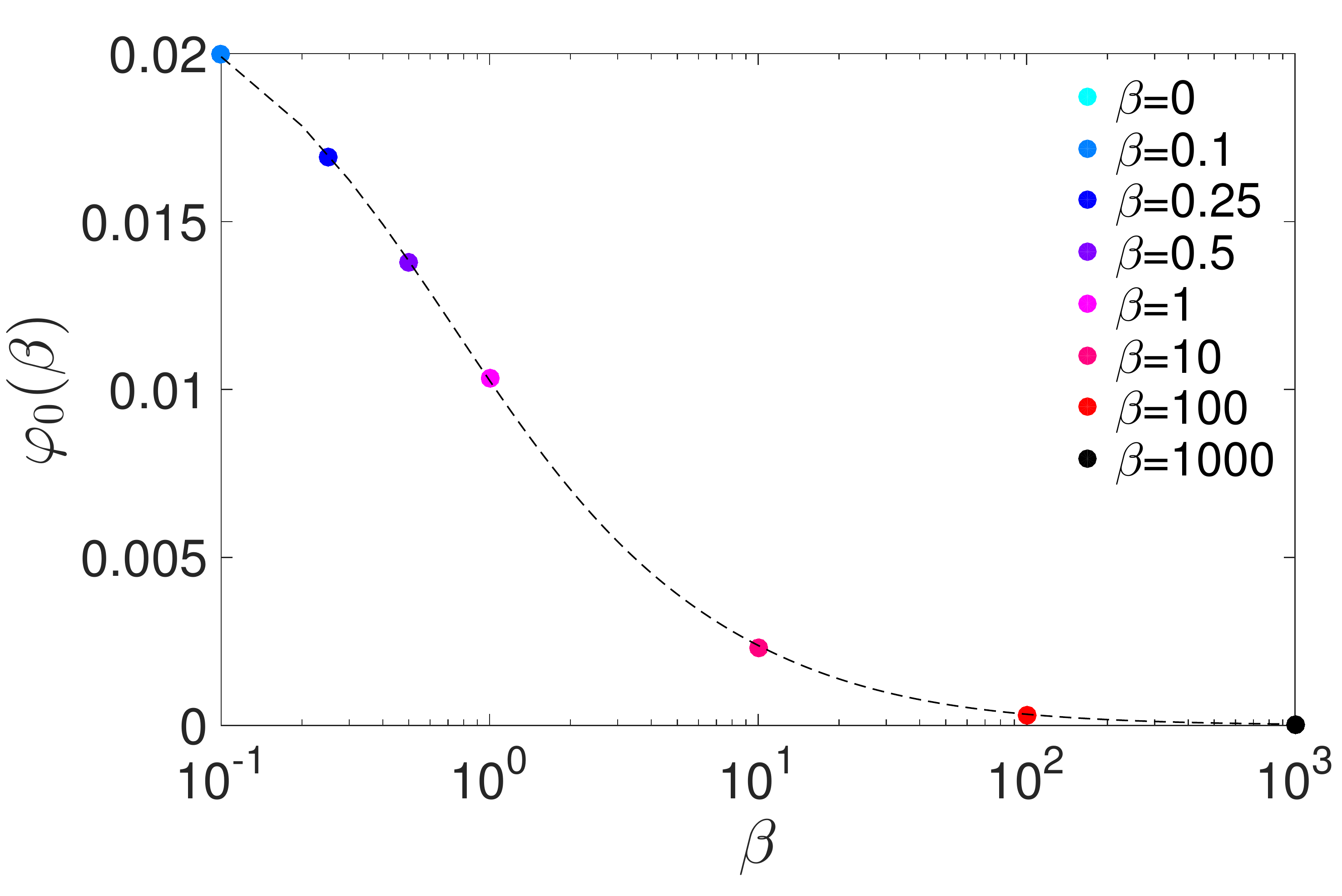}\includegraphics[width=0.49\textwidth]{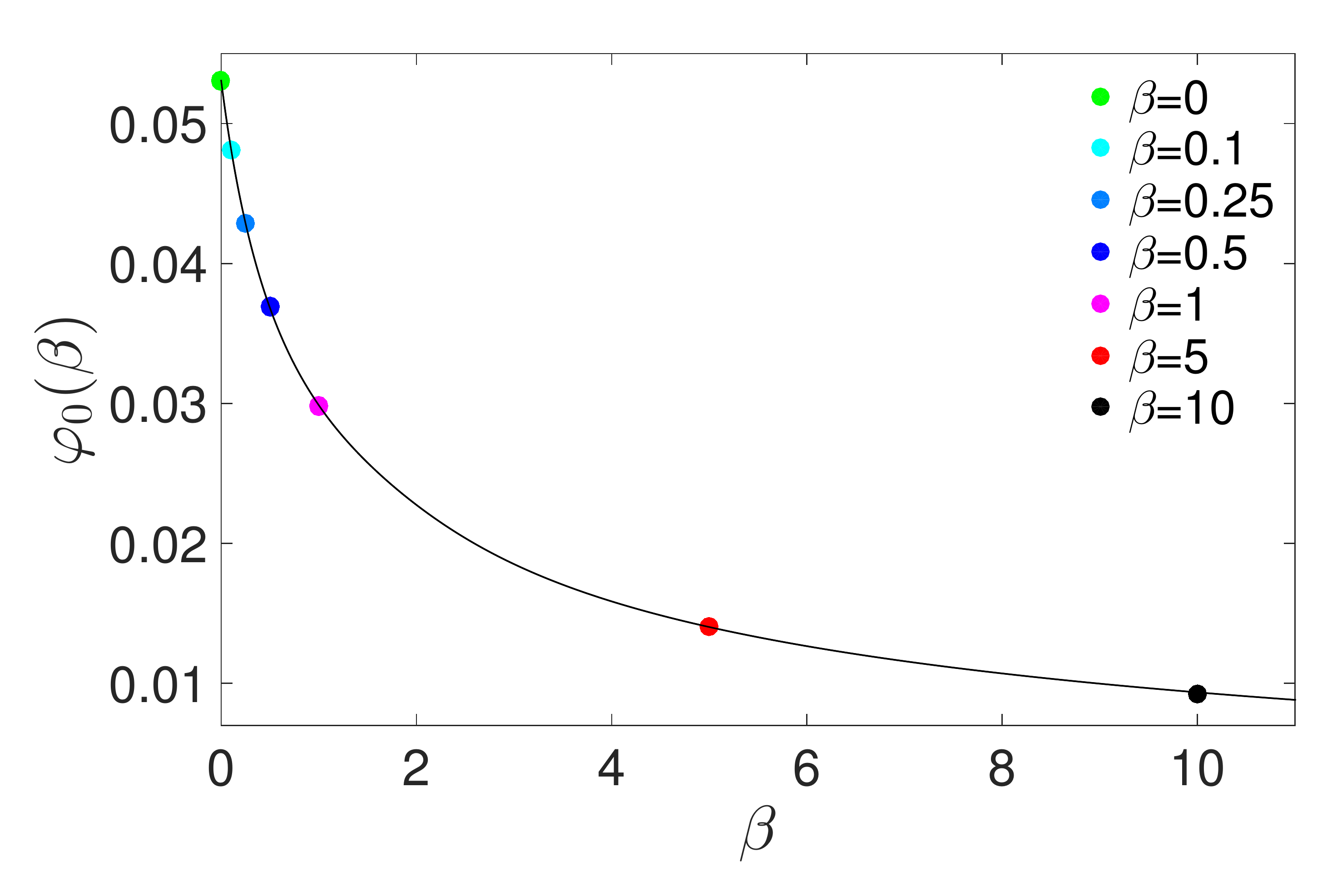}
\caption{\label{fig:f_0_3D}Empirical fit for $\varphi_{0}\left(\beta\right)$
as a function of $\beta$ for different vaues of $\beta$. (Left)
Three dimensions, $256^{3}$ grid. (Right) \label{fig:f_0_2D}\textbf{
}Two dimensions, $512^{2}$ grid.}
\end{figure}

For $r=0$, we use a fitting formula in agreement with (\ref{eq:f_beta_r_0_2D}),
\begin{eqnarray}
\varphi_{0}\left(\beta\right)=\frac{\eta f_{\beta}(0)}{\beta} & = & \frac{z_{0}+z_{1}\beta^{3}\log(\beta)}{1+z_{2}\beta+z_{3}\beta^{2}+z_{4}\beta^{4}},\\
g_{\beta}(0) & = & 0,\nonumber 
\end{eqnarray}
where $z_{1}-z_{4}$ are coefficients (obtained by fitting for each
kernel data over a range of $\beta$'s) and $z_{0}$ is fixed from
the inviscid condition (\ref{eq:f_0_inviscid}).\textbf{ }The empirical
fit for $\varphi_{0}\left(\beta\right)$ is shown in Fig. \ref{fig:f_0_2D}.
Note that for finite $\beta$, one must ensure that the system size
used to tabulate the values of $f_{\beta}$ and $g_{\beta}$ is sufficiently
large, $l\gg h\sqrt{\beta}$.

\begin{figure}[h]
\centering{}\includegraphics[width=0.49\textwidth]{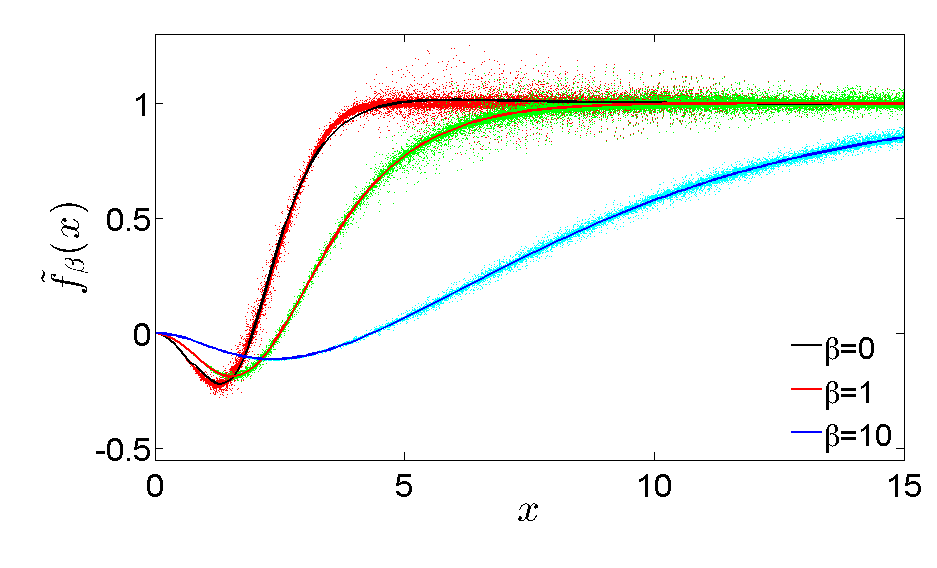}\includegraphics[width=0.49\textwidth]{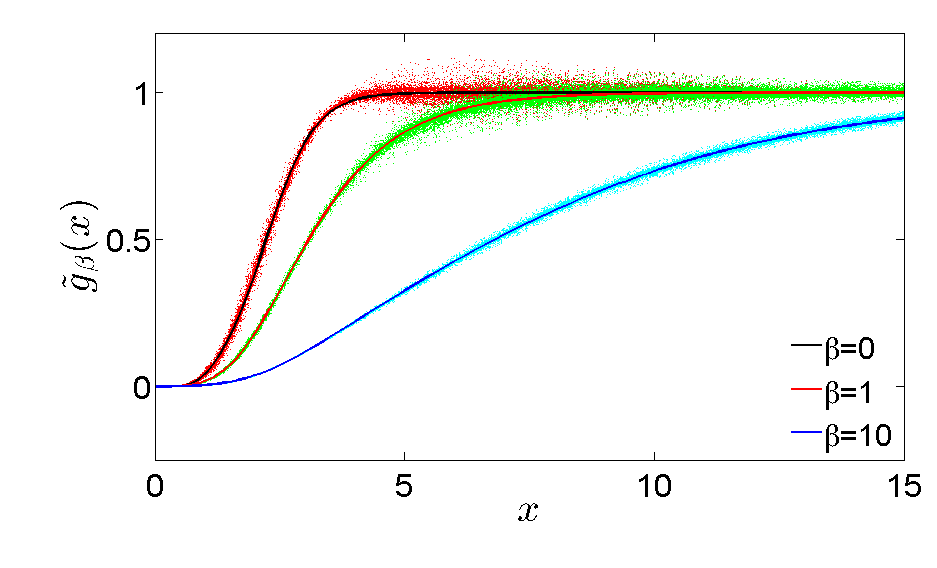}
\caption{\label{fig:fg_6pt_2d}Empirical fitting of $\tilde{f}_{\beta}(r)$
and $\tilde{g}_{\beta}(r)$ in two dimensions for different values
of $\beta$ for the 6-point kernel, $512^{2}$ grid.}
\end{figure}

For $r>0.1h$ we introduce normalized functions $\tilde{f}_{\beta}$
and $\tilde{g}_{\beta}$ via 
\begin{eqnarray}
f_{\beta}(r) & = & -\frac{\beta}{2\pi\eta x^{2}}\cdot\tilde{f}_{\beta}(x),\\
g_{\beta}(r) & = & \frac{\beta}{\pi\eta x^{2}}\cdot\tilde{g}_{\beta}(x),\nonumber 
\end{eqnarray}
where in the inviscid case we take $\beta/\eta=\D t/\rho h^{2}$,
and $x=r/h$ is the normalized distance. For finite $\beta$, we know
that $\tilde{f}_{\beta}(x\gg\sqrt{\beta})\approx\tilde{g}_{\beta}(x\gg\sqrt{\beta})\approx1$
according to (\ref{eq:fg_inv_2D}). The numerical data is fitted with
the empirical fitting functions 
\begin{eqnarray}
\tilde{f}_{\beta}(x) & = & \frac{x^{3}\ln(x)}{\beta(a_{0}+2x)}e^{-\frac{p_{1}x}{\sqrt{\beta}}}+\frac{a_{1}x^{2}+a_{2}x^{3}+a_{3}x^{4}}{1+b_{1}x^{2}+b_{2}x^{3}+a_{3}x^{4}},\\
\tilde{g}_{\beta}(x) & = & \frac{x^{3}}{\beta(c_{0}+4x)}e^{-\frac{p_{2}x}{\sqrt{\beta}}}+\frac{x^{3}}{e^{-p_{3}x}(c_{1}+c_{2}x+c_{3}x^{2})+x^{3}},\nonumber 
\end{eqnarray}
as shown in Fig. \ref{fig:fg_6pt_2d}. Here $a_{0}-a_{3},\ p_{1}-p_{3},\ b_{1}-b_{3},\ c_{0}-c_{3}$
are empirical coefficients, computed by fitting numerical dats for
$\beta$ in $\{0,\ 0.1,\ 0.25,\ 0.5,\ 1.0,\ 5.0,\ 10.0\}$. Intermediate
values in the range $0<\beta<10.0$ are interpolated using linear
interpolation, and larger $\mbox{\ensuremath{\beta}}$'s are handled
using the steady Stokes fit (\ref{eq:fg_steady_2d}).

\section{\label{AppendixConcentric}Stokes flow between two concentric spheres}

Consider steady Stokes flow around a rigid shell or sphere of radius
$a$, placed in a centered position inside another spherical shell
or cavity of radius $b=a/\lambda$. We consider the case when the
outer shell is moving with velocity $V$ and the inner shell is at
rest, and for simplicity set the viscosity to unity, $\eta=1$. Brenner
\citet{BrennerBook} gives the drag force on the inner sphere for
no slip boundary conditions as 
\begin{equation}
F=-6\pi aVK,\label{eq:drag_inner}
\end{equation}
where 
\[
K=\frac{1-\lambda^{5}}{\alpha}\text{ and }\alpha={1-\frac{9}{4}\lambda+\frac{5}{2}\lambda^{3}-\frac{9}{4}\lambda^{5}+\lambda^{6}}.
\]
Let us denote the constants 
\begin{eqnarray*}
A & = & -\frac{15V}{4a^{2}}\cdot\frac{\lambda^{3}-\lambda^{5}}{\alpha},\\
B & = & \frac{3Va}{2}\cdot\frac{1-\lambda^{5}}{\alpha},\\
C & = & \frac{V}{2}\cdot\frac{1+\frac{5}{4}\lambda^{3}-\frac{9}{4}\lambda^{5}}{\alpha},\\
D & = & \frac{Va^{3}}{4}\cdot\frac{1-\lambda^{3}}{\alpha}.
\end{eqnarray*}
The velocity in the region between the two spherical shells can be
obtained from the expressions given by Brenner as
\[
v_{r}=-\cos\theta\left(\frac{A}{5}r^{2}-B\frac{1}{r}+2C+2D\frac{1}{r^{3}}\right),
\]

\[
v_{\theta}=\sin\theta\left(\frac{A}{5}r^{2}-\frac{B}{2}\frac{1}{r}+2C-D\frac{1}{r^{3}}\right),
\]

\[
v_{\phi}=0,
\]
and the pressure is

\[
\p=\p_{\infty}+\mu B\frac{\cos\theta}{r^{2}}-2\mu Ar\cos\theta,
\]
where $\p_{\infty}=0$ since we impose that the pressure have mean
zero to remove the null mode for pressure. In spherical coordinates,
with the symmetry axes aligned with the direction of the flow, the
traction on the surface of the inner sphere, which is the jump in
the stress across the inner shell, is 
\begin{eqnarray*}
\V{\lambda}=\M{\sigma}\cdot\V n & = & \mu\cos\theta\left(2Ar-\frac{B}{r^{2}}\right)\hat{\r}+\mu\sin\theta\left(Ar+\frac{B}{r^{2}}\right)\hat{\V{\theta}},
\end{eqnarray*}
where $\hat{\r}=(\sin\theta\cos\phi,\quad\sin\theta\sin\phi,\quad\cos\theta)$
and $\hat{\V{\theta}}=(\cos\theta\cos\phi,\;\cos\theta\sin\phi,\;-\sin\theta)$.

\section{\label{AppendixPhysicalBCs}Imposing Physical Boundary Conditions}

The local averaging and spreading operators have to be modified near
physical boundaries, specifically, when the support of the kernel
$\delta_{a}$ overlaps with a boundary. A proposal for how to do that
has been developed by Yeo and Maxey \citet{ForceCoupling_Channel},
and an alternative proposal has been developed in the context of the
immersed boundary method by Griffith \emph{et al.} \citet{IBMDelta_Boundary}.
Here we have chosen to use the former approach because of its simplicity
and the fact that it is independent of the kernel, as well as the
fact that it ensures that the interpolated velocity strictly vanishes
at a no-slip boundary; this ensures that the mobility of a marker
is a monotonically decreasing function as it approaches a no-slip
boundary. Since the description in \citet{ForceCoupling_Channel}
is limited to steady Stokes flow and a single no-slip boundary, we
give here an algebraic formulation that extends to a variety of boundary
conditions; this formulation is implemented in the IBAMR library and
used in the examples in this paper in non-periodic domains.

The basic idea in the handling spreading and interpolation near boundaries
is to use the standard IB kernel functions in a domain extended with
sufficiently many ghost cells so that the support of all kernels is
strictly within the extended domain. For interpolation, we first fill
ghost cells and then interpolate as usual using the ghost cell values.
For spreading, we take the adjoint operator, which basically means
that we first spread to the extended domain including ghost cells
in the usual manner, and then we accumulate the value spread to the
ghost cell in the corresponding interior grid point, using the \emph{same}
weight (coefficient) that was used when filling ghost cells for the
purposes of interpolation.

This process requires a consistent method for filling ghost cells,
that is, for extending a (cell-centered or staggered) field $\V u$
from the interior to the extended domain. In general, this will be
an affine linear mapping of the form
\[
\V u_{\text{ext}}=\M E\M u_{\text{int}}+\V c,
\]
where $\M E$ is an extension matrix and $\V c$ encodes inhomogeneous
boundary conditions. Let us denote with $\J_{0}$ the standard IB
interpolation operator that interpolates an \emph{extended} field
at a position inside the interior of the domain. The interpolated
value in the presence of physical boundary conditions is then given
by the affine linear mapping
\[
\J_{BC}\left(\M u_{\text{int}}\right)=\J_{0}\V u_{\text{ext}}=\J_{0}\M E\M u_{\text{int}}+\J_{0}\V c.
\]
The corresponding spreading operation is defined to be the adjoint
of $\J_{BC}$ for homogeneous boundary conditions, as this ensures
energy conservation in the absence of boundary forcing. Specifically,
we use
\[
\S_{BC}=\M E^{T}\J_{0}^{\star}=\M E^{T}\S_{0}.
\]

The specific form of the extension operator $\M E$ used in our implementation
is based on linear extrapolation to a given ghost point based on the
corresponding value in the interior \emph{and} the values at the boundary
as specified in the boundary conditions. Specifically, for homogeneous
Neumann conditions we do a mirror image $u_{\text{ghost}}=-u_{\text{\text{int}}}$,
while for homogeneous Dirichlet boundary conditions, such as no slip
boundaries, we simply do a mirror inversion $u_{\text{ghost}}=-u_{\text{\text{int}}}$.
This makes our implementation exactly identical to that proposed by
Yeo and Maxey \citet{ForceCoupling_Channel} in the context of the
FCM method. One can think of this approach to no-slip boundaries as
taking an inverted mirror image of the portion of the kernel outside
of the domain \citet{ForceCoupling_Channel}. Note that the same $\M E$
is used to implement boundary conditions both in the fluid solver
and when interpolating/spreading near boundaries; this greatly simplifies
the implementation without lowering the second-order accuracy of the
fluid solver \citet{NonProjection_Griffith}. In our implementation,
we use simple transpose for spreading which is the adjoint $\M E^{\star}$
with respect to the standard inner product.


\begin{thebibliography}{10}

\bibitem{DirectForcing_Uhlmann}
M.~Uhlmann.
\newblock {An immersed boundary method with direct forcing for the simulation
  of particulate flows}.
\newblock {\em J. Comp. Phys.}, 209(2):448--476, 2005.

\bibitem{RigidIBAMR}
A.~P.~S. Bhalla, R.~Bale, B.~E. Griffith, and N.~A. Patankar.
\newblock {A unified mathematical framework and an adaptive numerical method
  for fluid-structure interaction with rigid, deforming, and elastic bodies}.
\newblock {\em Journal of Computational Physics}, 250:446--476, 2013.

\bibitem{IBM_Projection}
K.~Taira and T.~Colonius.
\newblock {The immersed boundary method: a projection approach}.
\newblock {\em J. Comp. Phys.}, 225(2):2118--2137, 2007.

\bibitem{RigidIBAMR_ZhangGuy}
Calvin Zhang, Robert~D Guy, Brian Mulloney, Qinghai Zhang, and Timothy~J Lewis.
\newblock Neural mechanism of optimal limb coordination in crustacean swimming.
\newblock {\em Proceedings of the National Academy of Sciences},
  111(38):13840--13845, 2014.

\bibitem{IBM_Rigid_Boundary}
Shen-Wei Su, Ming-Chih Lai, and Chao-An Lin.
\newblock An immersed boundary technique for simulating complex flows with
  rigid boundary.
\newblock {\em Computers \& fluids}, 36(2):313--324, 2007.

\bibitem{IBM_PeskinReview}
C.S. Peskin.
\newblock {The immersed boundary method}.
\newblock {\em Acta Numerica}, 11:479--517, 2002.

\bibitem{TetherPoint_IBM}
Joseph~M Teran and Charles~S Peskin.
\newblock Tether force constraints in stokes flow by the immersed boundary
  method on a periodic domain.
\newblock {\em SIAM Journal on Scientific Computing}, 31(5):3404--3416, 2009.

\bibitem{IBSE_Poisson}
David~B Stein, Robert~D Guy, and Becca Thomases.
\newblock Immersed boundary smooth extension: A high-order method for solving
  pde on arbitrary smooth domains using fourier spectral methods.
\newblock {\em Journal of Computational Physics}, 304:252--274, 2016.

\bibitem{FIISPA_Patankar}
O.~M. Curet, I.~K. AlAli, M.~A. MacIver, and N.~A. Patankar.
\newblock A versatile implicit iterative approach for fully resolved simulation
  of self-propulsion.
\newblock {\em Computer Methods in Applied Mechanics and Engineering},
  199:2417--2424, 2010.

\bibitem{FIISPA_Ardekani}
Arezoo~Motavalizadeh Ardekani, Sadegh Dabiri, and Roger~H Rangel.
\newblock Collision of multi-particle and general shape objects in a viscous
  fluid.
\newblock {\em Journal of Computational Physics}, 227(24):10094--10107, 2008.

\bibitem{CFD_Patankar}
S.V. Patankar.
\newblock {\em Numerical heat transfer and fluid flow}.
\newblock Hemisphere Pub, 1980.

\bibitem{IBM_Sphere}
T.T. Bringley and C.S. Peskin.
\newblock {Validation of a simple method for representing spheres and slender
  bodies in an immersed boundary method for Stokes flow on an unbounded
  domain}.
\newblock {\em J. Comp. Phys.}, 227(11):5397--5425, 2008.

\bibitem{heil2008solvers}
Matthias Heil, Andrew~L Hazel, and Jonathan Boyle.
\newblock Solvers for large-displacement fluid--structure interaction problems:
  segregated versus monolithic approaches.
\newblock {\em Computational Mechanics}, 43(1):91--101, 2008.

\bibitem{FluidStructure_FEM_AMG}
MW~Gee, U~K{\"u}ttler, and WA~Wall.
\newblock Truly monolithic algebraic multigrid for fluid--structure
  interaction.
\newblock {\em International Journal for Numerical Methods in Engineering},
  85(8):987--1016, 2011.

\bibitem{BrownianDynamics_OrderNlogN}
A.~Sierou and J.~F. Brady.
\newblock {Accelerated Stokesian Dynamics simulations}.
\newblock {\em J. Fluid Mech.}, 448:115--146, 2001.

\bibitem{StokesianDynamics_Rigid}
James~W Swan, John~F Brady, Rachel~S Moore, et~al.
\newblock {Modeling hydrodynamic self-propulsion with Stokesian Dynamics. Or
  teaching Stokesian Dynamics to swim}.
\newblock {\em Physics of Fluids}, 23:071901, 2011.

\bibitem{HYDROLIB}
K~Hinsen.
\newblock {HYDROLIB: a library for the evaluation of hydrodynamic interactions
  in colloidal suspensions}.
\newblock {\em Computer physics communications}, 88(2):327--340, 1995.

\bibitem{SphereConglomerate}
B~Cichocki and K~Hinsen.
\newblock Stokes drag on conglomerates of spheres.
\newblock {\em Physics of Fluids}, 7:285, 1995.

\bibitem{HYDROPRO}
A~Ortega, D~Amor{\'o}s, and J~Garc{\'\i}a~de La~Torre.
\newblock Prediction of hydrodynamic and other solution properties of rigid
  proteins from atomic-and residue-level models.
\newblock {\em Biophysical journal}, 101(4):892--898, 2011.
\newblock Code available at
  \url{http://leonardo.inf.um.es/macromol/programs/hydropro/hydropro.htm}.

\bibitem{HYDROPRO_Globular}
Jos{\'e} Garc{\'\i}a de~la Torre, Mar{\'\i}a~L Huertas, and Beatriz Carrasco.
\newblock Calculation of hydrodynamic properties of globular proteins from
  their atomic-level structure.
\newblock {\em Biophysical Journal}, 78(2):719--730, 2000.

\bibitem{RegularizedStokeslets}
Ricardo Cortez, Lisa Fauci, and Alexei Medovikov.
\newblock {The method of regularized Stokeslets in three dimensions: analysis,
  validation, and application to helical swimming}.
\newblock {\em Physics of Fluids}, 17:031504, 2005.

\bibitem{RegularizedStokeslets_2D}
Ricardo Cortez.
\newblock The method of regularized stokeslets.
\newblock {\em SIAM Journal on Scientific Computing}, 23(4):1204--1225, 2001.

\bibitem{RegularizedBrinkmanlet}
Ricardo Cortez, Bree Cummins, Karin Leiderman, and Douglas Varela.
\newblock Computation of three-dimensional brinkman flows using regularized
  methods.
\newblock {\em Journal of Computational Physics}, 229(20):7609--7624, 2010.

\bibitem{BD_IBM_Graham}
Yu~Zhang, Juan~J de~Pablo, and Michael~D Graham.
\newblock An immersed boundary method for brownian dynamics simulation of
  polymers in complex geometries: Application to dna flowing through a nanoslit
  with embedded nanopits.
\newblock {\em The Journal of Chemical Physics}, 136:014901, 2012.

\bibitem{BrownianDynamics_DNA2}
Richard~M Jendrejack, David~C Schwartz, Michael~D Graham, and Juan~J de~Pablo.
\newblock {Effect of confinement on DNA dynamics in microfluidic devices}.
\newblock {\em J. Chem. Phys.}, 119:1165, 2003.

\bibitem{StokesianDynamics_Wall}
James~W. Swan and John~F. Brady.
\newblock {Simulation of hydrodynamically interacting particles near a no-slip
  boundary}.
\newblock {\em Physics of Fluids}, 19(11):113306, 2007.

\bibitem{StokesianDynamics_Slit}
James~W Swan and John~F Brady.
\newblock Particle motion between parallel walls: Hydrodynamics and simulation.
\newblock {\em Physics of Fluids}, 22:103301, 2010.

\bibitem{StokesianDynamics_Confined}
James~W Swan and John~F Brady.
\newblock The hydrodynamics of confined dispersions.
\newblock {\em Journal of Fluid Mechanics}, 687:254, 2011.

\bibitem{BD_LB_Comparison}
Anthony~JC Ladd, Rahul Kekre, and Jason~E Butler.
\newblock {Comparison of the static and dynamic properties of a semiflexible
  polymer using lattice Boltzmann and Brownian-dynamics simulations}.
\newblock {\em Physical Review E}, 80(3):036704, 2009.

\bibitem{RegularizedStokeslets_Walls}
Josephine Ainley, Sandra Durkin, Rafael Embid, Priya Boindala, and Ricardo
  Cortez.
\newblock The method of images for regularized stokeslets.
\newblock {\em Journal of Computational Physics}, 227(9):4600--4616, 2008.

\bibitem{RegularizedStokeslets_Periodic}
Karin Leiderman, Elizabeth~L Bouzarth, Ricardo Cortez, and Anita~T Layton.
\newblock A regularization method for the numerical solution of periodic stokes
  flow.
\newblock {\em Journal of Computational Physics}, 236:187--202, 2013.

\bibitem{NonProjection_Griffith}
B.E. Griffith.
\newblock {An accurate and efficient method for the incompressible
  Navier-Stokes equations using the projection method as a preconditioner}.
\newblock {\em J. Comp. Phys.}, 228(20):7565--7595, 2009.

\bibitem{StokesKrylov}
M.~Cai, A.~J. Nonaka, J.~B. Bell, B.~E. Griffith, and A.~Donev.
\newblock {Efficient Variable-Coefficient Finite-Volume Stokes Solvers}.
\newblock {\em Comm. in Comp. Phys. (CiCP)}, 16(5):1263--1297, 2014.

\bibitem{ForceCoupling_Stokes}
S.~Lomholt and M.R. Maxey.
\newblock {Force-coupling method for particulate two-phase flow: Stokes flow}.
\newblock {\em J. Comp. Phys.}, 184(2):381--405, 2003.

\bibitem{ISIBM}
F.~Balboa Usabiaga, R.~Delgado-Buscalioni, B.~E. Griffith, and A.~Donev.
\newblock {Inertial Coupling Method for particles in an incompressible
  fluctuating fluid}.
\newblock {\em Comput. Methods Appl. Mech. Engrg.}, 269:139--172, 2014.
\newblock Code available at \url{https://code.google.com/p/fluam}.

\bibitem{IrreducibleActiveFlows_PRL}
Somdeb Ghose and R~Adhikari.
\newblock Irreducible representations of oscillatory and swirling flows in
  active soft matter.
\newblock {\em Physical review letters}, 112(11):118102, 2014.

\bibitem{ActiveSuspensions}
Donald~L Koch and Ganesh Subramanian.
\newblock Collective hydrodynamics of swimming microorganisms: Living fluids.
\newblock {\em Annual Review of Fluid Mechanics}, 43:637--659, 2011.

\bibitem{BrownianBlobs}
S.~Delong, F.~Balboa Usabiaga, R.~Delgado-Buscalioni, B.~E. Griffith, and
  A.~Donev.
\newblock {Brownian Dynamics without Green's Functions}.
\newblock {\em J. Chem. Phys.}, 140(13):134110, 2014.
\newblock Software available at
  \url{https://github.com/stochasticHydroTools/FIB}.

\bibitem{ForceCoupling_Fluctuations}
Eric~E. Keaveny.
\newblock Fluctuating force-coupling method for simulations of colloidal
  suspensions.
\newblock {\em J. Comp. Phys.}, 269(0):61 -- 79, 2014.

\bibitem{SELM}
P.~J. Atzberger.
\newblock {Stochastic Eulerian-Lagrangian Methods for Fluid-Structure
  Interactions with Thermal Fluctuations}.
\newblock {\em J. Comp. Phys.}, 230:2821--2837, 2011.

\bibitem{IBAMR_Fish}
Amneet Pal~Singh Bhalla, Rahul Bale, Boyce~E. Griffith, and Neelesh~A.
  Patankar.
\newblock Fully resolved immersed electrohydrodynamics for particle motion,
  electrolocation, and self-propulsion.
\newblock {\em Journal of Computational Physics}, 256:88 -- 108, 2014.

\bibitem{IBAMR}
B.E. Griffith, R.D. Hornung, D.M. McQueen, and C.S. Peskin.
\newblock {An adaptive, formally second order accurate version of the immersed
  boundary method}.
\newblock {\em J. Comput. Phys.}, 223(1):10--49, 2007.
\newblock Software available at \url{http://ibamr.googlecode.com}.

\bibitem{IBAMR_HeartValve}
B.E. Griffith.
\newblock Immersed boundary model of aortic heart valve dynamics with
  physiological driving and loading conditions.
\newblock {\em Int J Numer Meth Biomed Eng}, 28:317--345, 2012.

\bibitem{IBM_Staggered}
B.E. Griffith.
\newblock {On the volume conservation of the immersed boundary method}.
\newblock {\em Commun. Comput. Phys.}, 12:401--432, 2012.

\bibitem{ImprovedLeak}
Charles~S Peskin and Beth~Feller Printz.
\newblock Improved volume conservation in the computation of flows with
  immersed elastic boundaries.
\newblock {\em Journal of computational physics}, 105(1):33--46, 1993.

\bibitem{IBFE}
Boyce~E Griffith and Xiaoyu Luo.
\newblock Hybrid finite difference/finite element version of the immersed
  boundary method.
\newblock {\em Submitted in revised form}, 2012.

\bibitem{DivFreeIB}
Y.~X. Bao, C.~S. Peskin, B.~Griffith, D.~McQueen, and A.~Donev.
\newblock {C. S. Peskin}.
\newblock In preparation, 2016.

\bibitem{New6ptKernel}
Y.~X. Bao, J.~Kaye, and C.~S. Peskin.
\newblock {A Gaussian-Like Immersed Boundary Kernel with Three Continuous
  Derivatives and Improved Translational Invariance}.
\newblock Preprint ArXiv:1505.07529. Software available at
  \url{https://github.com/stochasticHydroTools/IBMethod}, 2015.

\bibitem{Brinman_Original}
HC~Brinkman.
\newblock A calculation of the viscous force exerted by a flowing fluid on a
  dense swarm of particles.
\newblock {\em Applied Scientific Research}, 1(1):27--34, 1949.

\bibitem{ReactiveBlobs}
A.~Pal~Singh Bhalla, B.~E. Griffith, N.~A. Patankar, and A.~Donev.
\newblock {A Minimally-Resolved Immersed Boundary Model for Reaction-Diffusion
  Problems}.
\newblock {\em J. Chem. Phys.}, 139(21):214112, 2013.

\bibitem{FiniteRe_3D_Ladd}
Reghan~J Hill, Donald~L Koch, and Anthony~JC Ladd.
\newblock Moderate-reynolds-number flows in ordered and random arrays of
  spheres.
\newblock {\em Journal of Fluid Mechanics}, 448:243--278, 2001.

\bibitem{FiniteRe_2D_Ladd}
Donald~L Koch and Anthony~JC Ladd.
\newblock Moderate reynolds number flows through periodic and random arrays of
  aligned cylinders.
\newblock {\em Journal of Fluid Mechanics}, 349:31--66, 1997.

\bibitem{VACF_Ladd}
A.J.C. Ladd.
\newblock {Numerical simulations of particulate suspensions via a discretized
  Boltzmann equation. II. Numerical results}.
\newblock {\em Journal of Fluid Mechanics}, 271(1):311--339, 1994.

\bibitem{ApproximateCommutators}
Howard Elman, Victoria~E Howle, John Shadid, Robert Shuttleworth, and Ray
  Tuminaro.
\newblock Block preconditioners based on approximate commutators.
\newblock {\em SIAM Journal on Scientific Computing}, 27(5):1651--1668, 2006.

\bibitem{IBM_Implicit_Fisher}
H.D. Ceniceros and J.E. Fisher.
\newblock {A fast, robust, and non-stiff Immersed Boundary Method}.
\newblock {\em J. Comp. Phys.}, 230(12):5133 -- 5153, 2011.

\bibitem{RPY_Brinkman}
L~Durlofsky and JF~Brady.
\newblock Analysis of the brinkman equation as a model for flow in porous
  media.
\newblock {\em Physics of Fluids}, 30(11):3329--3341, 1987.

\bibitem{RotnePrager}
Jens Rotne and Stephen Prager.
\newblock Variational treatment of hydrodynamic interaction in polymers.
\newblock {\em The Journal of Chemical Physics}, 50:4831, 1969.

\bibitem{RPY_FMM}
Zhi Liang, Zydrunas Gimbutas, Leslie Greengard, Jingfang Huang, and Shidong
  Jiang.
\newblock A fast multipole method for the rotne--prager--yamakawa tensor and
  its applications.
\newblock {\em Journal of Computational Physics}, 234:133--139, 2013.

\bibitem{RPY_Shear_Wall}
Eligiusz Wajnryb, Krzysztof~A Mizerski, Pawel~J Zuk, and Piotr Szymczak.
\newblock Generalization of the rotne--prager--yamakawa mobility and shear
  disturbance tensors.
\newblock {\em Journal of Fluid Mechanics}, 731:R3, 2013.

\bibitem{StaggeredIBM}
Alexandre~M Roma, Charles~S Peskin, and Marsha~J Berger.
\newblock An adaptive version of the immersed boundary method.
\newblock {\em J. Comput. Phys.}, 153(2):509--534, 1999.

\bibitem{PETSc}
Satish Balay, William~D. Gropp, Lois~Curfman McInnes, and Barry~F. Smith.
\newblock Efficient management of parallelism in object oriented numerical
  software libraries.
\newblock In E.~Arge, A.~M. Bruaset, and H.~P. Langtangen, editors, {\em Modern
  Software Tools in Scientific Computing}, pages 163--202. Birkh{\"{a}}user
  Press, 1997.
\newblock Software available at \url{http://www.mcs.anl.gov/petsc}.

\bibitem{Elman_FEM_Book}
Howard~C Elman, David~J Silvester, and Andrew~J Wathen.
\newblock {\em Finite elements and fast iterative solvers: with applications in
  incompressible fluid dynamics}.
\newblock Oxford University Press, 2014.

\bibitem{SmoothingDelta_IBM}
Xiaolei Yang, Xing Zhang, Zhilin Li, and Guo-Wei He.
\newblock A smoothing technique for discrete delta functions with application
  to immersed boundary method in moving boundary simulations.
\newblock {\em Journal of Computational Physics}, 228(20):7821--7836, 2009.

\bibitem{SmallRe_3D_Ladd}
Reghan~J Hill, Donald~L Koch, and Anthony~JC Ladd.
\newblock The first effects of fluid inertia on flows in ordered and random
  arrays of spheres.
\newblock {\em Journal of Fluid Mechanics}, 448:213--241, 2001.

\bibitem{DirectForcing_Balboa}
F.~Balboa Usabiaga, I.~Pagonabarraga, and R.~Delgado-Buscalioni.
\newblock {Inertial coupling for point particle fluctuating hydrodynamics}.
\newblock {\em J. Comp. Phys.}, 235:701--722, 2013.

\bibitem{MultiblobSprings}
Adolfo Vazquez-Quesada, Florencio Balboa~Usabiaga, and Rafael
  Delgado-Buscalioni.
\newblock A multiblob approach to colloidal hydrodynamics with inherent
  lubrication.
\newblock {\em The Journal of Chemical Physics}, 141(20), 2014.

\bibitem{Mobility2D_Hasimoto}
H~Hasimoto.
\newblock On the periodic fundamental solutions of the stokes equations and
  their application to viscous flow past a cubic array of spheres.
\newblock {\em J. Fluid Mech}, 5(02):317--328, 1959.

\bibitem{Tractions_Fauci}
Harvey~AR Williams, Lisa~J Fauci, and Donald~P Gaver~III.
\newblock Evaluation of interfacial fluid dynamical stresses using the immersed
  boundary method.
\newblock {\em Discrete and continuous dynamical systems. Series B}, 11(2):519,
  2009.

\bibitem{BD_LB_Ladd}
Rahul Kekre, Jason~E. Butler, and Anthony J.~C. Ladd.
\newblock {Comparison of lattice-Boltzmann and Brownian-dynamics simulations of
  polymer migration in confined flows}.
\newblock {\em Phys. Rev. E}, 82:011802, 2010.

\bibitem{BoundaryIntegral_Wall}
Oana Marin, Katarina Gustavsson, and Anna-Karin Tornberg.
\newblock A highly accurate boundary treatment for confined stokes flow.
\newblock {\em Computers \& Fluids}, 66:215--230, 2012.

\bibitem{HE_Spheres_TwoWalls}
S~Bhattacharya, J~Blawzdziewicz, and E~Wajnryb.
\newblock Hydrodynamic interactions of spherical particles in suspensions
  confined between two planar walls.
\newblock {\em Journal of Fluid Mechanics}, 541:263--292, 2005.

\bibitem{BoundaryIntegral_Periodic3D}
Ludvig Af~Klinteberg and Anna-Karin Tornberg.
\newblock Fast ewald summation for stokesian particle suspensions.
\newblock {\em International Journal for Numerical Methods in Fluids},
  76(10):669--698, 2014.

\bibitem{EmpiricalDrag_3D}
Sofiane Benyahia, Madhava Syamlal, and Thomas~J O'Brien.
\newblock Extension of hill--koch--ladd drag correlation over all ranges of
  reynolds number and solids volume fraction.
\newblock {\em Powder Technology}, 162(2):166--174, 2006.

\bibitem{BrownianDynamics_OrderN}
J.~P. Hernandez-Ortiz, J.~J. de~Pablo, and M.~D. Graham.
\newblock {Fast Computation of Many-Particle Hydrodynamic and Electrostatic
  Interactions in a Confined Geometry}.
\newblock {\em Phys. Rev. Lett.}, 98(14):140602, 2007.

\bibitem{BrownianDynamics_FMM}
Shidong Jiang, Zhi Liang, and Jingfang Huang.
\newblock A fast algorithm for brownian dynamics simulation with hydrodynamic
  interactions.
\newblock {\em Mathematics of Computation}, 82(283):1631--1645, 2013.

\bibitem{BoundaryIntegral_Pozrikidis}
Constantine Pozrikidis.
\newblock {\em Boundary integral and singularity methods for linearized viscous
  flow}.
\newblock Cambridge University Press, 1992.

\bibitem{FastHierarchicalSolver}
Sivaram Ambikasaran and Eric Darve.
\newblock An $o (n \log n)$ fast direct solver for partial hierarchically
  semi-separable matrices.
\newblock {\em Journal of Scientific Computing}, 57(3):477--501, 2013.

\bibitem{FastSolver_Ho}
Kenneth~L Ho and Leslie Greengard.
\newblock A fast direct solver for structured linear systems by recursive
  skeletonization.
\newblock {\em SIAM Journal on Scientific Computing}, 34(5):A2507--A2532, 2012.

\bibitem{HODLR_BDLR}
Amirhossein Aminfar, Sivaram Ambikasaran, and Eric Darve.
\newblock A fast block low-rank dense solver with applications to
  finite-element matrices.
\newblock {\em Journal of Computational Physics}, 304:170 -- 188, 2016.

\bibitem{iFMM}
Pieter Coulier, Hadi Pouransari, and Eric Darve.
\newblock The inverse fast multipole method: using a fast approximate direct
  solver as a preconditioner for dense linear systems.
\newblock {\em arXiv preprint arXiv:1508.01835}, 2015.

\bibitem{IBM_Implicit_Comparison}
E.P. Newren, A.L. Fogelson, R.D. Guy, and R.M. Kirby.
\newblock A comparison of implicit solvers for the immersed boundary equations.
\newblock {\em Computer Methods in Applied Mechanics and Engineering},
  197(25-28):2290--2304, 2008.

\bibitem{IBMultigrid_Guy}
Robert~D Guy and Bobby Philip.
\newblock A multigrid method for a model of the implicit immersed boundary
  equations.
\newblock {\em Communications in Computational Physics}, 12(2):378, 2012.

\bibitem{ImplicitIB_Projection}
Qinghai Zhang, Robert~D Guy, and Bobby Philip.
\newblock A projection preconditioner for solving the implicit immersed
  boundary equations.
\newblock {\em Numer. Math. Theor. Meth. Appl.(Special issue on Fluid Structure
  Interactions)}, 7(4):473--498, 2014.

\bibitem{RotnePrager_Periodic}
C.~W.~J. Beenakker.
\newblock {Ewald sum of the Rotne-Prager tensor}.
\newblock {\em J. Chem. Phys.}, 85:1581, 1986.

\bibitem{BrennerBook}
John Happel and Howard Brenner.
\newblock {\em Low Reynolds number hydrodynamics: with special applications to
  particulate media}, volume~1.
\newblock Springer Science \& Business Media, 1983.

\bibitem{ForceCoupling_Channel}
Kyongmin Yeo and Martin~R Maxey.
\newblock Dynamics of concentrated suspensions of non-colloidal particles in
  couette flow.
\newblock {\em Journal of Fluid Mechanics}, 649(1):205--231, 2010.

\bibitem{IBMDelta_Boundary}
B.E. Griffith, X.~Luo, D.M. McQueen, and C.S. Peskin.
\newblock Simulating the fluid dynamics of natural and prosthetic heart valves
  using the immersed boundary method.
\newblock {\em International Journal of Applied Mechanics}, 1(01):137--177,
  2009.

\end{thebibliography}

\end{document}